\documentclass[a4paper,10pt]{article}
\usepackage{amsfonts}
\usepackage{amssymb}
\usepackage[utf8]{inputenc}

\usepackage[english ]{babel}

\usepackage{bbding}

\usepackage[all]{xy}

\usepackage{tikz-cd} 
\usepackage[all]{xy}
\usepackage{amsfonts}
\usepackage{amssymb}
\usepackage{graphicx}
\usepackage{mathtools}
\usepackage{skak} 
\usepackage{foekfont}

\usepackage{calligra}
\usepackage[T1]{fontenc}

\usepackage{amsmath,mathdots}

\usepackage{stmaryrd}

\usepackage[utf8]{inputenc}

\usepackage{mathtools}

\usepackage{tikz}

\definecolor{amber(sae/ece)}{rgb}{1.0, 0.49, 0.0}
\newfont{\rsfsten}{rsfs10 scaled 1200}
\usepackage{MnSymbol}
\usepackage{tikz}
\usepackage{amscd}
\usepackage{MnSymbol}
\usepackage{wasysym}

\usepackage{mathrsfs}

\usepackage{pdfcolmk}

\usepackage{graphicx}

\newcommand*{\rom}[1]{\expandafter\@slowromancap\romannumeral #1@}

\usepackage{wrapfig}

\newcommand{\tightunderset}[2]{%
  \mathop{#2}\limits_{\vbox to .3ex{\kern-0.95ex\hbox{$#1$}\vss}}}
\newcommand{\tightoverset}[2]
{%
  \mathop{#2}\limits_{\vbox to .3ex{\kern-0.95ex\hbox{$#1$}\vss}}}

\newcommand{\oset}[2]{%
  {\mathop{#2}\limits^{\vbox to -.5\ex@{\kern-\tw@\ex@
   \hbox{\scriptsize #1}\vss}}}}
\usepackage[utf8]{inputenc}

\usepackage{stmaryrd}

\usepackage {ifsym}

\usepackage {url}
\title {Metric Inequalities with    Scalar Curvature. }

\author{Misha Gromov}

\begin{document}
\maketitle
\tableofcontents
\begin{abstract} We establish several inequalities for  manifolds with {\it positive scalar curvature}  and, more generally, for the scalar curvature bounded from below.

In so far   as geometry is concerned  these inequalities appear as  generalisations  of the classical   bounds on the {\sl distances between conjugates points} in surfaces with {\it positive sectional  curvatures}.

The techniques  of our proofs is based on the
  Schoen-Yau  {\sl descent method  via minimal hypersurfaces, }
while  
the overall  logic  of    our arguments  is inspired by and closely related to the {\sl  torus splitting argument} in  Novikov's   proof of the {\sl topological invariance of the rational Pontryagin classes.}

\end{abstract}

\section {Formulation of  the Key Inequalities.}

Our point of departure is the following  inequality for {\it torical bands} which are smooth  manifolds  homeomorphic to  {\it tori times intervals}. \vspace{1mm}

{\Large [$\varocircle$$_\pm$]} {\sf \textbf {Torical  $\frac {2\pi}{n}$-Inequality.}} Let $V$ be an  $n$-dimensional torical band, 
$V=\mathbb T^{n-1}\times[-1,+1]$,  where the  boundary is  
$$\partial(V)= \partial_- \cup \partial _+= \partial_-(V) \cup \partial _+(V)= (\mathbb  T^{n-1}\times \{-1\}) \cup(\mathbb T^{n-1}\times\{+1\}).$$

{\it Let $g$ be a smooth Riemannian metric on $V$, where  the  scalar curvature is bounded from below by a positive constant  $\sigma>0$,
$$Sc(g)\geq \sigma>0.$$
Then the distance between the two boundary components of $V$ satisfies 
$$dist_\pm=dist_g(\partial_-(V), \partial _+(V))\leq 2\pi\sqrt{\frac { n-1}{\sigma n}} \hspace {1mm}\left (<\frac {2\pi}{\sqrt\sigma}\right ).
 \leqno{  \mbox {
{\Large $[$}{$\varocircle$}{\small$ _\pm \leq  2\pi\sqrt\frac { n-1}{\sigma n}$}{\Large$]$}}      }$$}

{\sf  On Normalisation of $Sc$.} We use the customary  normalisation of the   scalar curvature, where    the unit spheres satisfy 
$$Sc(S^n)=n(n-1).$$
Thus, by scaling,   the  inequality  $\left[dist_\pm \leq  2\pi\sqrt\frac { n-1}{\sigma n}\right]$ for a non specified  $\sigma>0$  reduces to that for  $Sc(V) \geq n(n-1)$, where it reads 
$$dist_\pm=dist_g(\partial_-(V), \partial _+(V))\leq  \frac { 2\pi}{ n} 
\leqno{  \mbox {
{\Large $[$}{$\varocircle$}{\small$ _\pm \leq  \frac { 2\pi}{ n}$}{\Large$]$}}}  $$

In particular,

\vspace {1mm}

{\it  all torical  bands in the unit sphere,  satisfy
 $$dist(\partial_ -,\partial +)\leq \frac {2\pi}{n}.$$}
 
 This  is obvious for $n=2$,   where 
 $ [$$\varocircle$$_\pm\leq\frac {2\pi}{2}]$ is sharp as well as obvious.  One also    expects  a two line  proof 
of a stronger inequality  for all $n$,   but  to my surprise, I was    unable to directly prove even  the corresponding    inequality for {\it principal curvatures}  of  
 $(n-1)$-tori embedded to $S^n$ , where this  inequality is    formulated  below in terms  of {\it 
 focal coradii}   as follows.\vspace {1mm}

{\sf Normal Tubes, Normal  Bands and   $rad^{\bigodot}(Y)$.}  The {\it normal focal radius} of a smooth submanifold $Y$ in a Riemannian manifold $X$, denoted    $rad^{\bigodot}(Y)=rad^{\bigodot}(Y\subset X)$ is the  maximal $r$ such that the normal  exponential  map  
$$\exp:T_\perp(Y) = T(X){|Y}\ominus T(Y) \to X$$ 
is one-to-one\footnote {It would be  more in the spirit of "focal" to  require the normal  exponential map to be {\it locally} one-to-one, but this, probably,   makes no difference in the present context  for $X=S^n$.} on the subset of vectors  $\nu\in T_\perp(Y)$, such that  $||\nu||< r$.

In other words, this is the maximal $r$ such that  {\it the normal $r$-tube } around $Y$, called {\it normal $r$-band} if $codim(Y)=1$, that is   the   open $r$-neighbourhood $U_r(Y)\subset X$  for $r=   rad^{\bigodot}$, normally projects\footnote{This projection sends  each $x\in U_r(Y)$  to the  unique(!) nearest point in $Y$. } to $Y$ and  fibers  $U_r(Y)$ into   $r$-balls of dimension $dim(X)-dim(Y)$.
\vspace {1mm}

{\it Examples.} (a)  The normal  focal radii and the  geodesic   curvatures of sub-spheres  
$$S^m(\rho)=S_s^m(\rho )\subset S^n=S^n(1)\subset \mathbb R^{n+1}, \mbox { } \rho\leq1,$$  
centred at ponts  $s\in S^n$  are
$$r=rad_{S^n}^{\bigodot}(S^m(\rho))= \arcsin \rho \mbox { and }  curv_{S^n}(S^m(\rho))=  \frac {\sqrt{1-\rho^2}}{\rho}=\tan r.$$

(b) The Clifford torus $ \mathbb T^n_{Cl} \subset S^{2n-1}\subset (\mathbb R^2)^n$, that is the product of $n$ circles of radii   $\frac {1}{\sqrt n}$ in the plane, satisfies: 
  $$rad^{\bigodot}(\mathbb T^n_{Cl})=\arcsin \frac {1}{\sqrt n}.$$
{\sf Conjecturally,}  $\mathbb T^n_{Cl}$ has maximal  $rad^{\bigodot}$ among all $n$-tori smoothly embedded to $S^{2n-1}$.

\vspace {1mm}

{\sf \textbf { Normal Radius Inequality for  $\mathbb T^{n-1}\subset S^n$}}.   {\it If a  smooth hypersurface $Y$  in the unit $n$-sphere is homeomorphic to the $(n-1)$-torus, then 
$$ rad^{\bigodot}(Y)\leq \frac  {\pi}{n}. \leqno  { [{\bigodot}\leq \frac  {\pi}{n}].}$$} 

This inequality-- this will become clear later on  --  is  non-sharp. 

{\sf Conjecturally,} the  sharp constant must be asymptotic   for   $n\to \infty$  to 
$$ \frac{const}{n^\alpha }\mbox{   for some $\alpha > 1$ }. $$

{\it  On Sharpness  of {\small {\Large $[$}{$\varocircle$}{\small$ _\pm \leq \frac { 2\pi}{n}$}{\Large$]$}}. }
This  inequality  agrees with the obvious one in the $2$-sphere  (where the conventionally defined scalar curvature equals twice the sectional  curvature) where the widths of 
the bands between concentric circles as well as the  distances between opposite sides of    (all) quadrilaterals     are bounded by   
$\frac {2\pi}{2}=\pi=diam(S^2)$ and where 
 these inequalities become sharp for doubly punctured
 spheres  and for quadrilaterals  which degenerate 
 to {\it geodesic digons} joining opposite points in $S^2$.

And if $n\geq 2$, we shall see in the next section   that the {\it extremal} bands,      where  $dist_\pm= \frac {2\pi}{n}$,  also have {\it constant scalar curvatures} and   their opposite sides    {\it collapse to points}, but they {\it do not have constant sectional curvatures}  for $n>2$  anymore. 

\vspace{1mm}

 {\sf \textbf  {Quadratic   Decay   Theorem.}}  Let $X$ be a complete Riemannian manifold,  and let
 $$\min_{B(R)} Sc(X)$$ denote the minimum of the scalar curvature (function)  of $X$ on the ball  $B(R)=B_{x_0}(R)\subset X$  for some centre point $x_0\in X$. \vspace {1mm} 
 
 {\it If $X$ is  homeomorphic to $\mathbb T^{n-2}\times \mathbb R^2$, then there exists a constant $R_0=R_0(X, x_0)$, such that 
 $$  \min_{B(R)} Sc(X)\leq \frac{4\pi^2}  {(R-R_0)^2}  \mbox  {  for all } R\geq R_0. \leqno {[\asymp \frac {4\pi^2}  {R^2}]}$$}

 {\it Outline of the   Proof.} Let $X_0\subset X$ corresponds to the  torus  $\mathbb T^{n-2}\times  \{0\}\subset \mathbb T^{n-2}\times \mathbb R^2$ under the homeomorphism $\mathbb T^{n-2}\times  \mathbb R^2 \leftrightarrow X$ and let $R_0= diam_X(X_0)$.
 Then the $(R-R_0)$-neighbourhood $U_{ R-R_0}(X_0)\subset X$ is contained in the ball
 $B_{x_0}(R)$  for $x_0\in X_0$.
 
 If  $U_{ R-R_0}(X_0)$ is  homeomorphic to
  $ \mathbb T^{n-1} \times (-1,+1)$, then   $ [\asymp \frac {4\pi^2}  {R^2}]$   follows from 
the torical  $\frac {2\pi}{n}$-inequality     and if  the topology of    $U_{ R-R_0}(X_0)$ is more complicated, then we apply  a generalisation of  the   $\frac {2\pi}{n}$-inequality from section 4. 

\vspace{1mm}
 
 {\sf On Uniformly Positive Scalar Curvature.}  The obvious  corollary to  $[\asymp \frac {4\pi^2}  {R^2}]$ is \vspace {1mm}
 
 \hspace {8mm}{\it non-existence of complete metrics  with $Sc\geq \sigma>0$} on  $\mathbb T^{n-2}\times \mathbb R^2$. \vspace{1mm}
 
Notice that  there are similar  results   for other  manifolds $X$ proven with  {\it Dirac operators twisted with suitable   "almost flat" bundles  over $X$} [GL 1983],   [HaPS 2015].
 
 However,
for all I know, one can't rule out metrics with uniformly positive scalar curvature on    $\mathbb T^{n-2}\times \mathbb R^2$ with the present day Dirac operator methods. \footnote {Bernhard Hanke indicated to me  how one can do this, see section 7.} 
\vspace {1mm}

 {\sf  Examples of Metrics  on $\mathbb T^{n-2}\times \mathbb R^2$ with Quadratic Decay  of  Scalar Curvature.} Let $g =dt^2 +\varphi(t)^2 d\theta^2$, $t\in [0,\infty) , \theta\in [0, 2\pi] $, be   a   radial (rotationally  symmetric) metrics on $\mathbb R^2$. Then 
 $$Sc(g)(t)=-\frac {2\varphi''(t)}{\varphi(t)};$$
thus,  the metrics 
$$g_{fl}+dt^2 + t^{2\alpha} d\theta^2 $$   
on $\mathbb R^2\times \mathbb T^{n-2}$, where $g_{fl}$ are  flat on  $\mathbb T^{n-2}$  and where   $0<\alpha<1$,   do the job.

\subsection{On Generalisations and Proofs:  Admissions and and Acknowledgements}  

Simple generalisations of everything we stated so far is proven in section 2. 

Then, in the following sections, we formulate  and prove further  {\it generalisations and refinements} of these.
Also we   indicate   additional  {\it applications}  and articulate  several  {\it conjectures.} 
 
  Our approach is based  on   the Schoen-Yau   {\it    dimension descent argument}   [SY-Str 1979],  [SY 2017] accompanied by 
    {\it  torical symmetrization}   [GL 1983] and/or  {\it symmetrizartion by reflection}   [Gr 2014*].

However, some of our arguments have certain limitations  which are indicated below.

\textbf {1.} {\sf Problem with Singularities.}  Applications of minimal hypersurfaces $Y\subset X$ to  $Sc\geq 0$    depends on the regularity of these $Y$  which is  known to hold for  all   $Y$ if $n=dim(X)\leq 7$ and for generic ones for $n=8$ by a Nathan Smale theorem [NS 1993] More recently, Lokhamp [Loh 2016] and Schoen  and Yau [SY 2017] suggested  ways of  bypassing the singularity problem. 

As far as I understand, the  regularity  results  by Schoen and Yau in [SY 2017], such as   theorem 4.6,  suffice for  the needs  of the present paper and this is, probably, true about the corresponding results by Lokhamp. But since I have not studied   these papers in depth, I can  vouch  for the  validity of our proofs only for $n\leq 8$, where the singularity problem does not exist. \vspace {2mm}

\textbf {2.} {\sf Doubling,  Reflecting and Smoothing.} Some  results concerning closed Riemannian manifolds $X$ with  $Sc\geq 0$  generalise to manifolds $X$ with boundary
and, accordingly, to bands of the form $X\times [-1,1]$.

For instance if the boundary $\partial=\partial X$ is  {\it mean curvature convex}, i.e. if $mn.curv(\partial)>0$, then $X$ admits a full fledged theory of minimal hypersurfaces $Y\subset X$ with (free) boundary $\partial Y\subset \partial X$ and the Schoen-Yau descent method applies.

Alternatively, one may take the {\it double}  
 $$\tilde X=X\cup_\partial X$$
and show (see [GL 1980], and section 11.4) that \vspace {1mm}

 {\sl the natural {\it continuous} (but not necessarily smooth, not even $C^1$) 
Riemannian metric $\tilde g=g\&g$  on  $\tilde X$ can be $C^0$-approximated by smooth $C^\infty$-smooth metrics $\tilde g_\varepsilon$ with no decrease of the scalar curvature of $\tilde g$.} (The scalar curvature of $\tilde g$ is unambiguously defined away from $\partial$, which naturally embeds as a hypersurface in $\tilde X$.)

\vspace {1mm}

Thus 

{\sl all results for closed manifolds $X$, including those obtained with the Dirac 

operator methods,  extends to manifolds with mean curvature convex
 
 boundaries.}

(Analysis of {\it extremal cases}, e.g. showing that every  metric $g$  on $X=\mathbb T^{n-1}\times [0,1]$,  where $mn.curv_g(\partial (X))=0$, is {\it Riemannian flat}, needs an additional care.)

\vspace {1mm}

Next, let  $X$, be a manifold with corners, such, for instance, as the Cartesian  product of several manifolds with boundaries, such as  
the $n$-cube 
$$\square^n=\underset {n}{\underbrace {[0,1],\times ...\times [0,1]}},$$
or an $V$ diffeomorphic to $\square^n$.

 Schoen-Yau descent, probably,   applies to such  $X$,  where  one may  encounter  complications with  
   singularities of minimal hypersurfaces, $Y\subset V$ with $\partial Y\subset \partial X$,  at the  corners of $X$ on the boundary  $ \partial X$.  \footnote{Possibly, an adequate  regularity of minimising  hypersurfaces at the corner points  is  known to some people, but I could not locate such result   in the literature.} 

There is an  alternative approach applicable to those   $X$, which  serve as fundamental   
domains   of {\it refection groups} $\Gamma$ acting on manifolds $\tilde X$ without boundaries.

For instance, if $X$ is {\it diffeomorphic} to the $n$-cube $[0,1]^n$, then $\tilde X$ is {\it naturally  homeomorphic} to 
 the Euclidean space $ \mathbb R^n$ acted on by the 
group $\Gamma$  (isomorphic to $\mathbb Z^n \rtimes \mathbb Z^n_2$) generated by reflections of $\mathbb R^n$ in the hyperplanes $\{x_i=0\}$ and $\{x_i=1\}$, $i=1,...,n$, in $\mathbb R^n$.

This $\tilde X$ carries a unique  path metric  $\tilde g$, which  is equal to  the original Riemannian metric  $g$ on $X\subset \tilde X$, where  $ X$ is embedded to $\tilde X$ as a fundamental domain. 

In general, unlike the case where $X$ is isometric to $[0,1]^n$, rather than only diffeomorphic to it, the metric $\tilde g$ is  non-smooth at the boundary  $\partial X\subset \tilde X$
 (as well as at the boundaries of the $\Gamma$-translates of $X$).  In fact, this metric is continuous (but not, in general, $C^1$-smooth) at the smooth points of $\partial X$ which correspond to the interior points in the $(n-1)$-faces of $[0,1]^n$,  and it is only piecewise continuous     at the edges. 
 
 In fact the metric $\tilde g$ is Riemannian continuous if and only  the dihedral angles between the  $(n-1)$-faces of $X$ along the "edges" (corresponding to the  $(n-2)$-faces of $[0,1]^n$  equal  $\frac {\pi}{2}$ and  then $\tilde g$ is $C^1$ (only) if   the $(n-1)$-faces are totally geodesic with respect to this metric.

 \vspace {1mm}

It  follows from  the  Approximation/Reflection Lemma in section 4.9 in [Gr 2014*] that  \vspace {1mm}

{\sf if the the faces of $X$ are mean curvature convex and  if the dihedral angles between these faces  along the corners in $X$ are bounded by $\frac {\pi}{2}$, then $\tilde X$ carries  smooth(!) $\Gamma$-invariant metrics $\tilde g_\varepsilon$, $\varepsilon\to 0$,  which are, in some week sense converge  to the original metric $g$ on   $X \subset \tilde X$    and such  the scalar curvatures of  $\tilde g_\varepsilon$    are bounded from below by $Sc(g)$ on  $X\subset \tilde X$. }\vspace {1mm}

Granted this, many (all?) properties of closed (and complete)  manifolds with $Sc\geq \sigma$ generalise to manifolds  with suitable corners, e.g.

{\sl   non-existence of  cubical (i.e. diffeomorphic to $[0,1]^n$) domains in Riemannian manifolds  $X$ with $Sc(X)\geq 0 $ (e.g. in $X=\mathbb R^n$), 
 such that the $(n-1)$-faces   of these domains have $mn.curv>0$
and the dihedral angles between  these faces 
   are $\leq \frac {\pi}{2}$,  see [Gr 2014*]) . }

\vspace {1mm}

However -- this was pointed out by the referee of the present paper  --  the way it was written in the original  version of this paper was  unsatisfactory.

We discuss  in  section  11.8 what should be done about it. 
\vspace {1mm}

\textbf {3} {\sf Acknowledgements.} I want to express my gratitude  to the anonymous referee  for the time and  effort he/she has spent  reading this article. Besides pointing out   inadequacy  of the  first version of our treatment of smoothing multiple corners and of minimal hypersurfaces in non-compact manifolds, the referee has indicated  many other inconsistencies and    errors  in  the paper -- more than 60 in the first draft of the paper and more than  100 in the present version.  Also he/she has made    several  insightful  suggestions  which necessitated the appearance of  the    appendix (section 11) in this  article.\vspace {1mm}

\section {Bounds on Widths of Over-torical and  Related Riemannian  Bands.} 

A  {\it  band } is a manifold $V$  with two distinguished disjoint non-empty  subsets in the boundary $\partial(V)$, denoted 
 $$\partial _-= \partial_-(V)\subset \partial V \mbox { and } \partial _+= \partial_+(V)\subset \partial V.$$

A band is called {\it proper} if $\partial_\pm$ are unions of connected components of $\partial  V$  and 
$$\partial _-\cup \partial _+=\partial V.$$

 {\it Band maps}  $V\to \underline V$ are those  continuous ones  which respect these $\pm$-boundaries,  $\partial _\pm\to \underline \partial _\pm$.
 
 If $V$ is endowed with a Riemannian metric then the {\it width of a band}  is   the distance between   $\partial _-$ and $\partial _+$, that is the infimum of length of curves in $V$  between  $\partial _-$ and $\partial _+$.

 A compact proper  orientable  band is called {\it over-torical}  if it admits  a  band map to the toric band,
 $$f:V\to \underline V= \mathbb T^{n-1}\times [-1,+1], \hspace {1mm} n=dim (V),$$
  with {\it non-zero degree}.

Another way to put it is by saying that  the  relative fundamental class $[V]\in H^n(V, \partial V;\mathbb Q)$ decomposes to the product 
$$[V]=h_1\smile...\smile h_{n-1}\smile h_n$$
where $h_i$, $ i=1,...,n-1$, are (absolute) 1-dimensional cohomology classes,  $h_i\in H^1(V;\mathbb Q)$,  
and $h_n\subset H^1(V, \partial V;\mathbb Q)$   is  the   (relative) class,   of the   differential of   a function $V\to [-1,+1]$ such that $\partial_\pm \mapsto 
\pm 1$.

If $V$ is non-orientable, then {\sl overtorical} means that an orientable finite cover of $V$ is overtorical.\vspace {1mm}

{\sf \textbf  {Torical Symmetrization. }} {\it There exists a quasi-functorial symmetrizartion "operator" from Riemannian over-torical bands to torical ones 
$$ {\sf Sym} : V\leadsto \underline V$$
where $\underline V$ admits a free isometric action of the torus  $\mathbb T^{n-1}$  and such that 
$$width (\underline V)\geq width( V)$$
and  
$$Sc(V)>\sigma  \Rightarrow Sc(\underline V)>\sigma.$$} 

{\it Proof.}  This is  proven  in a  slightly different form in  [GL 1983] for $n\leq7$ by induction  as it is   explained below.  

(Earlier,   such  symmetrization for $n=3$ was used by  Fisher-Colbrie and Schoen  [FCS 1980], while     the proof for $n=8$ is essentially the same  as for $n\leq 7$ due to  Nathan's Smale generic regularity theorem.)

 \vspace {1mm}

{\sf Induction Step.} Let $V_k$ be a $\mathbb T^k$ invariant Riemannian band, $k=0,...,n-2$, which admits a 
 $\mathbb T^k$-equivariant band map to the torical band
 $$f_k: V_k\to \mathbb T^{n-1}\times [-1,+1]$$ 
 where  $\mathbb T^k$ acts on $\mathbb T^{n-1}\times [-1,+1]$ via the standard   (coordinate) embedding $\mathbb T^k\subset  T^{n-1}$  and such that $deg(f_k)\neq 0$.

 Let $Y_k\subset V_k$ be a volume  minimising hypersurface which is homologous  to the $f_k$-pullback of 
 $$\mathbb T^{n-2}\times [-1,+1]\subset  \mathbb T^{n-1}\times [-1,+1]$$
for the torus $\mathbb T^{k+1}\supset \mathbb T^{k}$, where  "homologous" refers to the relative  group   $H_{n-1}(V_k; \partial V_k= \partial_-\cup \partial _+)$.

It is easy to see that this   $Y_k$ is $\mathbb T^k$-invariant and that the lowest eigenfunction $\phi(y)$ of the second variation operator $L$ on $Y_k$,
$$ L =-\Delta+\frac{1}{2}(Sc(Y_k) -Sc(V_k|Y_k) -||curv_{V_k}(Y_k)||^2),$$
 is also  $\mathbb T^k$-invariant. (Here,  $\Delta$ is the Laplacian on $Y_k$, that is $\sum_i \frac  {\partial ^2}{dy_i^2}$ and $curv_X(Y)$ denotes the second fundamental  form of $Y_k\subset V_k$.) 

Then we let $V_{k+1}=Y_k\times \mathbb T^1$  with the metric $dy^2+\phi^2dt^2$, where a simple computation shows that if the scalar curvature of $V_k$ restricted to $Y_k$ is $\geq \sigma$, then the scalar curvature of $V_{k+1}$ is also bounded from below by $\sigma$.

It is also clear that  $V_{k+1}$ admits a $\mathbb T^{k+1}$-equivariant band  map of degree $=deg(f_k)$ to the torical band and that  $width(V_{k+1})\geq width(V_k)$. 

Thus, the inductive step is completed and the existence of torical symmetrization follows.
(See [GL 1983] for details).
\vspace {1mm}

{\it Remark on Singularities.}  $\mathbb T^{k}$-invariant  minimal hypersurfaces in $V_{k}$ correspond to 
hypersurfaces in  the quotient  manifolds $V_k/\mathbb T^k$, which are minimal with respect to the quotient metrics with  obvious conformal weights.  Then theorem 4.6 in [SY 2017] says, in effect, that even if some  hypersurfaces $Y_k$ were singular, say for $n-k\geq 8$, the final $\mathbb T^{n-1}$-symmetric  $V_{n-1}= \underline  V$ are non-singular. 

 (Schoen and Yau formulate their theorem for closed manifolds  but the needed  regularity for manifolds $V$ with boundaries  trivially  reduces  to that for doubles of $V$.)

\vspace {1mm} 

{\sf \textbf  { $\frac {2\pi}{n}$-Inequality  for Over-Torical Bands.}} {\it Overtorical bands with scalar curvatures $\geq n(n-1) (=Sc(S^n))$ satisfy
$$width(V)\leq \frac { 2\pi }{ n}. \leqno{\mbox { $\left [dist_ \pm \leq \frac {2\pi}{n}\right ]$}}$$}

{\it Proof.} Torical symmetrization reduces the general  case to that of $\mathbb T^{n-1}$-invariant metrics $g$, on torical  bands,  where 
$$g=dt^2 +\sum_i\varphi_i(t)^2d\tau_i^2,  \hspace{1mm} i=2,3,...,   n.$$
Then 
one easily  computes  
$$Sc(g)(t, \tau_2,...,\tau_n)= -2\sum_i \frac {\varphi''(t)}{\varphi_i(t)}-2\sum_{i<j}\frac {\varphi_i'(t) }{\varphi_i (t)}\frac {\varphi_j'(t) }{\varphi_j (t)}$$
and shows that the the longest $t$-interval where this function remains  defined for  $Sc(g)\geq \sigma>0$  is achieved with $\varphi_2=...=\varphi_n=\varphi$, where the proof follows by simple computation  on  p.401 in [GL 1983] which is
 reproduced below in the description of optimal (maximal)  torical bands with $Sc \geq \sigma$.

\vspace {1mm}

{\it Proof of  Propositions  from Section 1.}  The  inequality $\left [dist_ \pm \leq \frac {2\pi}{n}\right ]$
 implies everything we have stated so far, where in the case of the  quadratic decay theorem one needs to observe that the domains  $U_{ R-R_0}(X_0)\subset X$  (defined following the statement of this theorem) are,   in an obvious sense, {\it open  overtorical bands}  to which  the above   $\frac {2\pi}{n}$-Inequality applies.

Notice at this point that this argument automatically delivers the following\vspace {1mm}

{\sf Generalisation of The Quadratic Decay Theorem.} {\sl If a complete orientable  Riemannian $n$-manifold $X$  admits a proper continuous map $X\to \mathbb T^{n-2}\times \mathbb R^2$  of non-zero degree,  then the minima of the scalar curvature of $X$ over concentric $R$-balls in $X$ satisfy
 $$  \min_{B( R)} Sc(X)\leq \frac{4\pi^2}  {(R-R_0)^2}  \mbox  {  for  some $R_0\geq 0$ and  all } R\geq R_0. \leqno {[\asymp \frac {4\pi^2}  {R^2}]^\ast}$$}

 \vspace {1mm}
  
{\sf \large  Optimality of $\frac {2\pi}{n}$}.  {\it Every smooth  manifold $V=Y\times[-1,1]$ admits a Riemannian metric $g=g_\varepsilon $ with $Sc(g)\geq n(n-1)$ and the $g$-distance 
between the two boundary components  $Y\times \{-1\}$  and $Y\times \{1\}$  in $V$   equal $2\pi/n-\varepsilon$  for a given $\varepsilon>0$.}
 
 \vspace {1mm}

 For instance, if $Y=S^1$, then the {\it spherical suspension} 
 $V=${\large $\Theta$}$(Y)$ serves this purpose  for all $\varepsilon>0$.
 
 More generally, given  a Riemannian  metric $g_0$ on $Y$ and a  real function $\varphi(t)$,  let $ g=dt^2+\varphi(t)^2g_0$ be the metric on $Y\times [-l,l]$. If $g_0$ is flat then
 $$\sigma=Sc(g)=-2(n-1)\frac {\varphi''}{\varphi}-(n-1)(n-2)\frac {\varphi'^2}{\varphi^2},\leqno {\RIGHTcircle}$$
or 
$$\frac{\sigma}{n-1}= -2(\frac {\varphi''}{\varphi}-\frac {\varphi'^2}{\varphi^2}) -n\frac {\varphi'^2}{\varphi^2}, $$
 that is $$-2f' -nf^2=\frac{\sigma}{n-1}\mbox {   for } f=\frac {\varphi'}{\varphi}.$$ 
 Now let $\sigma=Sc(S^n)=n(n-1) $ and rewrite the above as 
 $$\frac {f'}{1+f^2} =(\arctan f)' =\frac {n}{2}$$
 and 
 $$f=f(t)=\tan \frac {n}{2}t$$
 which  is  a function defined on the $\frac{2\pi}{n}$-interval $(- \frac{ \pi}{n}, + \frac{ \pi}{n})$.

This settles the matter for flat manifolds  $Y$ and the general case follows by rescaling general  metrics in $Y$   with a large constants.

\section{Toric Bands  in Spheres and Lower Bounds on  Lipschitz  Constants of  Map $X\to S^n$ in terms of  $Sc(X)$.}

Suppose, there is  a {\it toric  band of width $d$}  in the  unit  $n$-sphere $S^n$ that  is a domain  $\underline  V\subset S^n$   which is homeomorphic to $\mathbb T^{n-1}\times [-1,1]$  and such that  the distance between the two boundary components $\partial_\pm(\underline V)$  of $\underline V$ is equal to $d$  and let 
$f$ be  a continuous   map of non-zero  degree 
 from an oriented  Riemannian $n$-manifold $X$ to $S^n$.

Recall  that  saying "degree" presupposes that $f$ is locally constant at infinity, i.e. constant on each boundary component of $X$  and, if $X$ is non-compact,  on every component of the complement to some (large) compact subset in $X$, 
and let us   additionally assume  that the (finite)  $f$-image   of the so defined infinity    {\it does not intersect $V$}. (This is relevant only if   $\partial X$ is disconnected and/or if $X$ is disconnected at infinity.)  

Then the pullback $V=f^{-1} (\underline V)\subset X$  is 
a Riemannian over-torical  band, such that  the distance between the  two parts $\partial_\pm(V)$  of its boundary is $\geq \lambda^{-1}d$, and 
 the inequality  
 $$d=width(V)\leq 2\pi\sqrt{\frac { n-1}{\sigma n}}$$
 (this is $ [dist_\pm \leq \frac {2\pi d}{n}]$ formulated   for $\sigma=n(n-1)$  in  the previous section) shows that
 $$ [Sc(X)\geq \sigma] \mbox  { \large $\Rightarrow$}\left [ Lip(f)\geq \frac {d}{2\pi}\sqrt \frac {\sigma n}{n-1}\right ], \leqno \bigcirc_{Lip}^n$$
 where, recall,
 $$Lip(f) =\sup_{x_1\neq x_2}  \frac {dist_{S^n}(f(x_1), f(x_2))}{dist_{X}(x_1,x_2)}.$$


Notice that the $\frac {2\pi}{n}$-inequality  $[dist_\pm\leq \frac {2\pi}{n}]$, that is  (essentially) 
 $\bigcirc_{Lip}^n $  applied to the identity map,  shows that torical bands in $S^n$ have widths $d\leq \frac {2\pi}{n}$.

{\sf Conjecturally}, the maximal widths  $d$ for large $n\to \infty$ must be asymptotic to 
 $\frac {1}{n^{1+\alpha}}$  { for some $ \alpha>0.$

\vspace {1mm}

{\it  Round Tori  and in  $S^n$ and  in $\mathbb R^n$}.  
Let us show that this $\alpha$  must be $\leq \frac {1}{2}$ by exhibiting 
 embedded tori $\mathbb T^{n-1}$ with bands of width $\approx \frac {1}{ n^{\frac {3}{2}}}$
around them, where we use the following terminology.\vspace {1mm}

{\sf \textbf  {Over-Torical Width  $width_{\hat{\cal T}}(X)$.}}   This is defined for Riemannian manifolds $X$ as the supremum of numbers $d$, such that $X$ admits  an equidimensional   locally isometric  (not necessarily globally one-to-one) immersion from an overtorical Riemannian band of width $d$.

For instance, it is obvious that 

$$ width_{\hat{\cal T}}(S^2)=\pi.$$

More significantly,  since the     Clifford torus in $S^3$ has $rad^{\bigodot} =\pi/4$, (see section 1)
$$width_{\hat{\cal T}}(S^3)\geq \pi/2$$
and consequently, \vspace {1mm}

{\it all continuous  maps $f$  from a Riemannian (possibly incomplete)    $3$-manifold $X$ with  $ Sc(X) \geq 6=Sc(S^3)$  to $S^3$ 
 which are  constant at infinity and have $deg(f)\neq 0$, satisfy
 $$Lip(f)\geq \frac {3}{4}.$$}

This improves the inequality    $Lip(f) \geq  \frac {3}{8\pi} $ from [GL 1983] but falls short of the {\sf conjectural bound} $Lip(f)\geq 1$.

 Another natural {\sf conjecture}  is the  {\it equality}  
 $$width_{\hat{\cal T}}(S^3)=\pi/2.$$ Moreover, one expects that 
\vspace {1mm} 

\textbf {\sl all  (possibly incomplete) 3-manifolds $X$ with sectional  curvatures $\geq 1$  
satisfy
$$width_{\hat {\cal T}}(X)\leq \pi/2.$$}

\vspace {1mm} 

Starting from $n=4$, codimension one tori in $S^n$   can't be  rotationally invariant any more; we construct certain   "roundish" ones   with relatively large  focal coradii    $r=rad^{\bigodot}$, i.e. with  the normal  exponential maps of these tori is one-to-one within distance  $\leq r$ from them.

We  construct these tori in the unit  Euclidean $n$-balls  (rather than in the unit  spheres) by induction as follows. 

 Given  codimension one tori $Y_1\subset  B^{n_1}\subset  \mathbb R^{n_1}$, and  $Y_2\subset  B^{n_2}\subset \mathbb R^{n_2}$ with focal coradii $r_1$ and $r_2$, take  $c_1, c_2>0$,
such that 
$$ c_1^2+c_2^2=1\mbox {  and } c_1r_1= c_2r_2,$$
 observe that the product of the $c_i$-scaled $Y_i$ in  $\mathbb R^{n_i}$ is contained in the unit ball
$$Y_\times =c_1Y_1\times  c_2Y_2\subset B^{n_1+n_2}\subset \mathbb R^{n_1+n_2}=\mathbb R^{n_1}\times \mathbb R^{n_2}. $$
and 
$$rad^{\bigodot}(Y_\times)= r_\times=c_1r_1=c_2r_2.$$

Then let $$Y_{\times_{+}}=\frac{1}{1+\delta}(Y_\times )_{+\delta}\subset B^{n_1+n_2}$$
be the $\frac{1}{1+\delta}$-scaled boundary of the $\delta$-neighbourhood of $Y_\times$ in  $\mathbb R^{n_1+n_2}$    with
$\delta=\frac {1}{2} r_\times $  and observe that

  $$rad^{\bigodot} (Y_{\times_+})=rad^{\bigodot} (Y_{\times})\cdot \frac{1}{2}\cdot\left ( \frac{1}{{1+\frac{1}{2}r_\times}}\right)=\frac {r_\times}{2+r_\times}. $$
 

In particular,  if a torus  $Y=Y_1=Y_2=Y(n)\subset B^n$ has normal focal radius $r=r(n)$,  the resulting  
$Y(2n)= (Y\times Y)_+\subset B^{2n}$ satisfies
 $$r(2n)=  rad^{\bigodot} (Y\times Y)_+ =  \frac {r(n)}{ 2\sqrt 2+ r(n)} $$
and the normal focal radius of 
$$Y(2n+1) =\left ((c_1Y(n))\times (c_2Y(n+1))\right)_+\subset B^{2n+1}.$$
satisfies a similar relation.

Then,  starting from $Y(2)=S^1\subset \mathbb R^2$ with $r(2)=1$  one  
 obtains $Y(4)$, $Y(4)$, ... , such that 
 $$r(4) \geq \frac {1}{ 2^\frac {3}{2}+ 1}>\frac {1}{4}= \frac {2}{4^\frac {3}{2}},$$
 $$r(8) >    \frac {1}{ 8\sqrt 2+1}=  \frac {1} {\frac {1}{2} 8^\frac {3}{2}+1}>\frac {1}{13}, $$
and, in general, 
 $$r(n) \geq   cn^{-\frac{3}{2}},$$
  where a (very) rough estimate is  
  $c>1$  for $n=2^i$
  and 
  $ c>  \frac{1}{3}$ for all $n$.

 Eventually,  since the normal bands around these  tori $Y(n)$, can be   transported from $B^n$ to $S^n$  by the obvious expanding map $B^n\to S^n$, we conclude that  
 $$width_ {\hat{\cal T}}(S^n)\geq  2c\cdot n^{-\frac{3}{2}}$$ 
 which  combined with $\bigcirc^n_{Lip} $ implies the following.

  \vspace {1mm}

{\sf \textbf {Spherical Lipschitz Bound Theorem.}} {\sl If the scalar curvature of a  (possibly incomplete) Riemannian  $n$-manifold   
 is bounded from below by $n(n-1)=Sc(S^n)$, then  all continuous maps $f$  from $X$ to 
  the sphere $S^n$ (and also    to the hemisphere to $S^n_+$) of 
  non-zero degrees\footnote{  Here such a map  $X\to S^n$  is supposed to be constant at infinity, including $\partial X$ and to be proper from the interior of $X$ to that of $S^n_+$.}  satisfy 
 $$Lip(f)> \frac  {c}{\pi\sqrt n}\mbox { for the above  }  c>\frac {1}{3}.$$} 
 (This $c$ is not optimal; but since this inequality is unlikely to be qualitatively sharp anyway there is no point  in fiddling with constants.)

{\it Remarks.}  (a) If $X$ is a {\it complete spin}\footnote{In fact, it suffices to have the universal covering of $X$  spin -- we return to this later on'; here  we recall that   an orientable  smooth manifold $X$ is {\it spin}  if the restrictions of the tangent bundle $T(X)$ to all surfaces $S\subset X$ are trivial bundles.}  manifold, then  the    {\it sharp} spherical Lipschitz bound $Lip(f)\geq 1$ is known to hold for these maps $f: X\to S^n$  by the work of  Llarull [Ll 1998].    This is accomplished by  carefully analysing the  {\it  algebraic Schroedinger-Lichnerowicz-Weitzenboeck formula} for  {\it  the Dirac operator on $X$ twisted with the spin bundle 
$\mathbb S^+(S^n)$    pulled back to  $X$ }  and applying the index theorem. 

In fact, this  Dirac operator proof rules out smooth proper  maps $f:X\to U\subset S^n$ of non-zero degrees,  which {\it strictly decrease areas of surfaces} $S\subset  X$ (such $f$ may have $Lip(f)>>1$) and  where  the complements to the   (open) subsets $U\subset S^n$  are {\it  zero dimensional}, or, more generally, where  {\it all connected subsets $A\subset S^{n}\setminus U$ are trees and/or closed curves with trivial (i.e. identity) Levi-Civita monodromy transformations  around them} (see section 10).

(b) It remains unknown: \vspace {1mm}
 
$\bullet $   {\it if  the   spin  condition is essential  for ruling out maps $f$ for which $area(f(S))< area(S)$,
 
$\bullet $   if the completeness condition is essential for $ Lip(f)\geq 1$,

$\bullet$ if one may allow  closed curves in  $S^n\setminus U$ with nontrivial Levi-Civita monodromies  even if $X$ complete and spin}. (See section 10 for further questions of this kind.)

\vspace {1mm}

(c)  The  above  inequality $Lip(f)\geq \frac {c}{\pi\sqrt n}$  (which applies to incomplete non-spin manifolds)  improves  upon  
$Lip(f) \geq   \frac {n}{2^n\pi}$ in [GL 1983].\footnote {The inequality $Lip(f)\geq \frac {c}{\pi\sqrt n}$ with $c=1/3$  gains over  $Lip(f) \geq   \frac {n}{2^n\pi}$ only for $n\geq 6$  but a more precise    evaluation of  $rad^{\bigodot}(Y(n))$  shows  this gain for all $n$.}
 
 This  gains in significance   as  $n\to \infty$,  where the proof   for $n\geq 9$  depends on the controlled singularity results by Lohkamp and Schoen-Yau, which the present author has not studied in detail.

 \vspace {1mm}

(d) The above  estimates of torical width  of $S^n$ and of focal radii of  tori in $S^n$   raise a multitude of questions concerning  $width_{\hat {\cal T}}(X)$, $rad^{\bigodot} (Y\subset X)$ and their generalisations for various $X$ and $Y$. These  will be  briefly discussed in section  7.

\section {$\frac {4\pi}{n}$-Bound on Width  and   Related Inequalities for Iso-Enlargeable Bands.}

{\sf Hypersphericity and Iso-Enlargeability.}  An oriented Riemannian manifold  $X$  is called {\it  hyperspherical}  if it  admits  continuous maps $f$  to $S^n$, $n=dim(X)$ with arbitrarily small $Lip(f)>0$,  which are   {\it constant at infinity}  which have  {\it non-zero degrees}.

A  Riemannian manifold  $X$ is called {\it iso-enlargeable} if there exists  a sequence  of Riemannian manifolds $\tilde X_i$ of dimension $n=dim(X)$  and of  locally isometric  maps
$\tilde X_i\to X$, such that $ \tilde X_i$ admit continuous  maps constant at infinity 
$$f_i:\tilde X_i \to S^n, $$
such that  
$$deg(f_i)\neq 0 \mbox {  and $Lip(f_i)\to 0$   for $i\to \infty$}.$$  
 
{\it Examples.} (a)  The archetypical hyperspherical manifolds are   the Euclidean spaces  $\mathbb R^n$. 

(b) Complete simply connected manifolds $X$  with {\it non-positive sectional curvatures $\kappa$} are also 
hyperspherical. 

This follows from (a),  since the the inverse exponential maps $\exp^{-1}: X\to \mathbb R^n=T_{x_0}(X)$ satisfy 
$Lip(\exp^{-1}\leq 1 $  for $\kappa(X)\leq 0$.

(c) If  a compact manifold $X$ is fibered over an $\underline X$, where  $\kappa(\underline X)\leq 0$  and  where the fibers also admit metrics with $\kappa\leq 0$ then the universal covering of $X$ is hyperspherical by an easy argument.

(d) Compact locally symmetric spaces  $Y$ that have no (local) factors isometric to real and/or complex hyperbolic spaces are enlargeable but {\it not overtorical}, since  the homology groups  $H_1(Y)$  are   {\it finite} for these $Y$.

Instances of such $Y$ are  compact quotients $H^n_\mathbb H/\Gamma$  of  {\it quaternion hyperbolic spaces} (here the sectional curvature  $\kappa(Y)<0$)
and compact quotients   

\hspace {-6mm}$SO(n)\backslash SL(n)/\Gamma$ , $n\geq 3$  (here $\kappa(Y)\leq 0$).

\vspace {1mm}

{\it Remark/Question.}  If the, locally isometric maps    $\tilde X_i\to X$ in the definition of iso-enlargeability  are required to be {\it covering maps}, which is equivalent to {\it completeness} of $\tilde X_i$ in the case where $X$ itself is complete (e.g. compact), then $X$ is called {\it enlargeable}, see   [GL 1983],   [Dr 2000],  [DFW 2003], [HaSch 2006],  [BH 2009],   [Han 2011].
 
It is obvious that 
 $$enlargeable \Rightarrow iso\mbox{-}enlargeable,$$

Also one may expect that  the reverse implication  holds  for   compact manifolds,
since sequences $\tilde X_i$  
 (sub)converge in a natural way to some $\tilde X$, where the maps $\tilde X_i\to X$ (sub)converge to a covering map  $\tilde X \to X$ and where  properly scaled maps $\tilde X_i\to S^n$ (sub)converge to Lipschitz maps  $\tilde f_i:  X\to S^n$. 
 
 But, in general, these $\tilde f_i $ are neither constant at infinity nor do they have non-zero degree, at least not in the ordinary sense (even if  $\tilde f_i:\tilde X_i\to X$  were covering maps to start with).
 Thus  \vspace {1mm}
 
 {\sf enlargeability of compact iso-enlargeable manifolds remains problematic even 
 
 for  compact {\it aspherical}\footnote{ A manifold is called   {\it aspherical} if its universal covering is contractible.} manifolds $X$.  } \vspace {1mm}
 
 (Examples of enlargeable manifolds with non-hyperspherical {\it universal} coverings exhibited  in [BH 2009] tilts one  toward accepting  a  possibility  of  iso-enlargeable but non-enlargeable  compact manifolds $X$.)

  \vspace {1mm}
 
 On the other hand, there is the following   relation between iso-enlargeability and the  overtorical width $width_{\hat {\cal T}}(X)$  which was defined in  the previous section.

\vspace {1mm}

{\it If $X$ is compact, then 
$$ [width_{\hat {\cal T}}(X)=\infty]\Leftrightarrow [X \mbox {\sf is iso-enlargeable}].$$}

In fact,   the (quantitative form of the  obvious)  implication "$\Leftarrow$"  has been already established the previous section.

Now, to prove "$\Rightarrow$", we  observe that  the maps $f:V\to \mathbb T^{n-1}\times [-1,1]$
used in the definition of "over-torical" can be assumed Lipschitz, where, moreover,  the  corresponding maps (coordinate projections)    $ V\to [-1,1]$ can be arranged to have their Lipschitz constants   equal to $$\frac {2}{width (V)}.$$
  
  These $f$,   by passing to  the   $\mathbb Z^{n-1}$-coverings 
   $\tilde V\to V$, become Lipschitz  maps $\tilde f: \tilde V\to \mathbb R^{n-1}\times [-1,1].$ which, by scaling  $\varepsilon: \mathbb R^{n-1} \to \mathbb R^{n-1}$,  turn to maps 
   $$\tilde f_\varepsilon: \tilde V\to \mathbb R^{n-1}\times [-1,1]$$  
   with  Lipschitz constants arbitrarily close to
   $ \frac {2}{width (V)}$  and which remain proper with degrees $\neq 0$.
   
 Finally, we compose   these $\tilde f_\varepsilon$ with the obvious map   $\mathbb R^{n-1}\times [-1,1]\to S^n$ of degree
 one and $Lip= \pi$ and obtain  maps  
 $$\tilde F_\varepsilon: \tilde V\to S^n, \mbox  { where  $deg(\tilde F)\neq 0$ and $Lip(\tilde F)\leq  \frac {2\pi}{width (V)}+\varepsilon'$ }$$
 with arbitrarily small $\varepsilon'$, and the implication  
  $$ [width_{\hat {\cal T}}(X)=\infty]\Rightarrow [X \mbox {\sf is iso-enlargeable}].$$}
is thus established.\vspace {2mm}

  \hspace {30mm}{\large {\it  \textbf  {$\cal V$-Width and  $\cal IE$-Width.}}} \vspace {1mm}

  Given a class  $\cal V$  of Riemannian bands   $V$ define   $width_{\cal V}(X)$   of a Riemannian manifold $X$ as we did it for $width_{\hat{\cal T}}$, namely, as  
  
  {\sf the supremum of numbers $d$, such that $X$ admits  an equidimensional   locally isometric  (not necessarily globally one-to-one) immersion from a   band  $V\in \cal V$ with $width(V)= d$.}
  
 Here, we are concerned  with  the class of {\it iso-enlargeable orientable    bands} $V$ which admits  proper maps (i.e. boundary to boundary)  $f:V\to Y\times [-1,1]$,  where $Y$ must be compact orientable iso-enlargeable  manifolds without boundaries\footnote {A more general definition of iso-enlargeability for bands  with no reference to closed manifolds is given  in section 11.7.}  and  where  $deg(f)\neq 0$.

  \vspace {1mm}
  
  {\sf \textbf {Iso-enlargeable  $\frac {4\pi}{n}$-Inequality.}} {\it The iso-enlargeable widths of   $n$-dimensional  Riemannian  manifolds $X$ are bounded by the over-torical widths as follows.
 $$width_{\hat {\cal T}}(X)\leq   width_{\cal IE}(X) \leq 2width_{\hat {\cal T}}(X).$$
 
Consequently, if $Sc(X)\geq \sigma>0$, then
$$width_{\cal IE}(X) \leq 4\pi \sqrt{\frac {n-1}{\sigma n}}.$$}
  
  {\it Proof.} The inequality $width_{\hat {\cal T}}\leq   width_{\cal IE}$ is obvious.
 
 To prove  $width_{\cal IE} \leq 2width_{\hat {\cal T}}$ let us  show that iso-enlargeable bands $V$ with  width $d$ contain over-torical ones with width $d/2$.

 In fact, since the above  $Y$ is iso-enlargeable, there exist locally isometric immersions of  $(n-1)$-dimensional  over-torical  bands  $Y_D$ to $Y$ with  $width (Y_D)  \geq D$ for all $D>0$. Then the pullbacks\footnote {Even if  $g: A\to B$ is a 
 non-injective map, we  speak of the $f$-pullback of $A$ for   a map $f:C\to B$, where $f^{-1}(A)$ is understood as 
 the set  of pairs $\{a,c\}_{g(a)=f(c)}\in A\times C$    which comes with  the   map  $f^{-1}(A)\to C$,     $(a,c)\mapsto c$  which has the  same kind of (non)-injectivity as $g: A\to B$.} 
 of  $Y_D\times [-1,1] $ under the maps $f:V\to Y\times [-1,1]$ come with natural maps 
 $$f^{-1} (Y_D\times [-1,1])\to [0,D]\mbox {  and } f^{-1} (Y_D\times [-1,1])\to [0,d],$$
 both with $Lip\leq 1$.
 
Then the  pull back of the circle of radius $d/2$ in $[0,D]\times  [0,d]$ under the pair of these maps (which may be assumed smooth and transversal to this circle)  
  serves as the required overtorical  band of width $\geq d/2$. 
  
 Finally, we recall  the $\frac {2\pi}{n}$-inequality  for over-torical bands in section 2 and obtain 
our $\frac {4\pi}{n}$-inequality for iso-enlargeable bands.\vspace {1mm}
  
 \vspace {1mm} 
  
  {\sc Improvement.} {\sl The   $\frac {4\pi}{n}$-inequality  for  compact iso-enlargeable bands $V$   
 can be upgraded to
 $width_{\cal IE} (V)\leq 2\pi \sqrt{\frac {n-1}{\sigma n}}$}   as it is explained   at the end of   section 11. 6.\footnote{The proof 
in the original version of this paper was incorrect, as it was pointed out to me by the referee.}
Accordingly, the curvature  decay estimate below can be improved by the factor of 2.

 \vspace {1mm}

  {\sf \textbf Iso-enlargeable   {$[\asymp \frac {8\pi^2}  {R^2}]$ -Decay   Theorem.}}  Let a manifold $X$ admit a proper map of non-zero degree to  the total space $\underline  X$ of a   two dimensional  vector bundle  $\underline  X\to Y$ 
  where 
 $Y$ is  a compact 
  iso-enlargeable (e.g. admitting a metric with non-positive curvature)  manifold.

  {\it If   the bundle $\underline  X\to Y$ is trivial then
   the  scalar curvatures of  all complete Riemannian metrics  $g$ in $X$ restricted to  concentric balls  $B(R)=B_{x_0}(R)\subset  X$ satisfy
   $$  \min_{B(R)} Sc(X)\leq \frac{8\pi^2}  {(R-R_0)^2}  \mbox  { for some $R_0=R_0(X,g,x_0) $ and   all } R\geq R_0. \leqno {[\asymp \frac {8\pi^2}  {R^2}]}$$}\vspace {1mm}

  {\it Proof.}  This  follows word for word the argument  for  the  quadratic  decay   theorem in section 1 and its generalisation in section 2     with "iso-enlargeable " for  "over-torical".

   {\vspace {1mm}

   {\large \sf What happens to nontrivial bundles $\underline X\to Y$?} 
   The above argument applies to non-trivial bundles, where the (total spaces of the) corresponding 
   {\it circle bundles are iso-enlargeable}, which is so, for instance by the above (c)  for   $Y$ which  admit metrics with non-positive sectional  curvatures.
   
 In general, the examples  in [BH 2009] indicate a possibility of  non-enlargeable circle bundles over enlargeable $Y$; yet, it seems hard(er) to find such  examples, where  the corresponding $\underline X$ would admit complete metrics with $Sc\geq \sigma>0$. \vspace {1mm}

 
 {\it Remarks and Questions.}   
  ($\square$) Let $V$ be an $n$-dimensional  manifold with sectional curvatures $\kappa(V)\geq 1$ which   admits a proper  map $\Phi:V\to [0,d]^n$  given by $n$ functions $\phi_i(v)$ 
  with $Lip(\phi_i)\leq 1$ and such that  $deg(\Phi)\neq 0$.  
  
  Is  the maximal $d$ for these $V$   achieved  by the regular cube $\square$  in    the hemisphere $S^n_+$  with the boundary $\partial\square\subset S^{n-1}=\partial S^n_+$? 
  
  ($\triangle$) The same question for spherical simplices $\triangle\subset S_+^n$ with 
  $\partial\triangle\subset S^{n-1}$: 
  
  {\sl do these simplices have maximal distances between opposite faces among all simplices with $\kappa\geq 1$?}
  \vspace {2mm}
  
  \hspace {8mm} {\large {\it  \textbf  { From $\cal V$-manifolds to  to $\cal V$-Enlargeable ones.}}}\vspace {2mm}
 
{\large $[ \cal V \leadsto \cal VE]$:}
 Given a "natural"  class  $\cal V$ of manifolds one defines an, a priori larger,  class $\cal VE$  of $\cal V$-enlargeable manifolds $X$ by the condition
$$ width_{\cal V}=\infty.$$

Thus, for instance the class $\hat{\cal T}$ of over-torical manifolds   leads  to the class $ { \hat {\cal T}\cal E}\supsetneqq \hat{\cal T}$ of  $\hat{\cal T}$-enlargeable  manifolds, which, as we know, is equal to the class $\cal IE$ of  
iso-enlargeable   manifolds, 

On the other hand,  if we depart from the class $\cal IE={ \hat {\cal T}\cal E}$, then 
the new class $\cal IEE$ defined by $width_{\cal IE}=\infty$ will coincide with $\cal IE$.

In the following section, following Schoen-Yau and Schick, we define class $\cal SYS\supsetneqq \hat{\cal T} $, where the corresponding class $\cal SYSE$
 of $\cal SYS$-enlargeable manifolds is {\it strictly greater than the class of iso-enlargeable ones.}

\section{ Schoen-Yau-Schick Manifolds  and $\cal SYS$-Bands.}

 {\it  Schoen-Yau Definition.}  [SY-Str 1979], [SY 2017].   A compact orientable  $n$-manifold $X$ is {\it SYS}, if there exist  $n-2$  integer homology  classes $h_1,h_2,...,h_{n-2} \in H_1(X)$,  such that their intersection consecutive
  $$h_1\frown h_2\frown...\frown h_{n-2}\in H_2(X)$$
   is {\it non-spherical}, i.e. it  is {\it not contained} in the image of the {\it Hurewicz homomorphism} $\pi_2(X)\to H_2(X)$, or, equivalently, it doesn't lift to the universal covering of $X$. 
 \vspace{1mm}
  
{\it Schick Definition.} [Sch 1998]  A homology class $h\in H_n(K)$, where   $K= K(\Pi,1)$ is the  
Eilenberg-MacLane space for an Abelian group $\Pi$, is called $SYS$, if its consecutive  {\it cap-producs} 
with some   cohomology classes $h_1, h_2,...,h_{n-2}\in H^1(K,
\mathbb Z)$ are non-zero,  
$$(...((h\cap h_1)\cap h_2)\cap....\cap h_{n-2}) =  h \cap(h_1\smile,...,\smile h_{n-2})\neq 0\in H_2(K).$$
(Geometrically speaking,  generic  2-dimensional  intersections of the  $n$-cycles $C\subset K$  representing $h$ with
$(n-2)$-codimensional pullbacks of generic  points of,  some, say piecewise linear,   maps $K\to\mathbb T^{n-2}$ are non-homologous to zero.)

Then a manifold $X$ is SYS if the Abel  classifying map $X\to K(\Pi,1)$ for $\Pi=H_1(X)$ sends the fundamental class $[X]\in H_n(X)$ to a SYS class in this $K(\Pi,1)$.

(Recall that, by definition,  the spaces $K(\Pi,1)$ have {\it contractible} universal coverings and fundamental groups {\it isomorphic} to $\Pi$.  The standard finite dimensional approximations to these $K$ are products of tori and {\it lens spaces} 
{\sf $L_i$}$ =S^N/\mathbb Z_{l_i}$, where the latter, observe, carry natural metrics with $Sc>0$.

 Abel's $X\to K$ maps, which  are unique up-to homotopy,  are characterised by inducing isomorphisms on the 1-dimensional homology groups.)

\vspace {1mm}

 {\it Historical Remark.}  In 1979 Schoen and Yau proved that SYS manifolds (defined slightly differently in [SY - Str 1979]  with incorporation of some  spin  manifolds) of dimensions $n\leq 7$    {\it carry no metrics with $Sc>0$.} Then, in the recent paper  [SY 2017], they published the proof for   all $n$.
 
 Meanwhile, Schick [Sch 1998] has shown that {\it no available Dirac operator methods can rule out $Sc>0$ on these manifolds.}
\vspace {1mm}

\hspace {40mm}  {\Large \sf  Examples.}\vspace {1mm}

$\bullet_1$ Overtorical manifolds are $SYS$.\vspace {1mm}

$\bullet_2$  Let $X$ be obtained by a  surgery applied on a closed curve $C$ in the $n$-torus as in   [Sch1998].

If $n\geq 4$, then 
 $X$ is SYS if and only if   $C$ represents  a {\it divisible} homology class in $H_1(\mathbb T^n)$.  

(Such an $X$  is over-torical if and only if $C$ is {\it homologous to zero}.)\vspace {1mm}

$\bullet_3$ If a compact orientable manifold $X$ admits a map $f$  of degree one to a SYS manifold that 
$X$ is SYS.

But if $deg(f)>1$ then $X$ is not necessarily SYS, unlike the case of the overtorical and iso-enlargeable manifolds. For instance if the curve $C$ in  $\bullet_2$  is $m$-divisible, than the 
some $m$-sheeted covering of $X$ is non-SYS.

{\sf Probably, these non-SYS coverings carry metrics with $Sc>0$.}\vspace {1mm}

$\bullet_4$ Products of SYS manifolds by overtorical ones are SYS.

But products SYS $ \bigtimes$  SYS  and 
SYS $\bigtimes$ [{\sf iso-enlargeable}] are, in general, not SYS. \vspace {2mm}

{\it SYS-Bands.}
A band $V$ is called  $SYS$  if it  admits a  band map ($\partial_\pm \to \partial_\pm$)  
of degree $\pm1$ to $Y\times [-1,1]$  where $Y$ is a compact SYS manifold.\footnote {A more general definition of SYS  for bands  with no reference to closed manifolds is given  in section 11.7.}

  Accordingly, define the $\cal SYS$-width  of $width_{\cal SYS}(X)$  of Riemannian manifolds $X$ based on the class $\cal SYS$ as we did it  for $\cal IE$ in the previous section.

\vspace{1mm}

 {\sf \textbf { $\frac {4\pi}{n}$-Inequality for SYS-Bands.}} {\it All Riemannian manifolds  $X$ with $Sc(X)\geq \sigma>0$ satisfy
$$width_{\cal SYS}(X)\leq  4\pi\sqrt{\frac {n-1}{\sigma n}}.$$

Consequently,    compact manifolds without boundaries, which have   
 $$width_{\cal SYS}(X)=\infty$$
  admit no metrics with positive scalar curvatures.}

\vspace {1mm}

{  Proof.}
By symmetrising   a   $SYS$-band $V\to Y\times [-1,1]$ as   in  the proof of the overtorical  $\frac {2\pi}{n}$-Inequality in section 2  (now  $Y$ plays the role of the torus $\mathbb T^{n-1}$ in section 2)  we arrive at $V_{\circ^{n-3}}$  with  $\mathbb T^{n-3}$-invariant metric with $Sc\geq \sigma$, such that\vspace {1mm}

{\sl the quotient space $\underline V^3=V_{\circ^{n-3}}/\mathbb T^{n-3}$
is an orientable  $3$-manifold with the boundary decomposed into  two (possibly disconnected) disjoint  parts say
$$\partial \underline V^3= S_-\cup S_+,$$  
 where
$$dist_{\underline V^3}(S_-,S_+)\geq d$$
for $d$ equal to the distance between the  two boundary components in $V$,}\vspace {1mm}

\hspace {-6mm} and where the Schoen-Yau-Schick property of $Y$  implies that  \vspace {1mm}

{\sl if a closed surface $S\subset \underline V^3$ separates $S_-$  from  $S_+$,  then the homomorphism 
$$\pi_1(S)\to \pi_1(\underline V^3)$$
has infinite image.}

Therefore  the $d/2$-equidistance  surface to $S_-$  (or to $S_+$) contains a circle $C$ which has infinite order in $\pi_1(\underline V^3)$  and, by the Poincar\'e duality,  the covering   $\tilde {\underline V^3}$ of $\underline V^3$ with the cyclic  $\pi_1(\tilde {\underline V^3})$  generated by the (homotopy class of) $C$ contains  a relative  $2$-cycle $\tilde C^\perp$ \footnote {Relative means relative to $\partial \tilde {\underline V^3} +\partial_\infty\tilde {\underline V^3}$  where  $\partial_\infty $ stands for the complement of a large ball in $\tilde {\underline V^3}$.}  with non-zero intersection index with the lift  $\tilde C$ of $C$ to  $\tilde {\underline V^3}.$

Take the pull back  of the cycle  $\tilde C^\perp$  to the corresponding covering $\tilde V_{\circ^{n-3}}$ of $$V_{\circ^{n-3}}=(V_{\circ^{n-3}}/ \mathbb T^{n-3})\times  \mathbb T^{n-3},$$ write this  pullback cycle     as  
$$\tilde C^\perp\times \mathbb T^{n-3}\subset \tilde V_{\circ^{n-3}}, $$ 
and symmetrize the minimal  cycle in the 
$(n-1)$-homology  class of  $\tilde C^\perp\times \mathbb T^{n-3}.$

Since  $dist (\tilde C, \partial \tilde {\underline V^3}) =dist (C, \partial   \underline V^3)\geq d/2$, the quotient surface of the resulting $\tilde V_{\circ^{n-2}}$ contains a point within distance $\geq d/2$ from its boundary, which implies (compare p. 310 in [GL 1983]) that
$$d/2 \leq 2\pi\sqrt{(n-1)/\sigma n}.$$
QED.

\vspace {1mm}

{\large \sf Question. } Can one   replace the above $4\pi\sqrt{\frac {n-1}{\sigma n}}$ by $2\pi\sqrt{\frac {n-1}{\sigma n}}$? 

\section {$\cal SYS$-Enlargeable Manifolds and Codimension Two Depth Inequalities.} 

A Riemannian manifold  $X$ is called {\it $\cal SYS$-Enlargeable}  if it has {\it infinite} 
$\cal SYS$-width.

 For instance, SYS manifolds  and iso-enlargeable manifolds  are {\it $\cal SYS$-Enlargeable}.

What 
is more interesting  is that  \vspace {1mm}

{\sl if an  $n$-manifold $X$ admits a proper Lipschitz  map}  $\phi$ ({\it Lipschitz} means  $Lip(\phi)<\infty$)   {\sl to an  iso-enlargeable manifold  of dimension $n-2$, say $\phi:X\to  \underline X$,
such that the homological  pullback   $\phi^![\underline x]\in H_2(X)$,  $[\underline x]=1\in H_0(\underline X)=\mathbb Z $,  is  non-spherical (as in the first  definition of SYS in  the previous  section), then $X$  is $\cal SYS$-enlargeable.} \vspace {1mm}

  Therefore, by  the above  $\frac{4\pi}{n}$-inequality,  \vspace {1mm}

 {\it If such an $X$ is compact, then it  admits no metric with $Sc>0$.}\vspace {1mm}

\vspace {1mm}

Thus, for example,

{\it products $X$ of SYS manifolds  by compact iso-enlargeable ones (e.g. those which  admit  metrics with   $\kappa(X_2)\leq 0$)     admit no metrics with positive scalar curvatures.}

(These $X$, in general, are  neither iso-enlargeable nor SYS.)
\vspace {1mm} 

{\sf \textbf { $\frac {8\pi}{n}$-Inequality for SYSE-Bands.}}  Denote by  $\cal SYSE$  the class of SYS-enlargeable manifolds,  say that  a compact band  $V$ is  SYSE  if it admits a map of  degree $\pm1$ to $Y\times [-1,1]$, where $Y$ is SYSE and accordingly define $width_{SYSE}(X)$ for   Riemannian manifold $X$ (see {\large $[ \cal V \leadsto \cal VE]$} in section 4).   

Then by arguing as in the proof of  the 
iso-enlargeable $\frac {4\pi}{n}$-inequality in section 4
we conclude that 
$$ width_{\cal SYS}(X) \leq  width_{\cal SYSE}(X) \leq 2width_{\cal SYS}(X)$$
for all Riemannian manifolds $X$.

Consequently, 

{\it if $Sc(X)\geq \sigma>0$ then
 $$width_{\cal SYSE}(X)\leq 8\pi\sqrt \frac {n-1}{\sigma n}.$$}

 {\it Question.} Can one improve $8\pi$ to $2\pi$ or, at least, to $4 \pi$?
 \vspace {1mm}
 
 {\sf \large \it  Depth Inequalities.} Define {\it the depth of a homology class} $h$ in a Riemannian manifold $X$ with boundary as the supremum of  $d\geq 0$ such that $h$ can be represented by a cycle positioned  within distance $\geq d$ from the boundary of $X$ . (If $X$ is incomplete, we include  the   points obtained by completion of $X$ in  the boundary of $X$.)

Let $Y$  be a closed  $(n-2)$-dimensional  manifold  and  $p:\underline X\to Y$ be    a  disc bundle, 
e.g. the trivial one  $\underline X=Y\times  B^2$.

Let  $X$ be a compact $n$-manifold with boundary  and  $f: X\to \underline X$  be a  proper continuous map   where {\sl proper}, means    {\sl boundary $\to$ boundary}. Let $h=f^!([Y]\in H_{n-2}(X) $ be the homology pull-back  of the  homology class of the zero section $Y=Y_\mathbf 0\subset \underline X$.

Let  $\underline X_{-\varepsilon}\subset \underline X$ be the complement of the open   $\varepsilon $-neighbourhood of $Y_\mathbf 0$ in  $\underline X$ and observe that the boundary of $  \underline X_{-\varepsilon} $  consists of two components, call them  
$\underline \partial_\pm$  which are canonically homeomorphic to the total space of the circle bundle  associated  to $\underline X\to Y$,   denoted $p_\circ: Y_\circ\to Y$.

 Let $ \partial_\pm= \partial_\pm(X)\subset \partial X$ be the two parts of the boundary of $X$ which are sent by the map  $f:X\to \underline X$ to $\underline \partial_+ $ and to $ \underline \partial_-$ correspondingly.
 
 Observe that 
 
 $\bullet$  $\underline X=Y_\circ\times [\varepsilon,1]$;
 
 $\bullet$ If $Y$ is iso-enlargeable then  $Y_\circ$ is also iso-enlargeable.

 $\bullet$ if the fibration  $p: \underline X =Y$ is trivial,   $\underline X =Y\times B^2$, and

  if  $Y$  is over-toric then also
 $Y_\circ$ is over-toric,  
 
 if 
  $Y$ is  SYS then 
 $Y_\circ$ is also SYS, 
 
 if $Y$ is  SYSE then 
 $Y_\circ$ is also SYSE.\vspace {1mm}
 
Now \vspace {1mm}

{\it let  the fibration $p:\underline X\to Y$ be  trivial}
and let \vspace {1mm}

\hspace {43mm}{\it $Sc(X)\geq n(n-1) =  Sc(S^n)$}. \vspace {1.6mm}

 Observe that 
the band-width $\frac {k\pi}{n}$-inequalities, (for $k=2, 4,8$ see  sections 2,4,5) imply the following bounds on  the depths of $h\in H_{n-2}(X)$
 by the argument that we have already used   several times, e.g. in the proofs of the  quadratic decay inequality in section 1.

\vspace {1mm}


{\sf  $[{\hat{ \cal T}}]_\circ$}  {\sl If   $Y$ is over-torical, i.e.  if it admits a map to the  torus $\mathbb T^{n-2}$  with degree $\neq 0$, then  
$$ depth (h)\leq \frac {2\pi}{n}.$$}
 {\sf This is the only case where our inequality is (known to be) sharp},

  \vspace {1mm}

   $[{{ \cal IE}}]_\circ$  {\sl  If $Y$ is iso-enlargeable, e.g. if it admits a metric with non-positive sectional curvature, then 
$$ depth (h)\leq \frac {4\pi}{n}. \footnote {This can  improved to $depth (h)\leq \frac {2\pi}{n}$ with the corresponding width inequality in section 11.7  , which, probably can be also used to improve the width  inequality in  the SYS case as well.}$$}
(Here the fibration $p$ need not be trivial.)
  \vspace {1mm}

 $[{{ \cal SYS}}]_\circ$   {\sl If  $Y$ is  SYS and if the map $f:X\to\underline X$  has $deg (f)=\pm1$, 
 then  
$$ depth (h)\leq\frac {4\pi}{n}.$$}
(The simplest  example of  a non-overtoric SYS manifold  $Y$  for $n-2\geq 4$  is obtained from the  $(n-2)$-torus by attaching a 2-handle based on a $k$-multiple of closed curve in this torus where $k\neq\pm1$. In this case one only need $deg(f)$ to be non-divisible by $k$.)

\vspace {1mm}

 $[{{ \cal SYSE}}]_\circ$ {\sl If   $Y$ is  SYSE and if the map $f:X\to\underline X$  has $deg (f)=\pm1$, 
 then  
$$ depth (h)\leq\frac {8\pi}{n}.$$}

(Recall, this was stated earlier, here as everywhere in this paper  the above inequalities  are established  unconditionally for $n\leq 8$,  while the case  $n\geq 9$  relies on the recent  partial  regularity results  by Lohkamp and by Schoen and Yau which the present author has not studied in detail.)
\vspace{1mm}

{\it On nontrivial bundles   $p:\underline X\to Y$.} Here,    similarly to  where we  addressed this issue  in section 4, one may drop the triviality of $p$  assumption, if, for instance, $Y$ admits a metric with 
$\kappa\leq 0$.

No reasonable   assumption of this kind,  however, seems in view for   SYS and SYSE manifolds.

In fact,   circle bundles over many  SYS manifolds, say on those obtained by surgery on closed curves in  $\mathbb T^n$ (see $\bullet _2$ in  section 5) are very likely to  carry metrics with $Sc>0$ and so the above inequality  can't hold  with any constant for non-trivial fibrations $p:\underline X\to Y$.
  \vspace {2mm}

{\it On  Complete Manifolds and Dirac Operators.}  The inequality
$depth(h)<\infty$ implies  that the interiors of the   manifolds $X$ in   $[{{ \cal IE}}]_\circ$ and    $[{{ \cal SYSE}}]_\circ$ \vspace {1mm}

\hspace {18mm} {\it admit no complete metrics $g$ with $Sc(g)\geq \sigma>0$}.\vspace {1mm}

(The inequality $depth(h)<\infty$ in the  remaining   cases follow from these two.)
 \vspace {1mm}
 
 Strangely enough,  even if $X$ is spin, this was proven  by the Dirac operator methods for {\it enlargeable and related} manifolds $Y$   [GL1983], [HPS2015]   {\it only under  additional geometric assumptions} on $X$ in  spirit of "bounded geometry".
 
(To be honest, I am not 100\% certain this is the case for  [HPS2015]  . The main result is stated in this paper for closed manifolds and  I   have not followed  the proofs in sufficient  details to understand what is actually proven there for complete non-compact  manifolds.\footnote {Bernhard Hanke told me that the results from   [HPS2015]  apply to complete manifolds by the $C^\ast$-arguments in Roe's partitioned index theorem and  in    [HS 2007].})

\vspace {1mm}


{\it Question.} {\sl  Do  all  products manifolds  $Y\times \mathbb R^2$, and, more generally, the  total spaces  of all $\mathbb R^2$-bundles  admit complete metrics $g$ with $Sc(g)\geq 0$? 

Do, for example, such metrics $g$ exist for  compact manifolds  $Y$ which admit metrics with strictly negative  sectional curvatures?}\vspace {1mm}

If there are no  such $g$ among rotationally symmetric warped product metrics,\footnote{Figuring  this out  does not seem hard,  but I have not tried  doing this.} then, probably,   no complete  metric $g$ on $Y\times \mathbb R^2$ has $Sc(g)\geq 0$, where  the best candidates  of this kind   of manifolds  with no complete  metrics  on them  with $Sc\geq 0$ are non-trivial $\mathbb R^2$-bundles over surfaces of genera $\geq 2$.
\vspace {2mm}

\section {External Curvature, Focal Radius and Depth  in Codimension$>2$}

 Observe that by Gauss theorema egregium  the scalar curvature of hypersurfaces $Y\subset S^n$, $n\geq 2$, with principal curvatures $c_i=c_i(y)$,  $y\in Y$, $i=1,..., n-1$, satisfies 
$$Sc(Y)=Sc(S^{n-1}) + \left(\sum_i c_i\right)^2 -\sum_i c_i^2 \geq (n-1)(n-2)-  \sum_i c_i^2.$$ 

It follows that if an $(n-1)$-dimensional manifold  manifold $Y$ {\it admits no metric with $Sc>0$}, 
that the suprema of the principal curvatures of all  smooth immersions from $Y$ to the unit sphere $  S^n$ satisfy 
 $$\sup_{i, y}|c_i(y)|\geq \sqrt {n-2}$$

This is significantly weaker then the  $\frac {\pi}{n}$-inequality for  the normal radius of  $\mathbb T^{n-1}\subset S^n$,  which implies that $\sup_{i, y}c_i(y)\geq \frac  {(1+\varepsilon_n)n}{\pi}$.
But it applies to such manifolds, for instance, as certain exotic spheres $Y$ of dimensions $8m+1$ and $8m+2$ which carry no metrics with $Sc>0$  by a theorem of Hitchin [H1974], yet are   immersible (but not embeddable!)\footnote{ According to   textbooks' terminology, a smooth map $A\to B$ is an {\it immersion} if it is {\it locally} one-to-one and the inverse map is smooth, while {\it embeddings} are  immersions which are globally one-to-one and, if $Y$ is non-compact,  are  additionally required to be homeomorphisms from $Y$ to their (possibly, non-closed in B) images.} to $S^n$ by Smale-Hirsch theorem.

Besides, this $\sup c_i$ inequality obviously  generalises to  $Y$ in $S^n$ of all codimensions $k$ where
it reads  
$$\sup_{i, j,y}|c_{i,j}(y)|\geq \frac  {\sqrt {n-k-1}}{k}, \hspace {1mm} i=1,...,n-1, \hspace {0.5mm}j=1,...,k,\hspace {0.5mm} y\in Y,$$
 for all $Y$ which admit no metric with $Sc>0$.

Then, obviously, the same holds true for Riemannian manifolds $X\supset Y$  with {\it sectional curvatures $\kappa\geq 1$.}

More interestingly, a similar inequality     holds for  immersions to unit  Euclidean balls  $B(1)\subset \mathbb R^n$
Namely,

\vspace {1mm} 

{\it if an $(n-k)$ dimensional $Y$ admits no metric with $Sc>0$, then the principal curvatures of  all smooth immersions 
$Y\to B(1)\subset  \mathbb R^n$ are bounded from below by  
$$\sup_{i,j, y}|c_{ij}(y)|\geq 
\frac {1}{const} \frac{\sqrt {n-k-1}}{k}  $$
for some universal positive  constant $const \leq 100$. }

\vspace {1mm} 

{\it Proof.} The Euclidean case     reduces to the spherical one, since the  standard projective map $\mathbb R^n \supset B(1)\to S^n$  distorts curvatures of the curves in $B(1)$ by a bounded amount.
(Compare Lemma (C$'$) in  3.2.3  in [Gr 1986]).

\vspace {1mm} 

{\it Remark.} This  $\sup c_{i,j}$-inequality also holds in the balls in the hyperbolic spaces with sectional curvature $\kappa=-1$.

Also, the following   weaker form of this inequality holds for the unit balls in  all  $n$-dimensional Riemannian manifolds $X$ with 
  $-1\leq  \kappa(X)\leq 1$.
$$\sup_{i,j, y}|c_{ij}(y)|\geq 
\frac {1}{const} \frac{\sqrt {n-k-1}}{k}-const'.$$

In fact -- this is obvious by   today's standards --    the exponential maps  $exp : T_x(X) \subset B(1)\to X$ in these $X$ can be approximated by maps with controlled distortion of curvatures of the curves in  $B(1)$.

\vspace {1mm}

{\it Discussion.} There is a huge gap between   the above  lower  bounds on the curvatures of submanifolds in $S^n$ (and/or in $B(1)\subset \mathbb R^n$) and the observed  curvatures in the available examples $Y\subset S^n$.

Probably,  certain  homogeneous submanifolds  $Y\subset S^n$, such as  \vspace {1mm}

{\sf $\bullet $ real and complex projective spaces {\it Veronese} represented by symmetric/Hermitian forms of rank one,
 
$\bullet $  Grassmannians  {\it Pl\"ucker} embedded to  exterior powers of linear spaces,  
 
$\bullet $  the same Grassmannians represented by  projectors in spaces of operators, \vspace {1mm}}

  \hspace {-6mm} give a fair idea of embeddings with economical $c_{ij}$.

For instance,  the   curvature of the obvious embedding  of the  product of spheres
$$Y= S^{n_1}\times S^{n_2}\times ... \times S^{n_k}\subset S^{n_1+n_2+... +n_j + j-1}=\partial B(1)\subset \mathbb R^{n_1+n_2+... +n_j + j}$$ 
has $\max c_{ij}= \sqrt k$  and it is plausible (?) that   \vspace {1mm} 

{\sl  no  embedding/immersion of this $Y$ to $S^n$ may have a  (significantly)  smaller curvatures $c_{ij}$.\vspace {1mm}}

Notice that above local bound  $\max c_{ij} \gtrsim \sqrt n/k$ is non-vacuous only if all  spheres are one dimensional, 
while the 
 only known  improvement of this bound 
  is  the inequality   $ \max c_{ij}  \gtrsim n$ which  was established  in the  previous section only    for codimensions 1 and 2  and only for   $\cal SYSE$-manifolds $Y$ (e.g. for $Y$ which admits metrics with $\kappa\leq 0$.) \vspace {1mm}

This, for instance,  leaves the following questions open. \vspace {1mm}

(a)  {\sl Does  the  torus 
 $$ {\underset {n-3}{\underbrace {S^{1}\times S^{1}\times ... \times S^1}}}$$
embed to $S^n$  with the  principal curvature $c_{ij}\leq 100/n$?} \vspace {1mm}

(b)   {\sl Does the product 
 $$ {\underset {n-3}{\underbrace {S^{1}\times S^{1}\times ... \times S^1}}}\times  S^2$$
embed  to $S^n$  with the  principal curvature $c_{i}\leq 10$? }\vspace {1mm}

In fact, we are more interested  in {\it depth of homology and cohomology classes}   in Riemannian manifolds $V$ rather than in their curvatures,  where,  \vspace {1mm}

{\it by definition},
 $depth(h)\geq d$ for an $h\in H^\ast(V)$ if {\it the restriction of $h$ to the subset $V_{-d}\subset V$}
of the points within distance $\geq d$ from the boundary of $V$ (including the infinity for non-compact $V$, as in the previous section)   {\it does not vanish}. \vspace {1mm}

{\large \sf Problem.}  {\sl Bound  "complexity" of an $h$ in terms of  $d= depth(h)$.}

For instance, let   the sectional curvature of $V$ be bounded from below by $\kappa(V) \geq 1$ and let $h$ be induced by a continuous map from the fundamental cohomology  class of a  product of spheres,  
$$h=f^\ast [Y]\mbox  {  for }  f: V\to  Y= S^{n_1}\times S^{n_2}\times ... \times S^{n_k}.$$

{\sl Does the  depth of  $h$ necessarily  tend to zero for $k\to \infty$? }

\section {Symmetrization of Riemannian Manifolds  with Point-wise Control  of the Scalar Curvature.}

{\large \sf   Step 1 in   Symmetrization by Reflections.} Let $X$ be  a compact Riemannian manifold,  with a (possibly empty) boundary, and let $(Y_0, \partial Y_0)  \subset( X, \partial X)$ be a    cooriented  hypersurface  which is {\it strictly  locally volume minimising},  i.e.    all sufficiently close to $Y$ hypersurfaces
$(Y,\partial Y)\subset (X, \partial X)$ different from $Y_0$ satisfy
 $$vol_{n-1}(Y)>vol_{n-1}(Y_0).$$ 

Let $U\subset X$  be a (small)  neighbourhood of $Y_0$ in $X$ which is divided by $Y_0$ 
into two "halves",  denoted $U_\pm\subset U$, and let $Y_{\pm \varepsilon}\subset U_\pm $ be hypersurfaces homologous to $Y_0$ in $U_\pm$ which minimise
the functionals
$$Y\mapsto vol_{n-1}(Y) -\varepsilon\cdot vol_n(U_{\pm\varepsilon})$$ where
$$U_{\pm\varepsilon} = U(Y_{\pm\varepsilon}) \subset U_\pm$$
denote the regions bounded by $Y_{\pm\varepsilon}$ and $Y_0.$

By the basic regularity theorems of  Simons-Federer-Almgren-Allard   these $Y_{\pm\varepsilon}$ do exist for small $\varepsilon\geq 0$ and they are smooth
 away from closed subsets of Hausdorff codimension $\geq 7$.

 Moreover, if $n=dim(X)=8$, then, according to  [NS 1993],\footnote{What we need is not formulated in  [NS 1993]  but the argument from [NS 1993] does apply.} these $Y_\pm(\varepsilon)$ everywhere  smooth for an open dense set of $\varepsilon>0$ [NS 1993].\footnote{I am not certain this is formulated in  [NS 1993]. }

The mean curvatures of all  these  $Y_{\pm\varepsilon}$ at the regular points satisfies
$$mean.curv(Y_{\pm\varepsilon})=\varepsilon$$
and  the dihedral angles between the tangent spaces to $Y_{\pm\varepsilon}$ and those    to $\partial X$  at all regular points of $Y_{\pm\varepsilon}$ on the boundary $\partial Y_{\pm\varepsilon}$  are $\leq \frac {\pi}{2}$.  (If the boundary $\partial X \subset X$ is totally geodesic then $Y_{\pm\varepsilon}$ is normal to $\partial X$.)

 \vspace {1mm}
 
 {\sf \large On non-strictly minimal $Y_0$.}  If $Y_0$ is {\it non-strictly} volume minimising, then there are   hypersurfaces in $X$  with the same volume as $Y_0$, which lie  $\varepsilon$-close to $Y_0$ for all small $\varepsilon$.  These do not intersect $Y_0$ and each of them lies   on one side  of $Y_0$, where it plays  the role of  $Y_{-\varepsilon}$ or of  $Y_{+\varepsilon}$.   
 
 Alternatively, one may slightly perturb the metric in $X$, such that $Y_0=Y_{0, \epsilon}$ becomes strictly minimising. Then this  $\epsilon$, which  goes through the following stages of symmetrisation,  is sent  to  $0$ at the end of the symmetrization process.

  \vspace {1mm}

Let $$U_{{[-\varepsilon,\varepsilon]}}= U_{-\varepsilon} \cup U_{+\varepsilon}$$
and let $\tilde U_{{[-\varepsilon,\varepsilon]}}$
be  obtained by reflecting $U_{{[-\varepsilon,\varepsilon]}}$ in   the two parts   $Y_{\pm\varepsilon}$ of the (relative) boundary of $U_{{[-\varepsilon,\varepsilon]}}$ in $X$ (i.e. with the exclusion of $U_{{[-\varepsilon,\varepsilon]}}\cap \partial X$), denoted 
$$\partial_\pm U_{{[-\varepsilon,\varepsilon]}}= Y_{-\varepsilon}\cup Y_{+\varepsilon}.$$

In other words,  $\tilde U_{{[-\varepsilon,\varepsilon]}}$  is a space, which is acted upon by the semidirect product group $\Gamma=\mathbb Z\rtimes \mathbb Z_2 $,  such that 
$$\tilde U_{{[-\varepsilon,\varepsilon]}}/\Gamma= U_{{[-\varepsilon,\varepsilon]}}, \leqno {\hspace {5mm}\bullet}$$

$\bullet$  there is an embedding  $E:U_{[-\varepsilon, \varepsilon]} \hookrightarrow \tilde U_{{[-\varepsilon,\varepsilon]}}$ which is inverse to the quotient map  $Q:\tilde U_{{[-\varepsilon,\varepsilon]}}\to U_{[-\varepsilon, \varepsilon]}$, 
$$ Q\circ E  =Id: U_{[-\varepsilon,\varepsilon]}\to U_{[-\varepsilon,\varepsilon]},$$

$\bullet$ the group $\Gamma$ is generated by two involutions  (reflections) of 
$\tilde U_{[-\varepsilon,\varepsilon]}$, one of them fixing   $E(Y_ {-\varepsilon} )\subset  E(\partial_\pm   U_{{[-\varepsilon,\varepsilon]}})$ and the other one $E(Y_ {+\varepsilon})\subset E( \partial_\pm   U_{{[-\varepsilon,\varepsilon]}}).$

Thus,  the action of  our $\Gamma=\Gamma_\varepsilon $  on  
$\tilde U_{{[-\varepsilon,\varepsilon]}}$ mimics   the action of the same  group   on the line $(\infty, \infty)$, which  is   generated by  the transformations 
 $$ t\mapsto \pm \varepsilon-t$$
 and where $\tilde U_ {{[-\varepsilon,\varepsilon]}}$  admits a $\Gamma$-equivariant map to 
 $(\infty, \infty)$, such that the pullback of ${[-\varepsilon,\varepsilon]}\subset (-\infty,\infty)$  is equal to  $ E(U_{{[-\varepsilon,\varepsilon]}})\subset \tilde U_{{[-\varepsilon,\varepsilon]}}$.
 
 In particular, the action of the group $\mathbb Z=2\varepsilon \mathbb Z\subset \Gamma$  on $\tilde U_{[-\varepsilon,\varepsilon]}$ is free  and    the quotient space is  equal  to the double of $U_{[-\varepsilon,\varepsilon]}$,
 $$ \tilde U_{{[-\varepsilon,\varepsilon]}}/2\varepsilon \mathbb Z=U_{{[-\varepsilon,\varepsilon]}}\bigcup_{\partial_\pm  U_{[-\varepsilon,\varepsilon]}} U_{{[-\varepsilon,\varepsilon]}}, $$
where     the boundary of this double is   the  double of the region  
 $U_{[-\varepsilon,\varepsilon]}\cap \partial X$   across the  boundary  of this region in $\partial X$,
 $$\partial \left (U_{{[-\varepsilon,\varepsilon]}}\bigcup_{\partial_\pm  U_{[-\varepsilon,\varepsilon]}} U_{{[-\varepsilon,\varepsilon]}}\right)= \left (U_{[-\varepsilon,\varepsilon]}\cap \partial X\right)\bigcup_{\partial'}  \left (U_{[-\varepsilon,\varepsilon]}\cap \partial X\right)$$
for $$\partial'= \partial \left (U_{[-\varepsilon,\varepsilon]}\cap \partial X\right)= \partial_\pm U_{[-\varepsilon,\varepsilon]}\cap \partial X.$$

The Riemannian metric on    $U_{[-\varepsilon,\varepsilon]}\subset X$, that is the restriction of 
the metric of $X$ to   $U_{[-\varepsilon,\varepsilon]}\subset X$,  naturally induces  a $\Gamma_\varepsilon$-invariant path metric in $\tilde U_{[-\varepsilon,\varepsilon]}$, call it $\tilde g_\varepsilon$; if  
 the hypersurfaces $Y_{\pm \varepsilon} \subset X$ are non-singular, e.g. if $n=dim (X)\leq  7$, then this metric is $C^0$-Riemannian.

And if, moreover,   the minimal hypersurface $X_0\subset X$ is {\it non-singular}, e.g. for  $n\leq 7$, then the spaces  $(\tilde U_{[-\varepsilon,\varepsilon]},  \tilde g_\varepsilon)$ Hausdorff converge to a smooth Riemannian manifold, which we denote
$$( \underrightarrow {\tilde U}_{ 0},\underrightarrow {\tilde g}_{ 0})=\lim_{\varepsilon \to 0} (\tilde U_{[-\varepsilon,\varepsilon]},  \tilde g_\varepsilon),$$ 
which is  isometrically acted upon by the  group  $\mathbb R\rtimes \mathbb Z_2 = \lim_{\varepsilon\to 0} \Gamma_\varepsilon$  for the above $\Gamma_\varepsilon \subset \mathbb R\rtimes \mathbb Z_2$, where the action of $\mathbb R \subset \mathbb R\rtimes \mathbb Z_2$ is free. 

In fact,  $\underrightarrow {\tilde U}_{ 0}=Y_0\times\mathbb R$ and 
$$\underrightarrow {\tilde g}_{ 0}=dy^2 +\tilde \phi(y)^2 dt^2$$
where $\tilde \phi$ is a smooth function on $Y_0$, which, probably,\footnote{This is, of course,  easy to check but we do not need it here. In any case, we prefer the Haussdorf  limit definition of  $\underrightarrow {\tilde g}_{ 0}$, since it  is less demanding on the regularities of $X$ and $Y_0$.}  is equal to $\phi$ from section 2, which was  defined there via the second differential (variation)  of   the function $Y\mapsto vol_{n-1}(Y)$.

In general, the spaces  $(\tilde U_{[-\varepsilon,\varepsilon]},  \tilde g_\varepsilon)$ Hausdorff converge away from the $\Gamma_\varepsilon$-orbits of the  $\delta$-neighbourhoods of the singularity of $E(Y_0)\subset  \tilde U_{[-\varepsilon,\varepsilon]} $,
  where $\delta=\delta_\varepsilon\to 0$ for $ \varepsilon\to 0$.

Then $( \underrightarrow {\tilde U}_{ 0},\underrightarrow {\tilde g}_{ 0})$ stands for the metric  completion of the resulting (smooth Riemannian)  limit  space 
 for $\varepsilon, \delta\to 0$.\vspace {1mm}

{\it Question. }Does  $\mathbb R$ act  {\it freely} on  $\underrightarrow {\tilde U}_{ 0}$  in the case where    the minimal hypersurface  $Y_0\subset X$ has {\it singularities}?\vspace {1mm}

{\large \sf  Consecutive Torical Symmetrization of  $(n-m)$-Overtorical Manifolds.} A  compact oriented   $n$-dimensional Riemannian manifold $X$ with (possibly empty) boundary is called 
 {\it $(n-m)$-overtorical}, if it comes along with    a continuous map   $f:X\to \mathbb T^{n-m}$,
 such that the   "homological pullback" of a point $\theta_0\in \mathbb T^{n-m}$, denoted $ f^\ast[\theta_0]\in H_m(X,\partial X)$,
 doesn't vanish.\vspace {1mm}

(To clarify, recall  that if  $f_o$ is a  smooth map homotopic to $f$, then the pullbacks $ f_o^{-1}(\theta)\subset X$ of {\it generic} $\theta\in \mathbb T^{n-m}$ are smooth oriented submanifolds of dimensions $m$ which are homologous to $ f^*[\theta_0]$, i.e.  $[f_o^{-1}(\theta)]= f^*[\theta_0]$. 

Alternatively, $ f^*[\theta_0]$ can be defined as the Poincare dual of the image of the fundamental cohomology class $[\mathbb T^{n-m}]^\ast \in H^{n-m}(\mathbb T^{n-m};\mathbb Z)$
 under the cohomology homomorphism $f^\ast:  H^{n-m}(\mathbb T^{n-m};\mathbb Z) \to H^{n-m}(X;\mathbb Z)$.
 
 For instance, $\mathbb T^{n-m} \times B$ is  $(n-m)$-overtorical  for all compact $m$-dimensional manifolds $B$, possibly  with a boundary.)
 \vspace {1mm}
 

{\sf If  $X$ is  $(n-m)$-overtorical and if  $Y_0\subset X$ represents the homology class in $H_{n-1}(X,\partial X)$  corresponding    to the (homological)  $f$-pullback of a codimension $1$ torus in $\mathbb T^{n-m}$, 
then the manifold $X(\varepsilon) =\tilde U_{{[-\varepsilon,\varepsilon]}}/2\varepsilon \mathbb Z$ is also 
{\it $(n-m)$-overtorical} for all sufficiently small $\varepsilon >0$.}\vspace {1mm}

In fact, this has nothing to do with minimality of $Y_0$. It is, obviously,  true for all hypersurfaces $Y_0\subset X$   in suitable homology classes in  $H_{n-1}(X,\partial X)$ and their small neighbourhoods $U_{[\pm \varepsilon]}\subset X$  say with smooth boundaries or even  with singular boundaries, provided the singular loci have codimensions $\geq 2$.\vspace {1mm}

Now let us iterate
$$X\mapsto X(\varepsilon_1) \mapsto  X(\varepsilon_1,\varepsilon_2)=X(\varepsilon_1)(\varepsilon_2)\mapsto  X(\varepsilon_1,\varepsilon_2,..., \varepsilon _i)\mapsto...,$$
where each step  
$$X(\varepsilon_1,\varepsilon_2,..., \varepsilon _i) \mapsto X(\varepsilon_1,\varepsilon_2,..., \varepsilon _i, \varepsilon _{i+1})$$
depends on the choice of a codimension 1 subtorus in the corresponding torus $\mathbb T^{n-m}$,
where each minimal 
$$Y_{0_i} \subset X(\varepsilon_1,..., \varepsilon _i)\mbox { and the corresponding }
Y_{\pm\varepsilon_i}\subset X(\varepsilon_1,..., \varepsilon _i)$$
are taken in the non-singular loci of $X(\varepsilon_1,..., \varepsilon _i)$ modulo their complements (which have positive codimensions.

It should be noted that the $(n-m)$-tori serving   $X(\varepsilon_1,..., \varepsilon _i)$ and  $X(\varepsilon_1,..., \varepsilon _j)$ with $j\neq i$  are {\it not canonically isomorphic.} 
However, it make sense of  taking a  generic infinite  sequence of these subtori, and then to send all $\varepsilon_i\to 0$.

 \vspace {2mm}

 Let us incorporate the relevant  properties of the resulting limit space, call it   $\tilde X_\infty$, in the following.
  \vspace {1mm}

{\sf \large  Definition of    $\mathbb R^{n-m}\rtimes O(n-m)$-Symmetrization.}    Given an $(n-m)$-overtorical  manifold $X$, call an oriented manifold  $\tilde X_{sm}$ with (possibly empty) boundary   an {\it   $\mathbb R^{n-m}\rtimes O(n-m)$-symmetrization of $X$} if it satisfies the following nine conditions.

  
  \vspace {1mm}

$\bullet_1$ $\tilde X_{sm}$ is   {\it isometrically acted upon  by  $\mathbb R^{n-m}\rtimes O(m)$}, that is  the isometry group of $\mathbb R^{n-m}$, such that {\sl the orbits of this  action are equal to the $\mathbb R^{n-m}$-orbits}.   \vspace {1mm}
 
$\bullet_2$  The action of  $\mathbb R^{n-m}$ on $\tilde X_{sm}$ is {\it free.}\vspace {1mm}
 
$\bullet_3$  The quotient map $\tilde X_{sm} \to \tilde X_{sm}/\mathbb R^{n-m}$ admits an {\it inverse},  say  $E_\infty :\tilde X_{sm}/\mathbb R^{n-m}\to \tilde X_{sm}$, where the image 
$$E_\infty (\tilde X_{sm}/\mathbb R^{n-m})\subset\tilde X_{sm}$$ 
is {\it normal} to the $\mathbb R^{n-m}$-orbits  in  $\tilde X_{sm}$.\vspace {1mm}

\vspace{1mm}
 
 $\bullet_4$  there exists  
 {\it a continuous map}   
  $$\tilde\varphi_\infty :\tilde X_{sm} \to X,$$
with the following properties.\vspace{1mm}

$\bullet_5$ The map $\tilde\varphi_\infty$ is {\it $\mathbb R^{n-m}\rtimes O(n-m)$-invariant}.\vspace{1mm}

$\bullet_6$   The corresponding  map  $\varphi_\infty$ from the quotient  manifold $ \tilde X_{sm}/\mathbb R^{n-m}$ to $X$,  
$$\varphi_\infty=\tilde \varphi_\infty/\mathbb R^{n-m}:   \tilde X_{sm}/\mathbb R^{n-m} \to X,$$
{\it sends the fundamental homology class 
$$[ \tilde X_{sm}/\mathbb R^{n-m}]\in H^{m}(\tilde X_{sm}/\mathbb R^{n-m}, \partial \tilde X_{sm}/\mathbb R^{n-m})\mbox { to } 
f^\ast[\theta_0]\in H^m(X, \partial X).$$}
 

$\bullet_7$ The map $\tilde\varphi_\infty$ is {\it 1-Lipschitz.}\vspace{1mm}

$\bullet_8 $ The map $\tilde \varphi_\infty$ is {\it scalar curvature non-increasing},
$$Sc(X)(\tilde \varphi_\infty (\tilde x))\leq Sc(\tilde X_{sm})(\tilde x), \hspace {1mm}\tilde x \in \tilde X_{sm},$$
$\bullet_9$ The map  $\tilde \varphi_\infty$ sends $\partial \tilde X_{sm} \to \partial X$ and it is  {\it  mean curvature non-increasing},
$$mn.curv(\partial X)(\tilde \varphi_\infty (\tilde x))\leq mn.curv(\tilde \partial X_{sm})(\tilde x), \hspace {1mm}\tilde x \in \tilde \partial X_\infty.$$

In short,  $\bullet_6$ -$\bullet_9$ say that  symmetrization $X\mapsto  \tilde X_\infty$ must be \vspace {1mm}

\hspace {20mm}{\sl the topology and distance in $X$ {\it non-increasing}}

 \hspace {-6mm}and, at the same time,  
 
\hspace {20mm} {\sl the scalar and mean curvatures {\it non-decreasing.}}
\vspace {2mm}

{\sf \large Symmetrization Theorem.} {\sl Every $(n-m)$-overtorical  manifold $X$ of dimension $n\leq 7$ admits an  $\mathbb R^{n-m}\rtimes O(n-m)$-symmetrization.} \vspace  {1mm}

{\it Proof.}   If $n\leq 7$, then  there is no serious  problem with singularities,  and  the above limit space $\tilde X_\infty$ satisfies the above conditions $\bullet_1$ - $\bullet_9$.

In fact, even though the natural metrics  on  $X(\varepsilon_1,..., \varepsilon _i)$ are only piecewise smooth, the limit space  $\tilde X_{sm}$ of the $\mathbb Z^{n-m}$-covers of   $X(\varepsilon_1,..., \varepsilon _i)$ for $\varepsilon_i\to 0$, $i=1,2,...$, is  a {\it smooth} Riemannian manifold with boundary.

Then   everything becomes fairly obvious, except for $\bullet_8$ and $\bullet_9$.

To prove these what we use is   that, at every symmetrization step $X\mapsto X(\varepsilon)$,
the Riemannian metric  on $X(\varepsilon)$, which is, a priori, only continuous and which may have a corner at the boundary, can be\vspace {1mm}

 {\sl $C^\infty$-smoothed with an  arbitrarily small decrease of the  scalar curvature of $X$
 as well as of the mean curvature of $\partial X$}.\vspace {1mm}
 
 \hspace {-6mm}  This is  explained in section 11.2.

Then the proof follows by   {\sl semicontinuity} of the scalar curvature and mean curvatures (the latter is essentially obvious)   under $C^0$-limits of Riemannian metrics, see [Gr 2014*] and [Bam 2016].

\vspace {2mm}

{\large \sf About $n> 7$.} It is not impossible that the Symmetrization Theorem remains valid for all $n$ but 
$\tilde X_\infty$  can be  used  for this purpose  only for $m\leq 6$.

For instance, suppose $X$ is homeomorphic to $\mathbb T^1\times S^7$. If  our minimizing hypersurface $Y\subset X$  
in the homology class of $S^7\subset X$ has a singularity, then, clearly,   the space $\tilde X_\infty$
is  also singular  because  $\tilde X_\infty/\mathbb R^1=Y$.

On the other hand if $X$ is homeomorphic to, say,   $\mathbb T^{n-m}\times S^{m}$ for $m\leq 6$, 
then  $\tilde X_\infty$  and $\tilde X_\infty/\mathbb R^{n-m}$ can be, a priori, non-singular.

This makes  the Symmetrization Theorem   plausible for all $n$ and $m\leq 6$.

And  the techniques (results?) from [SY 2017], probably, yield  the symmetrization theorem for all $n$ and $m\leq 2$.
\vspace{1mm}

{\it Almost  Symmetrization.} This means that the conditions $\bullet_8$ and $\bullet_9$ need to be satisfied up to an arbitrarily small error $\epsilon>0$:

$$Sc(X)(\tilde \varphi_\infty (\tilde x))\leq Sc(\tilde X_{sm})(\tilde x)+\epsilon$$
and 
$$mn.curv(\partial X)(\tilde \varphi_\infty (\tilde x))\leq mn.curv(\tilde \partial X_\infty)(\tilde x)+\epsilon$$

Now, since by [NS 1993] minimal hypersurface  in generic $8$-manifolds are non-singular, 
the argument used for $n\leq 7$,  implies the following.\vspace {1mm}

{\large \sf Almost Symmetrization Theorem for $n=8$.} {\sl All $(n-m)$-overtorical manifods of dimension $n\leq 8$ admit almost  $\mathbb R^{n-m}\rtimes O(n-m)$-symmetrizations.}

\section {Application of Symmetrization to   Manifolds with Positive and with Negative Scalar Curvatures.}

Let $V$ be a Riemannian band  and let   $Z_0\subset V$ be a closed hypersurface which separates $\partial_- V$ from $\partial_+ V$ and, thus, divides $V$ into two halves
$V_\pm \supset \partial_\pm(V)$. Let

$\bullet$   the mean curvatures  of $\partial_\pm V$  are bounded from below by some constants  $M_\pm$;  

$\bullet$ the scalar curvature of $V$ is bounded by  a given    function $\sigma=\sigma(d)$ of the signed distance $d=d(v)$  from $v$ to  $Z_0$, that is 
$$Sc(V)(v)\geq \sigma(dist_\pm(v, Z_0)),$$
where 
$$\mbox {$dist_\pm(v, Z_0)=dist (v, Z_0)$ for $v\in  V_+$ and $dist_\pm(v, Z_0)=-dist (v, Z_0)$  for $v\in  V_-$}. $$

\vspace {1mm} 
Since   $\mathbb R^{n-1} \rtimes O(n-1)$-symmetric metrics on  $\mathbb R^{n-1}\times [-l,l]$
cane be written as 
$$\hat g= \hat \varphi^2g_{Eu} +dt^2$$
where $g_{Eu}$ is the flat Euclidean metric on $\mathbb R^{n-1}$ 
and  where the scalar curvature of $\hat g$, which depends only on $t\in [-l,l]$, satisfies
 $$Sc(\hat g) =-2(n-1)\frac { \hat{\varphi}''}{ \hat\varphi}-(n-1)(n-2)\frac { (\hat{\varphi}')^2}{ \hat{\varphi}^2},\leqno {\RIGHTcircle}$$
the net  effect of $\mathbb R^{n-1} \rtimes O(n-1)$-symmetrization of  $V$ can be stated in concrete terms as follows. \vspace {1mm} 

 {\large \sf  Symmetrization of Bands.}  {\sl If the above band $V$ is  $\mathbb R^{n-1} \rtimes O(n-1)$-symmetrisable, then there exists  
 a smooth function  $\hat \varphi(t)=\hat \varphi_\sigma(t)$, on the segment $[-l,+l]$  
such that 
 $$\pm l\geq dist(Z, \partial_\pm V),$$
$$ \frac { \hat {\varphi}'(- l)}{ \hat \varphi(l)} \leq \frac {-M_-}{n-1}, \hspace {1mm}    \frac { \hat {\varphi}'( l)}{ \hat \varphi(l)} \geq  \frac {M_+}{n-1} $$
and 

$$-2(n-1)\frac { \hat{\varphi}''(t)}{ \hat\varphi}-(n-1)(n-2)\frac { \hat{\varphi}'(t)^2}{ \hat{\varphi}(t)^2}\geq \sigma(t),  $$
that is
$$-2f'(t) -nf(t)^2\geq\frac{\sigma(t)}{n-1}\mbox {   for } f(t)=\frac {\hat {\varphi}'(t)}{\hat\varphi(t)} \mbox { and all } t\in [-l,l].\leqno \hat{\RIGHTcircle_l}$$} 
 
\vspace {1mm}

\textbf {Symmetrization Corollary for $Sc\geq 0.$} Let $V$ be an isoenlargeable band 
and $Z_0\subset V$ be  a hypersurface which separates $\partial_-(V)$ from $\partial_+(V)$.
Let  $Sc(V)\geq 0$ and let 
{\it $$\mbox { $Sc(V)\geq \sigma_0>0$ on the $\delta_0$-neighbourhood of $Z_0$}.$$}
Then 

{\it the distance from $Z$ to the boundary $\partial V$  is bounded by a constant which depends only on the dimension of $V$, on $\sigma_0>0$ and on  $\delta_0>0$, 
$$ dist(Z, \partial V)\leq C=C_n(\sigma_0, \delta_0).$$}

Moreover,  this remains true 

{\it if the inequality $Sc(V)\geq 0$ is replaced by  $Sc(V)\geq -\varepsilon$  for a small positive $\varepsilon \leq \varepsilon_n(\sigma_0, \delta_0)>0$.} \vspace {1mm}

{\it Proof.} If   $ \sigma (t)\geq \sigma_0$ 
for $t \in [-\delta_0,\delta_0] \subset [-l,l]$ and $\sigma (t)\geq -\varepsilon$ for  all $t \in[-l,l]$, where $\varepsilon<< \delta_0, \sigma_0$,
then the inequality $\hat{\RIGHTcircle_l}$ implies that 
$$l\leq C= C_n(\sigma_0, \delta_0)$$
and the proof follows.

Notice that no  condition 
$mean.curv(\partial_\pm(V)\geq M_\pm$ has been  used at this point.

\vspace {1mm}
 
 \textbf {Sub-corollary for Complete Manifolds with $Sc\geq 0$.} {\it Open  isoenlargeable bands carry no  complete metrics with  scalar curvatures $Sc>0$.}\vspace {1mm}
 
 Moreover, \vspace {1mm}
 
 \hspace {5mm} {\sl  complete metrics with $Sc\geq 0$ on such bands  are Riemannian flat.}\vspace {1mm}
 
  In fact,   a deformation theorem by Kazdan and Warner  together with    the Cheeger-Gromoll splitting theorem  imply that if such a band admits no complete  metric with $Sc>0$ then every
complete  metric with $Sc\geq 0$ is  flat. \vspace {1mm}

\textbf {Representative Example.} {\sl If  a compact manifold  $Z$ admits a metric with negative sectional curvature,
then there is no complete  metrics with $Sc\geq 0$  on the  connected sums 
$$X= (Z\times \mathbb R)\#_i P_i$$  for  compact manifolds $P_i$.}

(A similar result is proven in 6.12 and 6.13 of  [GL 1983] for spin manifolds $X$.)

\vspace {2mm}

\textbf {Symmetrization Corollary for $Sc\geq \sigma<0.$} To get a perspective look at the following

{\it Model Example.} Let $V_{[-l,l] }$ be the band of width $2l$   between concentric horospheres in the hyperbolic space $H^n$ of constant curvature $-1$, which is the product 
$${V_{[-l,l] }}=\mathbb R^{n-1}\times [-l,l]\subset \mathbb R^{n-1}\times(-\infty +\infty)=H^n,$$
where the hyperbolic metric in these coordinates is   $g_{hyp}=e^{2t} g_{Eu}+dt^2$.
 
The scalar curvature of    $g_{hyp}$ in these coordinates,   in agreement with $\RIGHTcircle$,  is  $-n(n-1)$, while the mean curvatures of the boundaries $\partial_\pm V_{{[-l,l] }}= \mathbb R^{n-1}\times \{\pm l\}$ are 
$$mean.curv (\partial_\pm V_{[-l,l] })=\pm (n-1).$$
 
Such a band  becomes  compact   if divided  by  the action of $\mathbb Z^{n-1}\subset \mathbb R^{n-1}$ for 
 $$\mathbb R^{n-1}\times [a,b]/\mathbb Z^{n-1}=\mathbb T^{n-1}\times [a,b].$$

Now,

{\sl let $V$ be a  compact Riemannian  overtorical band,   where the scalar curvature   and the mean curvatures of the boundaries satisfy
$$Sc(V)\geq -n(n-1),$$
$$ mean.curv (\partial_-(V))\geq -(n-1),$$ 
$$ mean.curv (\partial_+(V))\geq (n-1).$$
then, in fact,
$$Sc(V)= -n(n-1),$$
$$ mean,curv (\partial_-(V))= -(n-1),$$ 
$$ mean.curv (\partial_+(V))= (n-1).$$}

{\it Proof.}  If either of the above three inequalities is {\it strict }  (i.e. ">") at some point, then,  by slightly  conformally perturbing the metric of $V$ (see section 11.2),   one can make {\it all three strict at all points}. 

This, by symmetrization,   would result in   a function $ \hat f=\frac{\varphi'}{\varphi} $ on some segment $[-l,l]$, $l>0$, such that
$$f(-l)< 1, \hspace{1mm}  f(l) > 1  $$
and 
$$ 2f' < -n(f^2-1)$$
which is, obviously, impossible. QED.

\vspace {1mm}

\textbf {Sub-corollary: Weak Rigidity  of  $H^n/\mathbb Z^{n-1}$.}   Let 
 $X=H^n/\mathbb Z^{n-1}$, where $H^n$ is the hyperbolic space  with  the sectional curvature $\kappa(g_{hyp})=-1$, and the group  $\mathbb Z^{n-1}$ discretely and  isometrically  acts on  $H^n$ by parabolic transformations, i.e. preserving a  horosphere in $H^n$.\footnote
 {If $n\geq 3$ then all discrete isometric actions of $\mathbb Z^{n-1}$ on $H^n$ are parabolic.}

 {\sl  If a Riemannian metric on $X$, which    coincides with the hyperbolic one  (descended from $H^n$ to $X$) outside a compact subset in $X$, satisfies 
$$Sc(g)\geq Sc(H^n)= -n(n-1)$$
 then 
$$Sc(g)= -n(n-1)$$
everywhere on $X$.}

\vspace {1mm}

{\large \sf Soap  Bubbles   and Rigidity of Bands.}   A sharper version of the above  sub-corollary, namely  the implication
$$Sc(g)\geq -n(n-1)\Rightarrow \kappa(g)=-1\leqno \mbox {{[$\varocircle$$_{-1}$]}  } $$
follows from the existence of {\it stable minimal bubbles}  in $X$, which are closed  hypersurfaces $Y$ which separate the two ends in $X$ and 
which minimise the functional 
$$ Y\mapsto vol_{n-1}(Y)- (n-1)vol (X_{\prec Y}),$$
where $X_{\prec Y}\subset X$ is the  part of $X$ which is  bounded by $Y$ and  which has $vol<\infty$. \vspace {1mm}

{\it Proof. } The relation   [$\varocircle$$_{-1}$]  for $n\leq 8 $ follows   from   $\varocircle$$[M]$ below by taking $M=n-1$.
\vspace {1mm}

$\varocircle$${[M]}$.  Let  a compact Riemannian band $V$ 
admit a function $M=M(v)$ such that 

$$\mbox { $M_{|\partial _-V}\leq -mean.curv (\partial _-V)$ and $M_{|\partial _+V}\geq mean.curv (\partial _+V)$}$$ 
and 
$$ \frac {n}{n-1}M^2 -2||dM|| +Sc(V)\geq 0.$$

If $n=dim(V)\leq 8$, then either 

  {\it there is a closed hypersurface $Y_\circ\subset V$, which separates $\partial _-V$
from $\partial _+V$ and which admits a metric with $Sc>0$},

 or 
 
  {\it $V$ decomposes into the (warped) product, $V= Y\times[-l,l]$ with the metric
 $$\varphi(t)^2g_Y +dt^2,$$
  where  the Riemannian metric $g_Y$ on $Y$ has zero Ricci curvature.}
\vspace {1mm}

This is shown for $n\leq 7$ in section $5\frac {5}{6}$ in [Gr 1996])   and  the case $n=8$ 
 can be taken care of with a help of    ideas  from [NS 1993].

But it  seems that the regularisation techniques of [Loh 2016] and/or of [SY 2017]
do not apply, at least not directly, to this case  and the validity of the above statement for $n\geq 9$ remains quite  problematic.

On the other hand [$\varocircle$$_{-1}$], where the implied $Y_\circ$ is the $n-1$ torus, must follow from these techniques  which are, in principle, applicable  whenever  torical symmetrization works.

\vspace {1mm}

{\large \sf On Min-Oo Rigidity Theorem.} By adapting an idea of  Witten to a "hyperbolically modified" Dirac operator, Min-Oo [Mi 1989] proved a version of the 
positive mass theorem for  $H^n$. In particular he has shown the following.\vspace {1mm}

{\sf \textbf { [MRT]}}  {\it If a complete spin manifold $X$ is isometric to $H^n$ outside a compact subset and if
$Sc(X)\geq -n(n-1)$ then $X$ is isometric to $H^n$.}
\vspace {1mm}

Since compact perturbations of $H^n$ can be  periodically extended by discrete actions of isometry groups $\Gamma$ on $H^n$, e.g. for the above parabolic  $\mathbb Z^n$, 

\vspace {1mm}

\hspace {33mm} {\sl  {\sf \textbf { [MRT]}}  follows from} {[$\varocircle$$_{-1}$]}. \vspace {1mm}

Thus, {\sf \textbf { [MRT]}} {\sl remains valid  without assuming $X$ is  spin}  (but with some reservations for $n\geq 9$).
\vspace {1mm} 

Moreover, this is shown in [AndMinGal 2007], that [MRT] combined with an argument from [Loh 1999] implies the positive mass theorem for $H^n$ and the spin condition   can be    disposed of  in the context of the full Min-Oo(-Wang Chru\'sciel-Herzlich) theorem  (unconditionally for $n\leq 8$).

 \vspace {1mm}

{\Large \sf  Proving  Rigidity by Symmetrization.}  The rigidity  of bands $V$ in the symmetrization context says that the universal coverings of

{\sl  the  extremal bands, where our $\frac {2\pi}{n}$-inequality becomes an equality,  must be $\mathbb R^{n-1} \rtimes O(n-1)$-invariant.}

We shall indicate below  the proof of this for  $n\leq 7$, where   the dimension $n=8$ needs a little  effort, and where   the  regularisation as developed in  [Loh  2016] and in [SY 2017] for $n\geq 9$ may  need an  additional  refinement to yield rigidity.

Now, assuming minimal varieties are non-singular, we observe that \vspace{1mm}

{\sf  the symmetrization process {\it strictly enlarges} the scalar curvature of  $V$,  unless the minimal 
hypersurfaces $Y\subset V$  used for this process are  {\it totally geodesic}.}\vspace{1mm}

In fact, by  the second variation formula in the form given to it in [SY - Inc 1979], 
 the corresponding operator $L$ from  section 2 is  {\it strictly positive}, which  implies  increase of the scalar curvature under symmetrization.
And this also work for symmetrization by reflection  in section 2 if one replaces the smoothing of edges argument by  an appeal to the  corresponding operator $L$.
 
 Thus, one represent all our  homology classes
 in $H_{n-1}(V, \partial V)$ by   totally geodesic submanifolds. This strongly restricts the geometry of $V$ but   does not, at least not obviously, imply the required 
  $\mathbb R^{n-1} \rtimes O(n-1)$-symmetry of $V$.
  
  However, by applying the same argument to the  soap bubbles  $Y_{\pm \varepsilon}\subset X$  which lie 
close to minimal $Y$ and minimise the functional 
$Y\to vol_{n-1}(Y)-\varepsilon \cdot vol_n U_{\pm\varepsilon}$ as in section 8
one  sees that no minimal $Y$ can be  locally strictly minimising  in either of the two halves it divides $V$ into.

This shows, that minimal $Y$ in all  homology classes, besides being totally geodesic,   are "freely movable" in $V$, namely, they serve as fibers of a fibrations of  $V$ over the circle. 

Then  the required $\mathbb R^{n-1} \rtimes O(n-1)$-symmetry of $V$ easily follows.\footnote{Since I have not written this down in detail, I might have missed some hidden difficulty in  this apparently quite  innocuous argument.}

\section {Comparison with Results Obtained with Twisted  Dirac Operators.}

Besides the method of   of  minimal hypersurfaces, a non-trivial information on geometry (and topology)  of Riemannian manifolds $X$ with $Sc (X) \geq \sigma$, $\sigma\in (-\infty, \infty)$, can be obtained by confronting

I:  Atiyah-Singer type  {\it index  theorems  for     Dirac operators}  which  yield {\it non-zero harmonic spinors} on $X$  
 
 with 
 
 II: the  {\it twisted    Schroedinger-Lichnerowicz-Weitzenboeck  formula} for manifolds with {\it lower bounds on their  scalar curvatures} which rules out, or significantly  restricts, such spinors.

 Comparison of  (partly overlapping)  results obtainable with minimal hypersurfaces and with Dirac operators exposes     {\it limitations of both  methods} and   exhibits  wide  gaps in our understanding of scalar curvatures;   this  begs for  a  new approach.\vspace {1mm}

Let us briefly demonstrate  this on   a few simple examples.\footnote{See [Ros 2007], [BER 2014],       [Gr2017] for  more elaborated techniques and  examples.} \vspace {1mm}

\textbf {(1)} {\large \sf Spin, Spinor Bundles and  Dirac Operators.}  Since the fundamental group of the special (i.e. orientation preserving) orthogonal group $SO(n)$  for $n\geq 3$ is $\mathbb Z/2\mathbb Z$,  there are  exactly {\it two} different orientable  bundles of rank $n\geq 3$  over   closed connected  surfaces. The trivial bundle  is, by definition, {\it spin} and the non-trivial one is {\it non-spin}.
 
 An orientable manifold $X$ is called {\it spin} if the restrictions of the tangent bundle $T(X)$ to all surfaces $S\subset X$ are spin (i.e trivial).\footnote{"Spin"  makes sense also for non-orientable bundles and manifolds but we do not need them at this point.}

For instance, all orientable hypersurfaces  $X^n \subset \mathbb  R^{n+1}$  are spin, all $3$-manifolds are spin  and 

{\it simply connected n-manifolds with trivial second homotopy groups are spin.}

The simplest non-spin manifolds are the complex projective spaces $\mathbb CP^n$ of even complex dimensions $n$ and connected sums of other manifolds with these  $\mathbb CP^n$.
 
 \vspace {1mm}

The {\it spinor bundle}  of a  Riemannian spin  manifold $X$ of dimension $n$, denoted   $ \mathbb S(X)$, is a unitary vector bundle of  vector bundle of rank  $2^n$ with a unitary connection associated  to the Levi-Civita connection in $T(X)$. If $n$ is even, the bundle  $ \mathbb S(X)$
splits, $\mathbb S=\mathbb S^+\oplus \mathbb S^-$.

{\it The Dirac operator}  is a canonically defined first order differential  operator  $D$  represented as a certain "natural"   linear combination of covariant derivatives
which  act on $ \mathbb S(X)$.  (See [BHMM 2015] for definitions,  basic results and geometric applications of the Dirac operator.)

When $n$ is even the Dirac operator splits:

  $D=D^+\oplus D^-$, where the operators  $D^+$ and  $D^-$
  are mutually adjoint for $D^+:C^\infty( \mathbb S^+) \to C^\infty( \mathbb S^-)$  and   $D^-:C^\infty( \mathbb S^-) \to C^\infty( \mathbb S^+)$ and  where $ind(D)=dim(ker D^+)-dim (kerD^-)$.
The solutions of $D(s)=0$ are called harmonic spinors on $X$.

\vspace {1mm}

{\it Twisted  Dirac operator.}  Given a complex vector bundle $L=(L,\nabla)$ with a linear connection, one naturally defines 
$$D^\pm_{\otimes L} : C^\infty(\mathbb S^\pm 
\otimes L) \to  C^\infty(\mathbb S^\mp \otimes L), $$
where the sections of  $S^\pm\otimes L$ in the kernel of $D_{\otimes L}=D^-_{\otimes L}\oplus D^+_{\otimes L}$ are called {\it $L$-twisted harmonic spinors.}

\vspace {1mm} 
 
 \textbf {(2)  \large \sf Chern character and Todd Genus.} The Chern character of a complex vector bundle $L$ over $X$  is a  certain  polynomial in the Chern   classes $c_i\in H^i(X;\mathbb Z) $  of $L$ in the rational  cohomology $H^\ast(X;\mathbb Q)$  starting from $c_0=rank(L)$,  
 $$ch(L)= c_0+ c_1+\frac {1}{2}(c_1^2+c_2 ) +\frac {1}{6} {c_1^3 +c_1c_2+3c_3}+...+ \frac {1}{i!} (c_1^i+...+k_ic_i)+...$$
where, observe, all $k_i\neq 0$. 
The basic  properties of $ch$ (which essentially define it) are  additivity and multiplicativity:
$$ ch(L_1\oplus  L_2)= ch(L_1) + ch (L_2) \mbox { and } ch(L_1\otimes L_2)=ch(L_1)\smile ch(L_2).$$

 {\it The  $\hat A$-genus} is another  polynomial, now in   the Pontryagin classes $p_i=p_i(T(X))\in H^\ast(X;\mathbb Q)$,
 $$\hat A(X)=1-\frac {1}{24}p_1+\frac {1}{5760}(-4p_2+ 7p_1^2)+  \frac {1}{967680}( -16p_3+44p_2p_1 -31p_1^3)+...$$
where again the coefficients at $p_i\in H^{4i}(X;\mathbb Z)$ are non-zero. 
 \vspace {1mm}

 \textbf {(3)  \large \sf Topological  Index  $I$.} 
   Let $\hat X$ be an oriented   Riemannian spin  $\Gamma$-manifold,  which means  $\hat X$
   is acted upon by a group $\Gamma$ and let 
    let $\hat L_1= (\hat L_1,\nabla_1)$ and $ \hat L_2 =(\hat L_2,\nabla_2)$ be complex   vector bundles with unitary connections such that $\Gamma$  acts on $\hat L_1$ and on $\hat L_2$ by fiber-wise unitary  (linear isometric) connection preserving   transformations compatible with the action  of $\Gamma$ on $\hat X$, such that
     the following conditions are satisfied.
 
(i)  The   action of $\Gamma$ on $\hat X$ is  { \it proper,  isometric  and orientation preserving}, where "proper" means that  there are at most finitely many $\gamma\in \Gamma$, such that  for all compact subsets $K\subset \hat X$ the intersections $K\cap \gamma(K)$  are empty for all but finitely many $\gamma\in \Gamma$.

preserves the connections in these bundles.

(ii) There exists  a unitary connection preserving  $\Gamma$-equivariant isomorphism between the bundles $\hat L_1$ and $\hat L_2$ at infinity, that is  on the complement to a $  \Gamma$-invariant subset $V \subset \hat X$ such that $V/\Gamma$ is compact.
 \vspace {1mm}
 
Let  
$$I=I(\hat X, \hat L_1\ominus \hat L_2) = (\underline{ \hat A\smile (ch(\hat L_1)-ch (\hat L_2))}[\hat X/\Gamma]$$
be  defined   by representing   $\hat A$,  $ch(\hat L_1)$  and $ch (\hat L_2)$ by the   Chern-Weil differential forms on $\hat X$, call them $\alpha, \lambda_1, \lambda_2 \in \bigwedge^\ast(\hat X)$ which, clearly, are $\Gamma$-invariant  and where  $\lambda_1-\lambda_2$ vanishes outside $V$.

Since   the form   $\iota=\alpha \wedge (\lambda_1 -\lambda_2) $  vanishes outside $V$,   since it     is $\Gamma$-invariant and since  the action of $\Gamma$ on $\hat X$  is proper,  $\iota$  descends to a form $\underline \iota$ on  the quotient space  $\hat X/\Gamma$, which vanishes outside a compact subset\footnote {We do not assume the action of $\Gamma$  on $\hat X$  to be free  and  the space  $\hat X/\Gamma$  may be singular but our forms are  defined on it anyway.} and defines  the cohomology class $(\underline{ \hat A\smile (ch(\hat L_1)-ch (\hat L_2))}$ of 
$\hat X/\Gamma$ with compact supports,
 $$[\underline \iota] = (\underline{ \hat A\smile (ch(\hat L_1)-ch (\hat L_2))}  \in H_{comp}^n(\hat X/\Gamma;\mathbb R).$$

  Then the index $I=\underline \iota[\hat X/\Gamma]$ can be defined as the integral
  $$\int_{\hat X/\Gamma}\underline \iota = \int_{\Delta}\iota $$
  for a fundamental domain $\Delta\subset \hat X$.

   \vspace {1mm}

\textbf  {(4)  \large \sf  Atiyah's $L_2$-Index Theorem.} Let the following conditions be satisfied.

(a) The manifold $\hat X$ is {\it complete}.

(b) The connections in $\hat L_1$ and  in $\hat L_2$ are {\it unitarizable}. This means  these  bundles admit  unitary structures, i.e. fiberwise Hermitian scalar products $ \langle... \rangle$,  preserved by the parallel transport in these connections.

(c) The above  unitary structures, (which are unique up to scaling) are {\it $\Gamma$-invariant}.

(d) The operators $D^2_{\otimes \hat L_1}$ and $D^2_{\otimes \hat L_2}$ are {\it uniformly positive at infinity$/\Gamma$}, where a differential operator $\cal D$ on sections $s=s(\hat x)$ of a unitary bundle on a manifold $\hat X$  with a $\Gamma$ action is called 
uniformly positive at infinity$/\Gamma$, if 
$$ \int_{ \hat X}\langle {\cal D}(s(\hat x)), s(\hat x)\rangle_{\hat x}d\hat x \geq c\cdot \int_{\hat X}  ||s(\hat x)||^2$$
for a constant $c>0$ and all sections $s$  with compact supports  outside  a certain  subset  $V\subset \hat X$ such that $V/\Gamma$ is compact.

  \vspace {1mm}

{\it If the topological index $I= I=I(\hat X,  \hat L_1\ominus \hat  L_2)$ does not vanish, then there exists  either 
an $\hat L_1$- or $\hat L_2$-twisted harmonic   square integrable spinor on $\hat X$.}

In fact, 

{\it the von-Neumann dimensions of the kernels $\hat K_{1,2}^\pm$ of the operators 
 $D^\pm _{\otimes \hat  L_{1,2}}$ satisfy
$$ dim_\Gamma (\hat K^+_1) -dim_\Gamma (\hat K^-_1)-dim_\Gamma (\hat K^+_2)+dim_\Gamma (\hat K^-_2)=I(\hat X, \hat L_1\ominus \hat L_2).$$}

  \vspace {1mm}

{\it About the Proof.} The equality
 $$ dim_\Gamma (ker D^+_{\otimes \hat L}) - dim_\Gamma (ker D^-_{\otimes \hat L_2})=I(\hat X, \hat L)$$ in the case of {\it compact} $\hat X/\Gamma$  is proven in [At  1976].

(If $\hat X/\Gamma$ is compact  only a single bundle $\hat L= \hat L_1$ is needed, since one may take the trivial bundle of rank zero for $\hat L_2$; then    the   conditions (a) and (d) are irrelevant.)

The case of {\it non-compact manifolds} with no $\Gamma$-actions is treated in [GL 1983]. 

The compatibility of the two arguments was pointed out in [Gr 1986], where one finds further references.
  \vspace {1mm}

{\it Suggestion.}  It would be interesting to  remove or  to relax  some of the  conditions in the formulation of the index theorem. 
  \vspace {1mm}

{\large \sf  $ \widetilde {spin}$-Example.} Let $\hat X$ be the universal covering $\tilde X$  of a manifold   $X$. If $ X$  is   spin then  the spin bundle $\mathbb S(X) $ and the Dirac operator in it are defined and lift $\Gamma$-equivariantly to $\hat  X=\tilde X$ for the Galois action of $\Gamma=\pi_1( X)$ on $\tilde X$.

But if   $X$ is non-spin, yet $\tilde X$ is spin, then the group which acts on  
$\mathbb S(\tilde X)$ is the semidirect product  $\mathbb Z_2 \rtimes  \Gamma$ where $\mathbb Z_2$
acts by the  $\pm 1$-involution on  spinors which  corresponds to the Galois involutive transformation on the double covering of the principal bundle associated  to the tangent bundle  $T(\tilde X)$.
 
 Thus,  \vspace {1mm}
 
 {\it  Atiyah's $L_2$-index theorem applies to the Galois  coverings $\hat X$ of {\it non-spin manifolds} $X$ whenever these  $\hat  X$  are {\it spin}.} 
 
 \vspace {1mm}

\textbf  {(5) \large \sf  Twisted    Schroedinger-Lichnerowicz-Weitzenboeck  Formula.} This formula relates the squares of $L$-twisted Dirac operators with the  {\it rough Laplacians $\nabla^\ast \nabla$} in the   bundles 
$L=(L,\nabla)$ on $X$ with unitary connections, where, recall,  the operators   $\nabla^\ast \nabla$ acts on sections of $L$; they are     (non-strictly) positive 
$$ \int \nabla^\ast \nabla=\int ||\nabla||^2;$$
thus their kernels consist of $\nabla$-parallel sections of $L$  and  $rank(ker(\nabla^\ast \nabla)\leq rank (L)$.

Here is the formula.
$$D^2_{\otimes  l}=  \nabla^\ast \nabla+\frac {1}{4} Sc(x)\cdot Id +\cal R,$$
where $\cal R$ is a linear self adjoint  endomorphism (zero order operator) of $\mathbb S\times L$ defined by the {\it operator valued  curvature form $R$ of $L$} coupled by the {\it Clifford multiplication} in 
 $\mathbb S$ as follows.
$${\cal R}(s\otimes l)=\frac{1}{2} \sum_{1\leq i<j\leq n}e_ie_js\otimes  R(e_i\wedge e_j)(l),$$
where $Id : \mathbb S\times L$ is the identity operator,      where $e_i\in T_x(X)\subset  T(X)$ is an orthonormal frame   at the point $x\in X$,  where the above formula applies and where $s\in \mathbb S_x$ and $l\in L_x$.

Since  the Clifford multiplication operators   $e_i: s\mapsto e_is$ are unitary,  
$$|{|\cal R} (s\otimes l)l| \leq\frac{n(n-1)}{4} ||R||\cdot ||s||\cdot ||l||$$ 
where $||R||$ is the supremum of the norms of the curvature operator over all unit bivectors in the tangent spaces $T_x(X)$.

It  follows then the norm of the operator $\cal R$ is bounded by 
$$ ||{\cal R}||\leq const_n ||R||$$
for
$$const_n =  \frac{n(n-1)}{4} \sqrt {rank(\mathbb S)} = n(n-1) 2^{n-2} .$$

\textbf {(6) \large \sf  } Let $\hat X$ be a $\Gamma$-manifold with $\Gamma$-invariant  bundles $\hat L_{1,2}$, such that  the  assumptions (a), (b) and (c) in the above  Atiyah's $L_2$-index theorem are satisfied. 

Let, moreover, the  norms of the curvature operators $ R_1$ and $R_2$ of the (unitary) connections in $\hat L_1$ and $\hat L_2$   be bounded by 
$$ Sc(\hat X)(\hat x) \geq\varepsilon+  4const_n \cdot max ( ||R_{1}||_{ 
\hat x} , ||R_{2}||_{\hat x} )$$
for the above $const_n= n(n-1) 2^{n-2}$,  some $\varepsilon>0$ and  all $\hat x\in \hat  X\setminus V$ for a subset $V\subset \hat X $ with compact quotient $V/\Gamma$.\vspace {1mm}

Then the above \textbf {(4}) and \textbf {(5)} yield the following.\vspace {1mm}

 {\sf \textbf{ Theorem.}} {\it If the topological index $$I=I(\hat X,\hat  L_1\ominus \hat L_2) = \underline{ \hat A\smile (ch(\hat L_1)-ch (\hat L_2))}[X/\Gamma]$$
 doesn't  vanish,
 then there exists a point $\hat x\in \hat X$, where 
$$Sc(\hat X)(\hat x)\leq  4 const_n \cdot max ( ||R_{1}||_{\hat x} , ||R_{2}||_{\hat x} ).$$}
\vspace{1mm}

\textbf  {(7) \large \sf  Area  Enlargeable Manifolds } Recall that an $n$-dimensional  Riemannian manifold $X$ is called  {\it area enlargeable} if it admits a sequence of orientable  coverings $\tilde X_i\to X$ and 
of smooth maps $f_i:\tilde X_i\to S^n$ which are 

$\bullet_1$ constant at infinity,

$\bullet_2$ have non-zero degree,

$\bullet_3$ contract the areas of the surfaces $\Sigma\subset \tilde X_i$ by  
$$area(f_i(\Sigma))\leq \alpha _i area  (\Sigma) \mbox  { for } \alpha_i\underset {i\to \infty}\to 0$$

Observe that  area enlargeability is a weaker condition than enlargeability,  where instead of $\bullet_3$ one requires $Lip(f_i)\to 0$ (see section 4),  and that  area enlargeability, similarly to enlargeability,  is a {\it homotopy invariant} of {\it compact}  manifolds $X$.

Let us show of that area enlargeability is incompatible with $Sc>0$.

\vspace {1mm}

$[\square]$  {\it  Complete  area enlargeable  manifolds $ X$ the universal coverings of which are spin can't have 
$Sc(X)\geq \varepsilon>0$.} \vspace {1mm}

{\it Proof.} Let's first assume that  $n=2m$ and let $ L$ be a complex vector bundle of some rank $N$ over $S^n$ with {\it non-zero} Chern class $c_m\in H^n(S^n)$. 

Let $X$ be  the universal covering $\tilde {X}$ acted upon by $\Gamma=\pi_1( X)$, let $L_1$ be  the trivial bundle $X\times \mathbb C^N$     and let   $L_i$ be induced from $L$ by the composed map 
$$X=\tilde X\to \tilde { X_i}\tightoverset {f_i}\to S^n.$$
It is easy to see that non-vanishing of $c_m$ implies non-vanishing of  the topological index $I$ and that  the curvature of $L_i$ tends to zero for $i\to \infty$ 

Therefore,  the above \textbf {(6)}  applies to $(X, L_1,L_i) $ for a sufficiently large $i$ and yields  the proof for even $n$, while  the  case of  $n= 2m-1$ reduces to $n=2m$ by taking   $X\times S^1$. \footnote{ "Area enlargeable" appears as  "$\Lambda^2$-enlargeable"  in 
 [GL 1983],    where the coverings $\tilde {V_i}$ are assumed spin.}

\vspace {1mm}

\textbf {  (8) \large \sf  Llarull's Rigidity   Theorem.} The above, as it is shown in [Ll 1998] can be rendered sharp by taking the positive (or negative)  spin bundle $\mathbb S^+(S^n)$ for $\underline L$.

The Chern character of $\mathbb S^+(S^n)$ for $n=2m$ is equal to  the fundamental cohomology  class $[S^n]\in H^n(S^n)$ and the norm of the  Levi-Civita connection in this bundle  equals $\frac {1}{2}$ -- all this is more or less obvious.

What is less obvious (see {Ll 1998], [Min-Oo 2002] is that 
 the lowest  eigenvalue of the operator $\cal R$ on $\mathbb S\otimes \mathbb S^+$ on $S^n$ is equal $-\frac {n(n-1)}{4} $, which, by\textbf {(6)} (and a trifle  of linear algebra)   implies   the following\vspace {1mm}

$ \bigcirc$   Let $ X$ be a  Riemannian manifold, such that 

$\bullet $ $ X$ is complete,

$\bullet $   $ Sc( X)( x) \geq \varepsilon>0 $ for all $\underline x$  outside a compact subset in $X$. \vspace {1mm}

$\bullet $ the universal covering $\tilde X$ of $X$ is spin.

Let  a continuous  map   $f: X\to S^n$ satisfy the following conditions.\vspace {1mm}

$(\ast_{/\infty})$ $f$ is constant at infinity (i.e. constant outside a compact subset in $X$);

$(\ast_{deg})$ $f$ has non-zero degree;

$(\ast_{C^1})$ $f$ is $C^1$-smooth;

$(\ast_{ar})$ The map $f$ (non-strictly) decreases  integrals of the scalar curvature of $X$ over all smooth surfaces $\Sigma\subset X$.
(Since $S^n$ has constant scalar curvature $n(n-1)$ this amounts to the inequality
$$\int_SSc(X)(\sigma) d\sigma\geq n(n-1) area(f(\Sigma). )$$

Then\vspace {1mm}

{\it The map $f$ is a homothety: there exists a constant $\lambda>0$, such that 
$$dist_{S^n} (f( x_1), f( x_1))=\lambda\cdot  dist_{X} ( x_1,   x_1)\mbox  { for all } x_1, x_2 \in X.$$}

\vspace {1mm}

{\it About  the  Proof.} Here,  the Dirac operator on $ X$ is twisted with  the bundles $L=L_1$, which is   induced by $f: X\to S^n$  from $\mathbb S^+(S^n)$, and where one takes  the trivial bundle  of the same rank as $L$  for $L_2$.

 In this case, 
the formula for ${\cal R}: \mathbb S\otimes L\to \mathbb S\otimes L$ 
from the above \textbf {(5)}, that is 
$${\cal R}(s\otimes l)=\frac{1}{2} \sum_{1\leq i<j\leq n}e_ie_js\otimes  R(e_i\wedge e_j)(l),$$
written in the frames  of vectors $ e_i\in T_{\underline x}$, which simultaneously diagonalize the Riemannian metric of $ X $ and the metric induced by $f$ from $S^n$  effectively describes the action of $\cal R$ on the  corresponding  (Clifford) basis  in
 $\mathbb S( X) \otimes f^!(\mathbb S^+(S^n)$, which is  $\{\underline e_{i_1}\underline e_{i_2}... \underline e_{i_n}\otimes \underline e_{j_1} e_{j_2}... \underline e_{j_n}\}$.
Then a straightforward computation in  [Ll 1998] (and/or a  more conceptual argument in  [Min-Oo 2002]) shows that the spectrum of $\cal R$ is bounded from below by $-\frac {n(n-1)}{4}$ and the above \textbf {(6)} applies.

The above settles the case of even $n$. 

 If   $n$ is odd one uses  area contracting maps   $X\times S^1(R)\to S^{n+1}$  for large $R$ where the corresponding $\cal R$ is still bounded by   by
$\frac {1}{4}(n(n+1)=  \frac {1}{4}Sc(S^n)$    because the natural splitting of  metric in $X\times S^1\to S^{n+1}$  (see    [Ll 1998]).

Alternatively, one can  construct (non-split)   metrics $g_\varepsilon$ for on   $X\times S^1\to S^{n+1}$, for 
all  $\varepsilon>0$,  with $Sc(g_\varepsilon)\geq (n+1)(n+2) - \varepsilon  = Sc(S^{n+1})- \varepsilon$, such that area 
non-increasing maps $X\to S^n$ suspend to area non-increasing maps 
$(X\times S^1,g_\varepsilon) \to S^{n+1}$.

\vspace {1mm}

{\it Generalisation.} It is shown in  [GS 2002]  that the above remain valid for  $S^n$ if the standard metric $g$ on $S^n$ is replaced by  $g'$  with   {\it positive curvature operator}. This, shows, in particular, that Llarull's theorem is {\it stable  under  small perturbations} of the spherical metric $g_0$.

 \vspace {1mm}

\textbf { (9)  \large \sf   Discussion.} There are two  drawbacks of the above results compared to what can be done with minimal hypersurfaces.  

{ \sf I. Spin.} In the original paper [Ll 1998] the manifold $ X$   was assumed spin, which we have    relaxed to requiring the universal covering of $X$ to be  spin. Yet, we still can't prove, 
$\bigcirc$ or even $\square $ for all complete manifolds.

{ \sf II. Completeness.} Neither $\bigcirc$ or  $\square $ hold true as they stand for incomplete manifolds  and it is  unclear what  their correct  reformulations should be. 

And even if the area decreasing condition for maps $f: X\to S^n$  is  strengthened  to   to  
$Lip(f)\leq 1$, one can't   get any bound on $Sc(X)$ with Dirac operator methods  for incomplete  $ X $,   while  minimal hypersurface do allow such bounds (see section 3). 
\vspace {1mm}

On the other hand,  the Dirac operator results also have two advantages over those achieved with minimal hypersurfaces.

 {\sf [i] Area Versus Length.}  Application of  minimal hypersurfaces depends on distance rather than area estimates of metrics involved.    

 {\sf [ii] Non-Abelian Symmetries.} Dirac operator effectively  accommodates symmetries of underlying (model) manifolds. 
 
 For instance, one can not prove with minimal hypersurfaces 
that {\it no metric} $g\geq g_0$ on $S^n$, where  $g_0$ is the standard metric with the sectional  
curvature 1, {\it can have $Sc(g)\geq n(n-1)= Sc(g_0)$. \footnote {Such a proof may be possible for $n=3$ with suitable boundary conditions for
minimizing surfaces.}}

\vspace {1mm}

{\large \sf Specific  Problem.}  Let $Z\subset S^n$ be a closed subset
  of codimension $k\geq 2$,  let  $X $ be  an orientable $n$-dimensional  Riemannian manifold and let 
$$f: X \to S^n\setminus Z$$
be  a smooth  proper map of non-zero degree which is  distance decreasing or, more generally, area decreasing. \vspace {1mm}

\hspace {7mm} {\large \sf  When and how can one bound the scalar curvature of $ X$?}\vspace {1mm}

{\it Example.}  If $Z$ is a piecewise smooth one-dimensional subset (graph)  with trivial Levi-Civita  holonomies along all it cycles, e.g. a disjoint union of  trees,   and if  $ X$ complete, then  -- compare with remark (a) in section 3,
$$ \inf_{x \in X}Sc( X)( x ) < n(n-1) =Sc(S^n).$$

{\it Proof.} Let  $\epsilon : S^n\to S^n$ be  an arbitrarily small perturbation of  the identity map which sends a small neighbourhood of $Z$ to $Z$. Then the  bundle $L$ on $ X$ which is  induced from $\mathbb S^+(S^n)$ by the composed map $ \epsilon \circ f: \underline X \to S^n$ is trivial at infinity and the above proof of $\bigcirc$ applies. 
 
More generally, the same argument applies to closed subsets  $Z\subset S^n$   admit  sequences of maps 
$$\epsilon_i:  S^n\to S^n$$
such that

$\bullet$ the maps $\epsilon_i $ send small neighbourhoods of $Z$ in $S^n$ to subsets $Z_i\subset S^n$ as  above, namely  i.e. piecewise smooth with  trivial holonomies over all cycles in $Z_i$;

$\bullet$ the maps $\epsilon_i $  converge, for $i\to \infty$, to the identity map in the $C^1$-topology.  

 {\it Questions.} (a) Can one more effectively describe these $Z$ e.g. those of the topological dimension zero? 

(b) Does the above inequality $ \inf_{x }Sc(X)( x ) < n(n-1)$  hold true for smooth closed curves $Z\subset S^n$, $n\geq 3$, with {\it non-trivial} holonomy?.

(c) Does $S^n$ minus a point admit a {\it incomplete} metric $g\geq g_0$  with $Sc(g)>n(n-1)=Sc(g_0)$ (where $g_0$ is the spherical metric)?


\vspace{2mm}

Let us generalise the class of overtorical manifolds $X$,  where non-zero multiples of the   fundamental cohomology classes, denoted  $[X]^\circ\in H^n(X;\mathbb Z)$,  decompose   
   into products of one dimensional classes, 
$$\mbox {$k[X]^\circ =h_1\smile...\smile h_n$,\hspace {1mm} $h_i\in H^1(X;\mathbb Z)$},$$
as follows.
\vspace {1mm} 

\textbf {(10) \large \sf  Oversymplectic Manifolds.} A compact  orientable  $n$-dimensional manifold $ X$ is oversymplectic if  a multiple  of the  fundamental cohomology class  of $X$, decomposes into product of one and two dimensional classes,
$$k\cdot [X]^\circ=h_1\smile...\smile h_m,  $$
and such an $X$ is called [{$\tilde \uparrow 0$}]-oversymplectic, if  \vspace {1mm}

{\it the  classes $h_i$   vanish in the cohomology of the  universal covering $\tilde X$ 

under the natural 
homomorphism 
$H^\ast(X)\to H^\ast (\tilde X)$.} \vspace {1mm}

Notice that  [{$\tilde \uparrow 0$}]  is automatic for 1-dimensional classes.

Also note that if $n=2m$,  then,  by   grouping 1-dimensional $h_i$ into pairs, one  can 
make all  $h_i \in H^2(X;\mathbb Z)$, $i =1,...,m$, and if $n=2m+1$ all but one among $h_i$ can be brought to $H^2(X;\mathbb Z)$.

Moreover,   the a priori different  2-dimensional classes  $h_i$,  can be replaced by a single one, namely by a generic linear combination $h$ of $h_i$, since 
$\smile_ih_i=k'\cdot h^{\smile m}$.

It follows that $X$ of dimension $n=2m$  is oversymplectic if and only if it admits a map of  {\it  non-zero degree}
to the complex projective
space $\mathbb CP^m$,  where the condition 
[{$\tilde \uparrow 0$}] says in effect  that {\it the  pull back of the symplectic (K\"ahler)  2-form on $\mathbb CP^m$  to the universal covering $\tilde X$ of $X$ is exact.} 

And if $n=2m+1$ is odd, there is such a map $X\times S^1\to \mathbb CP^{m+1}$.

Observe that [{$\tilde \uparrow 0$}]-oversymplecticity,   similarly to  overtoricity and to iso-enlrageability of manifolds $X$ is inherited by $X'$ which admit   maps $X'\to X$  of non-zero degrees and also
by the products $X'=X\times \mathbb T^k$.  

Still,  [{$\tilde \uparrow 0$}]-oversymplicity  seems significantly different from iso-enlargeability, and, probably, there are many  examples of  [{$\tilde \uparrow 0$}]- [{$\tilde \uparrow 0$}]-oversymplectic manifolds, even among projective algebraic ones,  which are not (iso)enlargeable.

 \vspace {2mm}
 
 The reason we brought forth  this  oversymplecticity is the  following proposition. \vspace {1mm}

({\Large $\star$}$\tilde \uparrow 0$) {\sl If $X$ is [{$\tilde \uparrow 0$}]-oversymplectic, then it admits  no metric with $Sc>0$, provided
the universal covering $\tilde X$ is spin.} \footnote{$\tilde X$ is spin if and only if  the restrictions of the tangent bundle $T(X)  $  to all 2-spheres in $X$ are trivial; if $n\geq 5$ this is the same as triviality of the normal bundles of embedded 2-spheres in $X$.} \vspace {1mm}

(This, as it was   mentioned earlier, implies  that  the only possibility for $Sc(X)\geq 0$ is $X$ being flat.

Also recall  that  vanishing of the second homotopy group $\pi_2(X)$ implies that  $\tilde X$ is spin  and observe that   $\pi_2(X)=0$  also implies [{$\tilde \uparrow 0$}].)

\vspace {1mm}

{\it Proof of} ({\Large $\star$)}. Let $n=2m$ and  $\tilde l$  be the lift of the canonical line bundle of 
$\mathbb CP^m$ to $\tilde X$. Because of  [{$\tilde \uparrow 0$}], this bundle is trivial  there are the  $p$-th order roots $\sqrt[p]{\tilde l}$  for all $p=1,2,...$, which are represented by  the $p$-sheeted coverings of the total space of the circle bundles associated to $\tilde l$. 

And albeit the Galois' actions of the fundamental group $\Gamma= \pi_1(X)$ on $\tilde X$ and on ${\tilde l}$ does not extend  to  $\sqrt[p]{\tilde l}$, the semidirect product $\mathbb Z_p \rtimes \Gamma$ does act on $\sqrt[p]{\tilde l}$.

Since $\tilde X$ is spin, the  twisted  Dirac operator $D_{\otimes \sqrt[p]{\tilde l}}$, i.e. $D$  with coefficients in $\sqrt[p]{\tilde l}$, is defined and the corresponding space  of harmonic $L_2$-spinors is acted upon by the group  $\mathbb Z_p\times \mathbb  Z_2 \rtimes \Gamma$.

Then an elementary computation shows that the topological index $D_{\otimes \sqrt[p]{l}}$ does not vanish for infinitely many $p$ and then, by 
the Atiyah  $L_2$-index theorem,  
$D_{\otimes \sqrt[p]{\tilde l}}$-harmonic $L_2$-spinors exist for arbitrarily large $p$.

But since the curvatures of the bundles $\sqrt[p]{\tilde l}$ tend to $0$ for $p\to \infty$, uniform positivity of the scalar curvature of $\tilde X$   would not allow such spinors for large $p$ according to  the twisted    Schroedinger-
Lichnerowicz-Weitzenboeck vanishing theorem. QED.

\vspace {1mm} 

\textbf {11  \large \sf  Continuation of Discussion.}  On the surface of things, 
({\Large $\star$}$\tilde \uparrow 0$) generalizes Schoen-Yau theorem on  non-existence 
of metrics with $Sc>0$ on overtorical manifolds, but...

(1)Here again  there is an annoying spin condition in the statement of  ({\Large $\star$}$\tilde \uparrow 0$}), which, for all we know must be  unnecessary.

(2) More seriously, we can say preciously     little  about  incomplete manifolds.

For instance,   \vspace {1mm}

{\sf one can't    bound  with the present day techniques  the width  of  product bands    $(Y\times \mathbb [-1,1],g)$
with 
metrics $g$ 
which have  $Sc(g)\geq \sigma>0$  for  [{$\tilde \uparrow 0$}]-oversymplectic manifolds $Y$.} \vspace {1mm}

(The same can be said about all  other  non-$\cal SYSE$-manifolds $Y$   which are  known not to  to admit  metrics with $Sc>0$).

Because of this, \vspace {1mm}

{\sf one is unable to    rule out complete  metrics with $Sc>0$ on $Y\times \mathbb R$ 
and  complete

 metrics with $Sc\geq \sigma>0$ on $X\times \mathbb R^2$ for [{$\tilde \uparrow 0$}]-oversymplectic  manifolds $Y$.}\vspace {1mm}

What is not hard to show,  however, is the following \vspace {1mm}

({\Large $\star$}$_{\times \mathbb R}$) {\it  Products  $X=Y\times \mathbb R$ carry no complete metrics  $g$ with $Sc(g) \geq \sigma>0$ 
for all [{$\tilde \uparrow 0$}]-oversymplectic manifolds $Y$ the universal coverings of which are spin.} \vspace {1mm}

{\it Sketch of the Proof.}  Since    $X= (X,g)$ is complete and two-ended, it admits a proper  $1$-Lipschitz function onto $ \mathbb R$, which we  call it $\phi:X\to \mathbb R$. 

Let $X'=X\times \mathbb R$ and let
$$\Phi_\varepsilon =(\varepsilon\cdot \phi, \varepsilon\cdot \phi'): X'\to \mathbb R^2 $$ 
where $\phi': X'=X\times \mathbb R  \to  \mathbb R$ is the  coordinate projection.

Let $ l'_\varepsilon$  be the $\Phi_\varepsilon$-pullback of   $ l_0$ to $X'$ and let 
$ l^\circ_\varepsilon$ be the formal difference between  $l'_\varepsilon$ and the trivial complex line bundle with the trivial connection. 

Let $l_0$  be a complex  line bundle over $ \mathbb R^2$ with a unitary connection, 
which is isomorphic to the trivial bundle outside a compact subset in $\mathbb R^2$ and such that the curvature $\omega_0$  of $l_0$ is  $\omega_0=p_0(t_1,t_2) dt_1\wedge dt_2$ for a non vanishing function $p_0\geq 0$.

Let $dim(Y)=2m$, let $l$ be the line bundle over  $X'$ induced by the composed map   $X'\to Y\to \mathbb CP^m$  from the canonical line bundle and let
$$ l^\circ_\varepsilon =l\otimes   l^\circ_\varepsilon.$$ 
  
Pass to the universal covering $\tilde X$  and observe as earlier, that  
Atiyah's  $L_2$-index theorem, applied  to   the Dirac operator twisted with $\sqrt[p]{\tilde l^\circ_\varepsilon}$  and  combined with   Schroedinger-
Lichnerowicz-Weitzenboeck vanishing theorem  for small $\varepsilon<< \sigma$ and  for $p\to \infty$,  rules out  $Sc\geq \sigma>0$ for  complete metrics  on $X$.   QED

  \vspace {1mm}

{\it Generalisation to Non-Compact $X$.}  The above  ({\Large $\star$)} generalises to complete [{$\tilde \uparrow 0$}]-oversymplectic manifolds $X$, where the fundamental class 
of $X$ in the cohomology   with {\it compact supports}, denoted $H^n(X,[\infty ])$,\footnote{ This $[\infty]$ stands for  the complement to a (large)  non-specified  compact subset in $X$.} decomposes into 1- and 2-classes also with compact supports   and where these classes must vanish in the cohomology 
$H^{1,2}(\tilde X , \widetilde {[\infty ]})$.

For instance, if $X$ of dimension $2m$  admits a proper map of non-zero degree to a complement of a subset $Z\subset \mathbb CP^m$, this condition is satisfied if $H_1(Z)=0$ and 
the symplectic form of $\mathbb CP^m$ vanishes on $Z$.

\vspace {1mm}
\textbf {12 \large \sf  Min-Oo - Goette - Semmelmann Rigidity Theorem.} A (very) special case of this theorem 
(2.10 in [GS 2001]) reads as follows.

Let $X$ b a compact orientable Riemannian manifold of dimension $2m$ and let $f: X\to \mathbb CP^m$ be a $C^1$-smooth area non-increasing {\it  spin map} of non-zero degree
where $f$ is called  {\sl spin} if  the restriction of the tangent bundle  $T(X)$ to  $\Sigma\subset X$ is trivial if an only the restriction  $T(\mathbb CP^n)_{|f(\Sigma)}$  is trivial for all surfaces $\Sigma\subset X$.

For instance, if $m$ is odd  and  $X$ is  spin then all maps $X\to \mathbb CP^n$ are spin.

\vspace {1mm}

{\Large $\ostar $}  {\it If $$Sc(X)(x)\geq Sc(\mathbb CP^m)(f(x)) \mbox {   for all } x\in X,$$ 
then $f$ is an isometry.}\vspace {1mm}

{\it About the Proof.}   The $\cal R$-term in the  Schroedinger-Lichnerowicz-Weitzenboeck  formula \textbf {(5)}) for $D$   twisted with {\it line bundles} $l$    shows (see  [Hit 1974])   that that  if the curvature form $\omega$ of an  $l$  
(where the cohomology class of $2\pi \omega$ equals the Chern class $c_1(l)$)
on a Riemannian manifold of dimension $2m$ diagonalises as
$$ \sum_{i=1,...,m}\lambda_ie_{2i-1}\wedge e_{2i}$$
for an orthonormal frame $e_1,e_2,..., e_{2m}$,
then 
$$||{\cal R}|| \leq  4\sum_i\lambda_i .$$
 
 If $l$ equals the (anti)canonical bundle $l_0=\bigwedge^m T(\mathbb CP^m)$
then, its curvature form  for the Levi- Civita connection of the Fubini Studi metric $g_0$ has 
$$\lambda_1=\lambda_2=... =\lambda_m=  m+1$$
and $g_0=m(m+1)$.

On the other hand,  an  easy homological computation shows that   the topological index $I(D,l_0)$ does not vanish on $\mathbb CP^m$  and since $deg(f)\neq 0 $  it doesn't vanish on $X$ either.  This  shows that $f$ can't be {\it  strictly}  area decreasing, while the equality case needs an additional  argument (see [GS 2001]).

The above applies, strictly speaking, to odd $m$, where $\mathbb CP^m$ is spin, and if $m$ is even,  one twists $D$ with the {\it virtual square root of} $l_0$ (see [Hit  1974],   [Min-Oo 1995], [GS 2001]).

\vspace{1mm} 

\textbf {(13) (\large \sf  Interpolating between  ({\Large $\star$}$\tilde \uparrow 0$)  with  ({\Large $\ostar$})} Unlike the (obvious) implication $\bigcirc\Rightarrow \square$
the (sharp) theorem  
 ({\Large $\ostar$})  by no means  implies (rough)  ({\Large $\star$}$\tilde \uparrow 0$). 
 
 But an obvious combination of the proofs of these theorems   brings the two  together  as follows.
 \vspace {1mm}

Let $X$  be a  complete  oriented Riemannian $2m$-manifold  and $f: X\to \mathbb CP^m\setminus Z$ be a proper  $C^1$-smooth  {\it area non-increasing} map of {\it non-zero} degree, where
$Z\subset \mathbb CP^m$ a smooth  submanifold on which the {\it symplectic form of $\mathbb CP^m$ vanishes} and which has $H_1(Z)=0$.

Let the composed map $\tilde X\to X\to    \mathbb CP^m\setminus Z\subset \mathbb CP^m $ from the universal covering of $X$ to  $\mathbb CP^m $, call it $\tilde f:\tilde X\to \mathbb CP^m$,  be spin.  \vspace {1mm}

 {\Large ($\ostar$}{\large $\tilde \uparrow\frac{1}{p}$)}  {\it If the $\tilde f$-pullback  of the generator  $c\in H^2( \mathbb CP^m);\mathbb Z)$, that is   $\tilde f^\ast(c)\in H^2(\tilde X;\mathbb Z)$, is divisible by a positive integer  $p$, then 
$$\inf_x Sc(X)(x) < \frac {1}{p} Sc(\mathbb CP^m),$$ 
unless $X$ is compact, $Z$ is empty, $p=1$ and $f$ is an isometry.} 

 \vspace {1mm}

One can only wonder if there is anything of this kind  that may come  from   minimal hypersurfaces.
\section {Appendix}

In this section, following the suggestions by the referee,  we  explain in a greater detail the following.

\vspace {1mm}

(a)  {\sl Smoothing  
 hypersurfaces  with  {\it no    decrease} of their mean curvatures.} \vspace {1mm}

(b)  {\sl Smoothing  {\sl  Riemannian metrics} with  with  {\it no    decrease} of their  scalar 
  
  curvatures.} \vspace {1mm}
 
In  the case (a)  the  principal  step of smoothing   consists of  "rounding corners" of piecewise smooth hypersurfaces which makes our hypersurfaces $C^1$-smooth.  

In the case (b) one  $C^1$-smoothes continuous piecewise smooth Riemannian metrics by
similarly "rounding them by bending" along their singular loci.

Both rounding constructions, however simple, depend on {\it specific} geometric properties of  the mean curvature and the scalar curvature correspondingly. 

Then, in both cases, the final step of smoothing   $C^1\leadsto C^\infty$ follows by  {\it   homotopy extension} construction   for solutions of {\it  general} differential inequalities which is  explained  in section 11.1 below.



\vspace {1mm}


 Next,   we   \vspace {1mm}
 
 (c)  {\sl elucidate  some    properties of  on minimal hypersurfaces in open manifolds needed for  the width inequalities}. \vspace {1mm}

 Finally, as it was  suggested by the  referee,   we  \vspace {1mm}

 (d) {\sl  summarise   topological obstruction for $Sc>0$ on closed manifolds    which follow from our inequalities }\vspace {1mm}
 
\hspace {-6mm} and  \vspace {1mm}

(e) {\sl highlight  several   conjectures mentioned in the main body of the article.}

\subsection {Universal Constructions of Smoothing  and Bending}

{\large \sf A. Linear Smoothing Operators.}  The most common  kind of   smoothing in linear  analysis is  achieved by applying    convolution-like  operators  to objects you want to smooth.  

For instance, let $X$ be a Riemannian manifold and 
  $Y\subset X$ be a compact   $C^{1}$-smooth 
cooriented hypersurface. 
Let $U\supset Y$ be a small  $C^\infty$-split neighbourhood of $Y$, say
$$U= \underline Y\times (-\delta, \delta),$$
where $Y$ is represented by a graph of $C^1$-smooth function $f( \underline y)$ and let $Y_\varepsilon$, $\varepsilon>0$ be the graphs of the   the $\varepsilon$-smoothed functions  $f_\varepsilon (\underline y)$, for
$$f(\underline y)\mapsto f_\varepsilon (\underline y) = \int_{\underline Y} K_\varepsilon(\underline y,\underline y') f(\underline y') d\underline y'.$$
Then, for the usual  $K_\varepsilon$,
$$Y_\varepsilon \underset {C^1}\to  Y \mbox  { for } \varepsilon\to 0$$  
and, since the mean curvature of $Y$  is {\it linearly} expressible in terms of the second derivatives of $f$,   the mean  curvatures of $Y_\varepsilon$ satisfy almost the same bounds as 
those of $Y$, whenever the latter are defined. 

For instance, if  $Y$ is  piecewise $C^2$, or more generally   if the above ($C^1$-smooth)  function $f$ is $C^{1,1}$,  i.e. $df$ is   Lipschitz, and if 
$$mn.curv(Y)(y)\geq \phi(y)$$
for a continuous function $\phi$ on $U\supset Y$ and almost all $y\in Y$, 
then 
$$mn.curv(Y_\varepsilon)(y_\varepsilon)\geq \phi(y_\varepsilon) +o(1), \hspace {1mm}  \varepsilon \to 0,$$
for all $y_\varepsilon \in Y_\varepsilon.$
(A slightly different format of this smoothing is suggested on    pp. 939 and 949 in section 3.4 and 5.7 in  [Gr 2014].) 

Similarly,  $C^{1,1}$-smooth Riemannian metrics $g$  on a manifold  $X$ can be $C^1$-approximated by 
$C^\infty$-smooth $g_\varepsilon$  such that if  $Sc(g)(x)\geq \phi(x)$, then 
$$Sc(g_\varepsilon)\geq \phi(x)+o(1),$$
because the operator $g\mapsto Sc(g)$  is linear in the second derivatives of $g$.

\vspace {2mm}

{\sf \large B. Smoothing by Local Bending.} Smoothing a function $f$ on $V$, where for instance, the second derivatives jump across a  double sided hypersurface $\Sigma\subset V$, can be achieved by deforming,  we call it {\it bending}, $f$  on one side of $V$, say to the "left"  of  $\Sigma$, such that the derivatives  on the left side  become  equal to the derivatives  on the right other side of $\sigma$ at all points $v\subset \Sigma$.

Linearity of the differential inequality we want to keep preserved by such bending,
e.g. of $mn.curv> \phi$ or $Sc> \phi$,
 is not indispensable.  What is essential for $C^1\leadsto C^\infty $ smoothing  is  the {\it connectivity} rather than convexity of the subsets in the sets of  values of derivatives of functions and/or of  metrics defined by the required  inequalities.

Let us formulate the relevant general bending property of solutions of such inequalities in terms of "cut-offs of deformations",  where  "connectivity", is hidden   the concept of "deformation/homotopy".

\vspace {1mm}

Let $Z\to V$ be a smooth fibration and let $Z^{[r]}$  be that space of the $r$-jets of  germs of $C^r$-smooth   sections $f:V\to Z$.  

 Let $J^r_f:V\to Z^{[r]}$ denote the $r$-th jet of $f$ and 
recall that by the definition of jets,

{\sl $J^r_{f_1}(v) = J^r_{f_2}(v)$ if and only the values of the sections  $f_1$ and $f_2$  as well as of  all their partial derivative of orders 1,2,...r,  in some local  coordinates, are equal at  $v$. }

In fact, this property {\it defines} jets as  as well as the spaces  $Z^{[r]}$.

Let $\Sigma\subset V$ be a piecewise smooth subset,\footnote{Probably, what follows holds true for all closed subsets $\Sigma\subset V$.} let $ U\supset \Sigma$ be a neighbourhood of $\Sigma$ in $V$ and let 
${\cal R}\subset  Z^{[r]}$  be an open subset.\vspace {1mm}

\textbf {$\bigstar$   Cut-off Homotopy   Lemma.}\footnote{This is a reformulation of {\sl the  weak flexibility lemma} given  as an exercise   on  p. 111 in [Gr 1986].}
 Let $f_t: V\to   Z$, $t\in  [0,1]$,  be a $C^r$-continuous family of smooth sections, such that  

(i)  $J^r_{f_t}(V)\subset \cal R$ for all $t\in [0,1]$

and 

(ii)  $J^{r-1}_{f_t}(v) $ is constant in $t$ for all $v\in \Sigma$. \vspace {1mm}

{\it Then there exists a smaller  neighbourhood   $\bar U\subset  U$ of $\Sigma$ and another   $C^r$-continuous family of smooth sections $\bar f_t:V\to Z$, such that, similarly to $f_t$, 
$$J^r_{\bar f_t}(V)\subset {\cal R} \mbox  { for all } t\in [0,1]$$
and

 \hspace {35mm}  $\bar f_t$ is equal to $f_t$ on $\bar U$\vspace {1mm} 

 \hspace {-6 mm}while  at the same time \vspace {1mm}

  \hspace {35mm} $\bar f_t$ is constant in $t$ outside $U$,\vspace {1mm}

 \hspace {-6mm}   i.e.  $\bar f_t(v)= f_0(v)$ for all $t$ and all $v\in V\setminus U$.

Moreover, if $f_t$ was  constant in the neighbourhood of a closed subset $V_0\subset V$ then $\bar f_t$  can be taken constant in $t$ on  $V_0$.}
\vspace {1mm}

{\it Remark.}  The general case of the lemma easily  reduces to that where 
$V$ is {\it  compact} and $\Sigma$  is a  {\it smooth submanifold} of codimension {\it one}.

 \vspace {1mm}
 
 {\sf \large  Warning and Perturbative Generalisation.} The conclusion of $\bigstar$ by no means holds true in general without assumption  (ii). However,

\hspace {5mm}{\it   "constant" in (ii) can be replaced by "almost constant" as follows.}

 \vspace {1mm}
 
\textbf {$\bigstar'$ } {\sl Let $f_{t, \theta}$ $t, \theta\in [0,1]$ be a two parameter $C^r$-continuous family
of sections $V\to Z$ where $f_{t,0}$ satisfies (ii).
Then there exists an $\varepsilon_0 >0$. such  that  the conclusion of  $\bigstar$ holds for  $f_{t, \theta_0}$ 
 for all $\theta_0\leq \varepsilon _0$.}
 
 \vspace{1mm}
  This follows from \textbf {$\bigstar$ }  and from what is called {\it microflexibility} of differential inequalities defined by {\it open} subsets in the jet spaces, which is, of course, fully trivial.
  
(More interestingly,   there are  classes of  {\it flexible}  maps, e.g. {\it smooth immersions } $V^n \to \mathbb R^{n+1}$, which are defined with certain   $\cal R$, where  extension of  homotopies $f_t$ (called {\it regular homotopies} for immersions) is possible  for {\it all} $f_t$.    But the proofs of this       in interesting cases don't, unlike $\bigstar$ and $\bigstar'$, reduce to    generalities but depend on  specific  constructions adapted to  specific properties    of particular  $\cal R$. See [Gr 1986] and references therein.)

 .

\vspace{1mm}

1-D -{\it Example.} Let $V=[0,\infty)$ and $Z=[0,\infty)\times \mathbb R\to [0,\infty)$ be the 
trivial fibration. Then sections of $Z$ correspond to real functions $f(v)$, $v\geq 0$, 
$$  Z^{[r]}=[0,\infty)\times \mathbb R^{r+1}$$
and the $r$-jets are maps
$$ J^{r-1}_{f} =\left (f, \frac {df}{dv},...,  \frac {d^rf}{dv^r}\right): [0,\infty) \to \mathbb R^{r+1}.$$

{\Large $\star$}  {\sl Given continuous functions  $a_0(v),..., a_r(v)$,   $b_0(v) > a_0(v),..., b_r(v)>a_r(v)$,
a number 
$\delta>0$ and a real function $f(v)$, $v\geq 0$, such that 
$$ a_i(v)< \frac {d^if(v)}{dv^i}<b_i(v) \mbox {  for $i=0,...,r$ and  all } v\geq 0,$$
there exists  a function $\bar f(x)$, which also satisfies these inequalities, 
$$ a_i(v)< \frac {d^i\bar f(v)}{dv^i}<b_i(v) \mbox {  for $i=0,...,r$ and  all } v\geq 0,$$
and such that

\hspace {19mm} $  \frac { d^r \bar f(0)}{dv^r}= c_r$ and  $\bar f(v)=f(v)$  
 for $v\geq \delta$,\vspace {1mm}

\hspace {-6mm}where $c_r$ is a given number in the interval $a_r(0)<c_r< b_r(0)$. }
\vspace{1mm}

The essential and essentially  obvious  case of this example is where $r=1$, $a_0=a_1=0$ and $b_0=b_1=\infty$. When this is understood, all of $\bigstar$ becomes obvious  as well.\footnote{This "obvious" presupposes familiarity with the basic  geometry of the jet  spaces.}

$\Leftcircle \hspace {-3mm}\Rightcircle${\it Mean Curvature Example.} Let $Y\subset X$ be a piecewise smooth $C^1$-hypersurface, where the singular locus is a smooth  hypersurface $\Sigma\subset Y$ (i.e. $dim(\Sigma)=dim(Y)-1$), where the  two smooth parts, say $V_1$ and $V_2$ of  $Y$ meet. This $Y$ can be obviously $C^\infty$-smoothed {\it along $\Sigma$} by deforming $V_1$ and $V_2$  near $\Sigma$  such that they would  $C^\infty$-match  at $\Sigma=V_1\cap V_2$ and  with almost 
no decrease of  their   mean curvatures.\footnote 
{This deformation at a point $\sigma\in \Sigma=V_1\cap V_2$  does not $C^2$-significantly  move  a  neighbourhood  $U_1\subset  V_1$  of $\sigma \in V_1$, unless $mn.curv(V_1)(\sigma)<  mn.curv(V_2)(\sigma)$ and the same applies to  $V_2$.}

Then $\bigstar$ -- here $r=2$ --   allows an  extension of these deformations/homotopies    to {\it all of $Y=V_1\cup V_2$}, such that the resulting  smoothed hypersurface, say $\bar Y=\bar Y_\varepsilon$, for given    positive  continuous  function $\varepsilon= \varepsilon(x)>0$ on $X$, satisfies.

$\bullet$  The $C^1$-distance between $\bar Y_\varepsilon$ and $Y$ is $\leq \varepsilon$.

$\bullet$ $\bar Y_\varepsilon$ coincides with $V$ outside the  $\varepsilon$-neighbourhood of $\Sigma$.

 $\bullet $ The mean curvatures of $\bar Y_\varepsilon$ satisfy the same, up to $\varepsilon$, inequalities as 
  the mean curvatures of $Y$. 
  
  In particular,  if $mn.curv(Y)(y)\geq \phi(y)$ for a continous function $\phi $ on $X\supset Y$, then $mn.curv(\bar Y_\varepsilon )(\bar y)\geq \phi(\bar y)+\varepsilon (\bar y)$. \vspace {1mm}
 
 Similarly  one can smooth  more general {\it piecewise smooth} $C^1$-hypersurfaces and also   {\it piecewise smooth $C^1$-Riemannian metrics} with {\it negligible  decrease} of their scalar curvatures.

\subsection {$\varepsilon$-Redistribution of Curvature} In  smoothing and bending constructions it is  easier  deal with   with strict inequalities, such as $Sc>0$, rather than with non strict ones, such as $Sc\geq 0$

Below are two simple (and, probably,  known)  propositions which allows one to relax "partially strict" inequalities.
\vspace {1mm}

{\sf \large  Redistribution of the Mean Curvature.} Let $X$ be a $C^\infty$-smooth Riemannian $n$-manifold and $V\subset X$  a  domain with  {\it cosimplicial corners}, i.e. each point at the boundary of $V$ admits a neighbourhood in $V$ which is diffeomorphic to  a neighbourhood  in the positive "octant"  $\mathbb R^n_+\subset \mathbb R^n$.

Let $\partial_i\subset \partial=\partial V$ denote  the $(n-1)$-faces of $V$ and $\partial_{ij}$ be the  $(n-2)$-faces
which we call {\it edges}.
Let  $\phi_i$ and $\psi_{ij}$ be smooth functions on $X$, such that \vspace{1mm}

(i)  the mean curvatures of the $(n-1)$-faces of $V$  satisfy 
$$mn.curv(\partial_i )(x)\geq \phi_i (x),$$
for  all $\partial_i $ and all $x\in\partial_i$; 

(ii) the dihedral angles between the $(n-1)$-faces  $\partial_i,\partial_j\subset \partial V$, 
satisfy
$$\angle_v(\partial_i,\partial_j)\leq \psi_{ij}(x)$$
for all  (non-empty) edges   
$\partial_{ij}=\partial_i\cap\partial_j$
and all $x\in\partial_{ij}.$
\vspace {2mm}

({\Large $ \star_>$})  If  $V$ is {\it compact} and the boundary $\partial V$ of $V$ is {\it connected}, then \vspace {1mm}

{\it either 
$$mn.curv(\partial_i )(x)= \phi_i (x)\mbox {   and } \angle_x(\partial_i,\partial_j)= \psi_{ij}(x)$$ 
for all points $x\in \partial V$} (where these equalities  make sense),\vspace {1mm}

\vspace {1mm} 

{\it or there exists an arbitrarily small $C^\infty$-perturbation $V'$ of $V$ (by a 
$C^\infty$-diffeomorphism close to the identity), such that 
 $$mn.curv(\partial'_i )(x)> \phi_i (x) \mbox {   and }\angle_x(\partial'_i,\partial'_j)< \psi_{ij}(x)$$
for all   $\partial'_i\subset \partial V'$  and  $ \partial'_{ij}=\partial'_i\cap\partial'_j$ and all $x\in \partial'_i$  and  $x\in \partial'_{ij}$  correspondingly.}

\vspace {1mm}

{\it Sketches  of Three Different  Proofs.} \textbf{(1)} {\sl Fredholm+Unique Continuation.} To get the idea, let $Y=\partial V$ be smooth  and  let $\phi_0(x)$ be a smooth function on $X$ which extends  the function $y\mapsto mn.curv(\partial V)(y)$ from $Y=\partial V$ to $X$.
Let $U_0\subset  X$ be a neighbourhood of a point  $y_0\in Y$.\vspace {1mm}

If a smooth function $\phi_0'$ is sufficiently $C^\infty$-close to $\phi_0$, then there exists 
\vspace {1mm}

({\Large $ \star_=$}) {\sl a $C^\infty$-perturbation $Y'$ of the hypersurface $Y=\partial V$, such that 
$$mn.curv(Y')(x)=\phi_0'(x)$$
for all $x\in Y' \setminus U_0$.}\vspace {1mm}

This follows by the implicit function theorem for the operator $Y\overset {\cal M} \mapsto mn.curv(Y)$, 
since

$\bullet_1$ the linearisation $L=L_{{\cal M}, Y} $ of  $\cal M$  at $Y$, being Fredholm, has {\it finite codimensional}  image;

$\bullet_2$ the adjoint operator of  $L$ (which happens to be equal to $L$) has the {\it unique continuation} property  for {\it connected}\footnote{This property (obviously) fails to be true if  $\partial V$ is disconnected and    ({\Large $ \star$})  doesn't have to  hold anymore.} hypersurfaces $Y$.

\vspace {1mm}

This  ({\Large $ \star_=$})  also holds  for  certain domains where the boundary $\partial V$ is non-smooth. For instance if $V$ has no corners, i.e. if  there is no triple intersections of faces,  $\partial _{ijk}=\partial_i \cap \partial_j\cap \partial_k $, then a version of 
 ({\Large $ \star_=$}),  which   is significantly stronger than ({\Large $ \star_>$}),  follows by 
 by     perturbing the faces $\partial_i$ one by one.

 But if  the linearized  boundary value problem  loses regularity at the corners (this, probably,  doesn't happen if  all dihedral angles $\angle (\partial_i,\partial_j)$ are $90^\circ$), then it becomes  unclear  
 if  ({\Large $ \star_=$}) remains true.

However,  ({\Large $ \star_>$}) is taken care of by the following.\vspace {1mm}

\textbf {Local Mnc-Lemma.}\footnote  {Albeit 100\%  elementary, this lemma heavily relies on the specifically    Riemannian/Pythagorean nature    of our problem.} Let $Y$ be a smooth hypersurface  in a Riemannian  manifold $X$ of dimension $n$. 

Then there exists a continuous   function $\varepsilon=\varepsilon _{X,Y}(y)>0$, on $Y$, which is  also continuous with respect to the $C^\infty$-topology in the space of hypersurfaces  $Y\subset X$
with the following property. 

{\sl  Let  $S_y(\epsilon)=S^{n-2}_y(\epsilon)\subset Y$ be the $\epsilon$-sphere around a point $y\in Y$ for some positive $\epsilon \leq \varepsilon$  and let $\delta>0$ be a positive number.

Then there  exists a diffeomorphism  $\phi :Y\to X$, such that \vspace {1mm}

$[\star_\delta]$  $\phi$ is $\delta$-close to the identity diffeomorphism $id:Y\to Y\subset X$ in the   $C^2$-topology; \vspace {1mm}

$[\star_\epsilon]$  the diffeomorphism $\phi$ is equal  to the identity outside a narrow band around  the sphere  $S_y(\epsilon)$,\vspace {1mm}

\hspace {40mm} $\phi(y)=y,$  \vspace {1mm}

 unless  
 $$v\in  B_y(\epsilon) \setminus ( S_y(\epsilon)\cup  B_y(\epsilon')), \mbox  { where $ \epsilon'=\epsilon  \left (1-\frac {1}{10^{10n}}\right)$,}$$
  and where $B_y(\epsilon)$ is the ball bounded by $S_y(\epsilon)$;  

 \vspace {1mm}

$[\star_>]$ the diffeomorphism $\phi$ strictly increases the mean curvature of $Y$  in the interior of   $B_y(\epsilon)$ close to $S_y(\epsilon)$, namely at all points  
$$v\in  B_y(\epsilon) \setminus ( S_y(\epsilon)\cup  B_y(\epsilon'')) \mbox  { for $ \epsilon''= \epsilon\left (1-\frac {1}{20^{10n}}\right)$.}$$}

\vspace {1mm}

{\it The  proof} of the lemma  is accomplished   by applying the  initial  stage of {\it bending the hypersurface}  $V\setminus B_y(\epsilon''))$ near its boundary},   which is  described in the next section  and where the existence of the {\it initial}  bending  we need here  is fully obvious.  \vspace {1mm}

\textbf{(2)}  {\sl Spread of Positivity.} The above lemma allows an extension of the strict inequality   
$mn.curv(Y)(y)>\psi(y)$ from (controllably small) balls  of radii $\epsilon'$ in  hypersurfaces $Y\subset X$ to larger balls of radii $$\epsilon=\epsilon'\left(1-\frac {1}{10^{10n}}\right)^{-1}$$
and, consequently, from  arbitrary open subsets $U'\subset Y$  to  larger $U\subset Y$, without moving the complements $Y\setminus U$.
\vspace {1mm}

\textbf {Localization Corollary.} The local lemma allows  perturbations $Y'$  required by   ({\Large $ \star_>$}) to be  localised in a given domain $U\subset Y$, which, in turn,  yields   ({\Large $ \star_>$}) 
for {\it non-compact} $Y$.




 

\vspace {2mm}

\textbf{ (3)} {\Large $\Leftcircle$} {\sl Variational Proof of ({\Large $ \star_>$}).} One can also construct a perturbation of $\partial=\partial V$   with a required control over $mn.curv(\partial)$,  say in the  smooth (no edges) case, by  minimizing the functional
$$V\mapsto  vol_{n-1}(\partial V)-\int_V\varphi'(x) dx$$
for  a suitably chosen function $\varphi'(x)$.

To see what such a $\varphi'(x)$ should be,    let $\mu(x)$ be a smooth extension of the mean curvature function from $\partial$  to $X$ 
and let $\delta(x)$ be  the signed distance function to $\partial$, i.e. 
$\delta(x)=dist(x,\partial)$  for  $x\in V$ and  $\delta(x)=-dist(x,\partial)$ for $x\in X\setminus V$.

Observe that the original  $V$   locally strictly and stably minimizes  the functional 
$$V\mapsto  vol_{n-1}(\partial V)-\int_V\varphi(x) dx  \mbox { for } \varphi(x)=\lambda\delta(x)-\mu(x)$$ 
where    $\lambda=\lambda (X, \partial)$    is a sufficiently large  constant, namely  $\lambda >>\sup_\partial  (curv^2(\partial)||+\sup_X|Ricci(X)|$,
see section 10.2 in [Gr 2012].

Thus,  the  functional $V\mapsto  vol_{n-1}(\partial V)-\int_V\phi'(x)$  has a unique local minimum 
$V'$ whenever  $\varphi'$ is  sufficiently close to the above $\varphi$, and where, observe,
$$mn.curv(\partial V')(x)=\varphi'(x)\mbox {  for all } x\in \partial V'.$$

Then it is not hard to arrange such a $\varphi'$  that would make  the mean curvature of $V$ increase  outside $U_0$ (compare section 11.9).

\vspace {2mm}

{\sf \large  Redistribution of the Scalar Curvature.} Let, besides the above functions $\phi_i$ and $\psi_{ij}$
on $ X$,  we are given a continuous function $\sigma(x) $ such that the scalar curvature of the Riemannian metric $g$ in $X$ satisfies
$$Sc(g)(x)\geq \sigma(x)  \mbox{ for all }x\in X. $$

If $V$ is connected  (non-compact is allowed) then  \vspace {1mm}

\hspace {-6mm}either 
$$ Sc(g)(x)=\sigma (x)\mbox { for all $x\in V$}, \hspace {1mm} mn.curv(\partial_i )(x)= \phi_i (x), \mbox {   and } \angle_x(\partial_i,\partial_j)= \psi_{ij}(x),$$ 
where the latter two equalities hold for all points $x\in \partial V$, where they  make sense,\vspace {1mm}

\vspace {1mm} 

Moreover, \vspace {1mm} 

{\it one can keep $g'=g$ on a given closed subset in $X$ on which $Sc(g)(x)>\sigma (x)$.}\vspace {1mm} 

{\it Sketches of two    Proofs.} \textbf{[1]} If $V$ is  compact with smooth boundary, 
the   (linearization of the)  operator $f\mapsto Sc(f^2g)$, $f>0$ is Fredholm with the unique continuation property  and  the above argument via linearization applies.
\vspace {1mm} 

\textbf{[2]} The initial  stage of what is called "intrinsic bending" in section  11.5  yields the 
following simple  proposition  that fully accomplishes our purpose. (As in the mean curvature case, 
the existence of {\it initial} intrinsic  bending is obvious.)

\vspace {1mm}

\textbf {Local Sc-Lemma.}  Let $X =(X,g)$ be a smooth $n$-dimensional Riemannian  manifold. 
Then there exists a continuous   function $\varepsilon=\varepsilon _{g}(x)>0$  on $X$, which is  also continuous with respect to the $C^\infty$-topology in the space of  Riemannian metrics  $g$ on $X$, with the following property. \vspace{1mm}

{\sl  Let  $S_x(\epsilon)=S^{n-1}_x(\epsilon)\subset X$ be the $\epsilon$-sphere around a point $x\in X$ for some positive $\epsilon \leq \varepsilon$  and let $\delta>0$ be a positive number.

Then there  exists a  smooth Riemannian metric, $g'$  on $  X$, such that \vspace {1mm}

$[\star_\delta]$  the metric $g'$ is $\delta$-close to $g$  in the   $C^3$-topology; \vspace {1mm}

$[\star_\epsilon]$  the metric  $g'$ is equal  to $g$ in a narrow band around  the sphere  $S_x(\epsilon)$,\vspace {1mm}

\hspace {40mm} $g'(x)=g(x),$  \vspace {1mm}

 unless  
 $$x\in  B_x(\epsilon) \setminus ( S_x(\epsilon)\cup  B_x(\epsilon')), \mbox  { where $ \epsilon'= \epsilon\left (1-\frac {1}{10^{10n}}\right)$,}$$
  and where $B_x(\epsilon)$ is the ball bounded by $S_x(\epsilon)$;  

 \vspace {1mm}

$[\star_>]$ the scalar curvature of $g'$ is strictly greater then that of $g$   in the interior of   $B_v(\epsilon)$ close to $S_x(\epsilon)$, namely at all points  
$$v\in  B_x(\epsilon) \setminus ( S_x(\epsilon)\cup  B_x(\epsilon'')) \mbox  { for $ \epsilon''= \epsilon\left (1-\frac {1}{20^{10n}}\right)$.}$$}

\subsection {Rounding and Smoothing Hypersurfaces with no Decrease of their   Mean Curvatures.}

Let $X=X^n$ be a $C^\infty$-smooth Riemannian manifold of dimension $n$
e.g.  $X=\mathbb R^n$, and let $Y=Y^{n-1} \subset X$ be a  {\it cooriented} hypersurface which is 

\hspace {2mm}{\sl $Y$ is locally $C^\infty$-diffeomorphic to  a convex polyhedral hypersurface in $\mathbb R^n$.}\vspace {1mm}

In other words,  $Y$ is  {\it piecewise $C^\infty$-smooth}  with 
all dihedral angles   between the tangent spaces of the smooth regions  $Y_i\subset Y$ at their meeting points  in  the singular locus  of $Y$ being {\it strictly less than $\pi$.} 
$$\angle _{ij}=\angle (T_y(Y_i),T_y (Y_j))<\pi.$$

\vspace {1mm}

For instance,
 the boundaries $Y$ of finite   intersections of  domains  $U_i\subset X$   with smooth boundaries which  all  intersect {\it  transversally} are of this kind. \vspace {1mm}



Let us agree that  the mean  curvatures of  cooriented hypersurfaces  are evaluated  with the {\it outward looking} normal  vectors, where this sign  convention  makes  mean curvatures of  boundaries of {\it convex} domain {\it positive}.

\vspace {1mm}

\hspace {39mm} \textbf {$[^\varepsilon\leftrightarrows_\varepsilon]$-Rounding.} \vspace {1mm}

\hspace {5mm} {\sf Move  $Y$  inward equidistantly  by $\varepsilon$ and then by the same $\varepsilon$ outward.}\vspace {1mm}

If $Y$ is a  {\it compact closed}  hypersurface and {$\varepsilon>0$ is sufficiently small, then the  resulting hypersurface, call it $Y\pm \varepsilon]$, is

 \hspace {20mm}  {\it $C^1$-smooth   and piecewise $C^\infty$-smooth}. \vspace {1mm}

To see this clearly,  let $Y_{-\varepsilon]} $ be the  inward  $\varepsilon$-equidistant hypersurface to $Y$, which, in the case $Y=\partial U$, is equal to    the boundary of the 
$\varepsilon$-neighbourhood of the complement $X\setminus U$. 

Observe that  if $\varepsilon$ is small, than

 \hspace {20mm} {\it $Y_{-\varepsilon]} $  has the same  corner pattern as $Y$.} 

Now the  piecewise structure of $Y\pm \varepsilon]$, that is the outward $\varepsilon$-equidistant to $Y_{-\varepsilon]}$,  can be  seen with    {\it the normal (nearest point) projection}  $ Y_{\pm \varepsilon]} \to Y_{-\varepsilon]} $,  call it $p_{-\varepsilon] }:Y\pm \varepsilon]\to Y_{-\varepsilon]} $,   where each $y\in Y_{\pm \varepsilon]}$ is sent by   $ p_{-\varepsilon}$  to  the unique ($\varepsilon$ is small!) nearest point in $Y_{- \varepsilon]}$:\hspace {1mm} 

 {\sl the smooth pieces of   \hspace{0.5mm} 
 $Y_{\pm \varepsilon]}$  are the  
$p_{-\varepsilon] } $-pullbacks of the ($(k-1)$-codimensional)  "faces" $Y_{-\varepsilon]}^{n-k} \subset  Y_{-\varepsilon]}$,  $k=1,2,...,n$.}

For instance, if $Y$ is the boundary of  the intersection of smooth domains with {\it transversally meeting} boundaries, then 

{\it these faces correspond to meeting points of $k$-tuples of such boundaries. } 
 \vspace {1mm}

Observe that the subset $Y_{\pm\varepsilon]}^{n-1}= p_{-\varepsilon] }^{-1}\left (Y_{-\varepsilon]}^{n-1}\right )\subset Y_{\pm \varepsilon]}$ satisfies
$$Y_{-\varepsilon]}^{n-1}= Y\cap Y_{\pm \varepsilon]}$$
and that     it is also  equal 
 to the set of points $y\in Y$,  such that $dist(y, sing(Y_{-\varepsilon]})>\varepsilon$, where  $sing(Y_{-\varepsilon]})$ is the set of points where $Y$ is non-smooth, i.e.  union of  what we call   {\it the edges} or {\it $(n-2)$-faces}   of  $Y_{-\varepsilon]}.$

Also observe that  that the mean curvatures of the  $C^2$-smooth  pieces  $Y_{\pm\varepsilon]}^{n-k}=p_{-\varepsilon] }^{-1}\left (Y_{-\varepsilon]}^{n-k}\right ) $ of $Y_{\pm\varepsilon]}^{n-k} $ for $k\geq 2$   are {\it large positive,} 
$$mn.curv\left ( Y_{\pm\varepsilon]}^{n-k}\right ) =  \frac{k-1}{\varepsilon}+O(1)\mbox { for $k\geq 2 $ and $\varepsilon\to 0 $}.$$

Thus, the normal (nearest point) projection  $p_{\pm\varepsilon]}: Y\to Y_{\pm\varepsilon]}$, which is defined for small $\varepsilon\geq 0$, is {\it mean curvature non-decreasing},
$$mn.curv\left (Y_{\pm\varepsilon]}\right)(p_{\pm\varepsilon]}(y))\geq mn.curv(Y)(y)$$
for all  $y\in Y$ and small $\varepsilon>0$.\vspace {1mm}

\hspace {40mm}{ \textbf { From $C^1$ to $C^\infty$.}} \vspace {1mm}

According to $\bigstar$ from  section 11.1,   the $C^1$-hypersurfaces $Y_{\pm\varepsilon]}$
can be $C^1$-approximated by $C^\infty$-smooth ones, call them $Y'_{\pm\varepsilon]}$, such that 

$\ast$   $Y'_{\pm\varepsilon]}$ coincide with $Y_{\pm\varepsilon]}$ outside the $\varepsilon'$-neighbourhood of the singular locus of $Y_{\pm\varepsilon]}$ (where $C^\infty$-pieces of $Y_{\pm\varepsilon]}$ meet) where $0<\varepsilon'<< \varepsilon$ can be taken
arbitrarily small.

$\ast$
The normal (nearest point)  projection $p'_{\pm\varepsilon]}:Y\to Y'_{\pm\varepsilon]}$  is defined for small $\varepsilon$ and it  moves  $Y$ {\it inward}. 

$\ast$ The mean curvature is almost non decreasing under this projection,
$$mn.curv\left (Y'_{\pm\varepsilon]}\right)(p'_{\pm\varepsilon]}(y))\geq mn.curv(Y)(y)+O(\varepsilon').$$

Finally, one can,  if one wishes,    $C^\infty$-approximate 
$Y'_{\pm\varepsilon]}  $  by $Y''$, where mean curvatures  $\geq$ than those at (suitably) corresponding points  of the original $Y$ and  where, moreover,  ({\large $\star_>$}) from the previous section   allows 
one to achieve  a   strict inequality $mn.curv(Y'')> mn.curv(Y) $,   unless (a connected component of)  $Y$ was smooth to start with.

For instance,  if the mean curvature of $Y$ were $\geq 0$  at the regular points of $Y$ then  one can  obtain    
a smooth approximation $Y''$ of $Y$  also with  $mn.curv\geq 0$.

In fact, with a little care, (arguing as in compare 11.6) one can arrange our  $C^\infty $-smooth  $Y'_{\pm\varepsilon]}$ itself,  
such that $$mn.curv\left (Y'_{\pm\varepsilon]}\right)(p'_{\pm\varepsilon]}(y))\geq mn.curv(Y)(y),$$
but this is not essential for the present paper.

\vspace {1mm}

\hspace {30mm}  \textbf {Smoothing non-compact $Y$.} \vspace {1mm}

If $Y$ is a non-compact hypersurface, then   instead of small constant $\varepsilon$ one takes a small and fast decaying function $\varepsilon=\varepsilon(y)>0$. 

A direct construction of  satisfactory  $Y_{\pm\varepsilon]}$ with variable $\varepsilon$,
however trivial, is  cumbersome. A better approach is via a local version of  
$[^\varepsilon\leftrightarrows_\varepsilon]$-{\sl rounding by bending}  procedure described below.\vspace {1mm}

\hspace {14mm}  \textbf {Bending Hypersurfaces Near their Boundaries }

\hspace {27mm} \textbf {and Localisation of Smoothing.} 

Smoothing  the edge, where two  smooth hypersurfaces  $Y_1$  and $Y_2$ in $ X$  meet, can be achieved   
by {\sl inward bending}  of  one of them say of $Y_1$ near its  boundary $\partial Y_1=Y_1\cap Y_2$, such that   

$\bullet$ the bending doesn't decrease the mean curvature of $Y_1$,

 $\bullet$  the bending increases the dihedral angle from the original $\angle (Y_1,Y_2)<\pi$ to $\angle(Y_{1\varepsilon},Y_2)=\pi$,

 where  $\varepsilon>0$, which can be chosen arbitrarily small, signifies that

  $\bullet$ the  bent hypersurface  $Y_{1\varepsilon}$ coincides with $Y_1$ within distance $\geq \varepsilon $ from $\partial Y_1$

and

 $\bullet$ $Y_{1\varepsilon}$  is everywhere   $\varepsilon$-close to $Y_1$.

\vspace {1mm}

The technical  advantage of such a bending is that it is easily   {\it localisable}: you need to bend $Y_1$  only at the points you want to.

Namely, let $Y_1, Y_2 \subset X$ be smooth cooriented hypersurfaces meeting at their common boundaries, 
 $ \partial Y_1= \partial Y_2 =Y_1\cap Y_2$, denote this intersection $Y_{12}$, and  let $\alpha: Y_{12}\to [0,2\pi ]$ be a smooth function, bounded from below by the dihedral angles between $Y_1$ and $Y_2$ at all points in $Y_{12}$,
$$ \alpha(y)\geq \angle_y (Y_1,Y_2), \hspace {1mm} y\in Y_{12}.$$

 Then, for all $\varepsilon>0$, there exists a smooth embedding $\phi=\phi_\varepsilon :Y_1\to X$ with the following properties.

{\sl $\bullet_1$ $mn.curv(\phi(Y_1),y)\geq mn.curv(Y_1, y)$ for all $y\in Y_1$;

$\bullet_2$ $\phi(y)=y$ for all $y\in Y_{12}$;

$\bullet_3$ $\angle_y(\phi(Y_1),Y_2)=\alpha (y)  $ for all  $y\in Y_{12}$;

$\bullet_4$ the map $ \phi$ is $\varepsilon$-close to the original embedding $Y_1\hookrightarrow  X$
(in the $C^0$-topology);

$\bullet_5$  the map $ \phi$ coincides with the embedding  $Y_1\hookrightarrow  X$ at the  points 
$y\in Y_1$ within distance $\geq \varepsilon$ from the the subset $Y_{12}'\subset Y_{12}\subset Y_1$, where 
$\alpha(y')>\angle_{y'} (Y_1,Y_2)$.

$\bullet_6$ The intersection of the  image $\phi(Y_1)\subset  X$ with $Y_2$ is equal to $Y_{12}$.
\vspace {1mm}}

{\it About the Proof.}  This easily accomplished by an isotopy of $Y_1$ which, moreover,  moves the equidistant hypersurfaces $Y_{12, \varepsilon} \subset Y_1$    (of dimensions $n-2$) by at most $\varepsilon$ in the $C^\infty$ topology.

{\it Remarks}.(a)  The condition $\bullet_5$ can be sharpened as follows.

$\bullet_5^\star$   {\sl the (open) subset $Y_{\phi \neq}\subset Y_1$ of the points $y\in Y_1$, where  $ \phi(y)\neq y $, lies   $\varepsilon$-close to $Y_{12}$ and the closure of  $Y_{\phi, \neq}$  intersects  the (closed) subset $Y_{=\angle}\subset Y_{12}$ , where $\alpha(y)=\angle_y(Y_1,Y_2)$,   only at the  boundary of $Y_{=\angle}$ in $Y_{12}$.}

 Achieving  this, which is  unneeded for our applications  anyway, requires a little bit of extra effort.

(b) If one doesn't insist on $\bullet_6$, then one can  bend/rotate $Y_1$ by an arbitrary "angle" 
in the interval $[\angle(Y_1,Y_2 ), \infty)$ in the spirit of  the  {\sl Frizzing Lemma} in [LM 1984].

\subsection  {Rounding Edges of Riemannian Doubles with no Decrease of their    Scalar Curvatures.}
Let $V=(V,g)$ be a smooth  Riemannian $n$-manifold with boundary  and let $W=V\bigcup_{\partial V} V$ be the double of $V$. This $W$ carries a natural continuous Riemannian  metric, call it $\tilde g=g\&g$, which equals $g$ on both  $V$-halves of $W$.

 Let  the boundary $\partial V$ has {\it positive mean curvature} and let us explain following [GL 1980] how \vspace {1mm}

 {\sl $\tilde g$ can be $C^0$-approximated by smooth metrics with
 their scalar curvatures 
 
 bounded from below by $Sc(g)$.}
 
 \vspace {2mm}



\vspace {1mm}
 
{\hspace {40mm}\textbf {$[^\varepsilon\leftrightarrows_\varepsilon]^{\hspace {-1mm}\Rightcircle}$-Rounding.}
\vspace {1mm}

 \hspace {1mm}  {\sl Let   $V_{-\varepsilon} \subset V$ be the complement of the $\varepsilon$-neighbourhood of 
 $\partial V\subset V$}
 
 \hspace {-3mm} and 
 
 \hspace {8mm} {\sl let  $W_\varepsilon \subset V\times\mathbb R $ be
 the boundary of the   $\varepsilon$-neighbourhood of

  \hspace {40mm}$V_{-\varepsilon} =V_{-\varepsilon} \times \{0\} \subset V\times \mathbb R$.}

\vspace {1mm}

This  $W_\varepsilon$ consists of two $\varepsilon$-equidistant copies of $V_\varepsilon$ and of
a semicircular part   $W_{\hspace {-1mm}{\small \Rightcircle}}= \partial V_\varepsilon\times S^1_+(\varepsilon)$, that is  a half of the boundary of  the  $\varepsilon$-neighbourhood of $\partial V_\varepsilon\subset V\times \mathbb R$,
as depicted in  figure 8 on p.  227 in  [GL 1980]).

If $V$ is compact and $\varepsilon \to 0$, then the principal curvatures $\lambda_1,...,\lambda_n$ of $W_{\hspace {-1mm}{\small \Rightcircle}}$ are evaluated in terms of the principal curvatures  $\mu_1,..., \mu_{n-1}$ of the boundary $\partial V \subset  V$as follows.\footnote {We correct here  a minor error from   [GL 1980].}
$$\lambda_i=(\mu_i +O(\varepsilon)) \cdot\cos \theta+o(\varepsilon) \mbox { for  $i=1,..., n-1$ and $\lambda_n=\varepsilon^{-1}+O(1)$,}\footnote{In fact, $\lambda_n=\varepsilon^{-1}+o(1)$, but this is unneeded for our present purpose.}$$
where $-\frac {\pi}{2}\leq \theta \leq  \frac {\pi}{2}$ denotes  the angular parameter of the (right) semicircle  
$S_+^1(\varepsilon)$.

Then the scalar curvature of $W_{\hspace {-1mm}{\small \Rightcircle}}$, which is expressed by the Gauss theorema egregium satisfies
$$Sc(W_{\hspace {-1mm}{\small \Rightcircle}})(v,\theta) = Sc(V)(v)+ (2\varepsilon ^{-1} mn(v)+O(1))\cdot \cos \theta+o(1),$$
where $mn(v)= mn.curv(\partial V)(v)$, where
$v\in \partial V$ and $-\frac {\pi}{2}\leq \theta \leq  \frac {\pi}{2}$.

Now we bring into play the inequality  
$$ mn.curv(\partial (V)>0$$
and see  that  the scalar curvature of $W_{\hspace {-1mm}{\small \Rightcircle}}$ is bounded from below by that of $V$ up to an error $\to 0$. And since
the part of the hypersurface  $W_\varepsilon$, which is parallel to $V$ has the same scalar curvature as $V$, the scalar curvature of  $W_\varepsilon$ is everywhere 
bounded from below, up to an error $\to 0$, by that of $V$. 

\vspace {1mm}
\hspace {40mm}{ \textbf { From $C^1$ to $C^\infty$.}} \vspace {1mm}

 The submanifold  $W_\varepsilon\subset V\times \mathbb R$ has the same type of regularity   as $Y_{\pm\varepsilon]} $ from the previous section, namely, it is {\sl $C^1$ and piecewise $C^\infty$}, and, similarly to  $Y_{\pm\varepsilon]}$, it can be  $C^\infty$-smoothed by applying $\bigstar$  from section 11.1 with a negligible decrease of the scalar curvature, where, moreover, this  smoothing can be performed equivariantly for   the obvious involution of $W=V\bigcup_{\partial V} V.$    (And as in the case of $Y_{\pm\varepsilon]} $, one can, with a little care,  do all this with no decrease of the scalar curvature at all.\footnote{In fact, the required smoothing   can be easily  done directly in this case, as it is indicated in [GL 1980],  where,  in truth -- this was pointed out by the referee  --   it is  only claimed that  $\inf Sc(W)\geq \frac {1}{2}\inf Sc(V)$ and removing "$\frac {1}{2}$" from this inequality {\it without an appeal to}   11.1$\bigstar$ requires   a bit of attention.}

\hspace {-1mm}\textbf {Non-Compact Manifolds, $C^0$-Approximation and $C^\infty$-Smoothing.}
\vspace {2mm}

Similarly to how that was done in the previous section, one extends  the above to {\it non-compact manifolds $V$} with a use of small positive {\it functions} $\varepsilon (v)$, $v\in \partial V$ instead of constant $\varepsilon$.\vspace {1mm}

Besides controlling  scalar curvature of $W_\varepsilon$, we want the metric on $W_\varepsilon$ to be $C^0$-close to $\tilde g$ on $W$ for the former brought to the latter by a suitable     diffeomorphism  $W_\varepsilon  \to W$.  

Firstly, such a diffeomorphism is constructed  from  $W_{\hspace {-1mm}\Rightcircle}$  to the $\frac {\varepsilon \pi}{2}$-neighbourhood  $U\subset W$of   $\partial V\subset W$ with a use of normal  decompositions 
 of $W_{\hspace {-1mm}\Rightcircle}$ and  $U$ as $\partial V\times [-\frac {\varepsilon \pi}{2},\frac {\varepsilon \pi}{2}]$ and then it is extended to all of  $W_\varepsilon$.

\vspace {1mm}




\subsection {"Intrinsic Bending" of Riemannian Manifolds along the Boundaries   with  no Decrease of their Scalar Curvatures.} 

The rounding construction from [GL 1980], which we presented in the previous section, was described in the normal coordinates in [Al 1985] and,  a more general form of it  appears in  [Mia 2002]   in [BMN  2010] and in [McFSzk 2012].\footnote{It was pointed out by the referee of the present paper that   the smoothing proposition   3.1  in  [Mia 2002],  due to  an extra term in there,  doesn't imply the corresponding edge smoothing results from 
[GL 1980] and in  [Al 1985], and, in the report on  the corrected version of this paper, the  referee   pointed out  that what we call "bending"  is proven in [BMN  2010].}
 
Let us  present a   construction of "bending" (Theorem 5 from  [BMN  2010])  following   [Gr 2014*].\vspace{1mm}

{\sf \large $\circlearrowleft_\varepsilon$-Family.}  Let $h$  be a smooth Riemannian metric on a  manifold $Y$, let  $A_{old},  A_{new}$ be smooth quadratic differential forms (i.e. symmetric 2-tensors) on $Y$ and let 
$$h_\varepsilon(t)=h + t A_{new} + \frac {t^2}{2\varepsilon}( A_{old}-A_{new}), \mbox { } 0\leq t\leq \varepsilon.\leqno {(++)}.$$
Then clearly, 
$$h_{\varepsilon}(0)= h\leqno {\mbox{\large $\ast_0$}}$$ 
while
 $$\frac {dh_\varepsilon(0)}{dt}=A_{new}\mbox {    and }\frac {dh_\varepsilon(\varepsilon)}{dt}=A_{old}.\leqno {\mbox{\large $\ast_1$}}$$

Assume at this point that 
$Y$ is compact. Then\vspace {1mm}

\hspace {-5mm}{\large $\ast_2$} \hspace {0mm}  
 the metrics  $h_\varepsilon(t)$  $C^\infty$-converge to $h$ 
for   $\varepsilon \to 0$,
 where this convergence is uniform in $t\in [0,\varepsilon].$
In fact,
$$h_\varepsilon(t)=h +o(\varepsilon),$$
which means that 
 the $C^r$-distances from $h_\varepsilon(t)$ to $h$ satisfy
$$||h_\varepsilon(t)-h||_{C^r}=o(\varepsilon)\mbox { for all } r=0,1,2,...\hspace {1mm} .$$
In particular, the deformation $h\leadsto h_\varepsilon(t)$ is almost circular in $t$ for small $\varepsilon$,
$$h_\varepsilon(\varepsilon) \to h\mbox { for } \varepsilon \to 0.$$

Also,
$$\frac {d^2h_\varepsilon(t)}{dt^2}=\frac {1}{\varepsilon}(A_{old}-A_{new}) \mbox { for all } (y,t)\in Y\times[0,\varepsilon].\leqno {\mbox{\large $\ast_3$}}$$



  Now, let us incorporate the family $h_\varepsilon(t)$ into family of Riemannian metrics $$g_{[0, \varepsilon]} = h_\varepsilon(t)+dt^2$$ 
   on the manifold $Y\times [0,\varepsilon]$.
   
   Then, clearly,
  
  \vspace {1mm}
    
 \hspace {-5mm}\textbf {[i]} \hspace {1mm}  {\sl the 
 second fundamental forms of the two boundary parts  $Y\times\{0\}$ and $Y\times\{\varepsilon\}$ in the manifold  $Y\times [0,\varepsilon]$ with  the metric $ g_\varepsilon$  are equal to $A_{new }$ and $A_{old}$  correspondingly, where these forms are  evaluated on the (same unit) vector field $\frac {d}{dt}$.}
   
   Furthermore, the scalar curvature of this metric satisfies,
   $$ Sc(g_{[0, \varepsilon]})(y,t)=-\frac {1}{\varepsilon}  trace(A_{old}-A_{new}) +O(1).\leqno {\mbox {\textbf {[ii]}} }$$

Indeed, the Ricci curvature $Ricci_\varepsilon$   of $g_{[0,\varepsilon]}$ at the field $\frac {d}{dt}$ is expressed in terms of 
$A_\varepsilon (t)= \frac {dh_ \varepsilon (t)}{dt}$ by Hermann Weyl's formula
$$Ricci_\varepsilon  \left (\frac {d}{dt},\frac{d}{dt}\right)=-trace\left( \frac {d}{dt} A_ \varepsilon(t)^\ast+ (A_ \varepsilon(t)^\ast)^2\right)$$ 
where $A_ \varepsilon^\ast(t)$ denote the (selfadjoint shape) operators associated to the quadratic forms $A_ \varepsilon(t)$ via the metrics $h_\varepsilon(t)$ on $Y$. 
 (Recall that  the eigenvalues of $A_ \varepsilon^\ast(t)$ are equalto the principal curvatures of the hypersurface $Y\times \{t\}\subset Y\times [0,\varepsilon]$ for the metric $g_{[0,\varepsilon]}$ in $Y\times [0,\varepsilon]$.)

Therefore, in view of  {\large $\ast_2$} and {\large $\ast_3$},
$$Ricci_ \varepsilon \left (\frac {d}{dt},\frac{d}{dt}\right)=-trace (A_{old}-A_{new}) +O(1)$$
and {\textbf {[ii]}}  follows by the Gauss theorema egregium for the hypersurfaces  $Y\times \{t\}$  in the manifold  $Y\times [0,\varepsilon]$ with the  metric $g_\varepsilon$.

\vspace{1mm}

{\sf \large Bending Lemma.} Let $X=(X,g)$ be a smooth  (possibly non-compact) Riemannian manifold with boundary $Y=\partial X$ and  let $\overrightarrow A_{old}$  denote the second fundamental form of $Y$ with respect to the {\it outward} normal field. (Notice that the boundaries  of {\it convex} domains have {\it positive  definite} second forms for  such fields.)  

Let $\overrightarrow A_{new}$ be another smooth quadratic form on $Y$. 
\vspace {1mm}

{\it If $$trace(\overrightarrow A_{new})< trace (\overrightarrow A_{old})$$
then there exits a family of $C^\infty$-smooth metrics $g_\epsilon =g_{new, \epsilon}$, $\epsilon>0$,  on $X$, such that}\vspace {1mm} 
 
[I] {\sl The restrictions of the Riemannian  metrics $g_\epsilon$ to $Y=\partial X\subset X$ are equal to such restrictions of $g$  for all $\epsilon>0$.} 

[II]  {\sl The second fundamental forms of $Y\in X$ with respect to $g_{new, \epsilon}=g_\epsilon$,  are equal to $\overrightarrow A_{new}$  for all $\epsilon>0$.}

[III] {\sl The scalar curvatures of $g_\epsilon$ are bounded  from below by those of $g$, 
$$ Sc( g_\epsilon)\geq Sc(g).$$}

[IV] {\sl The metrics $g_\epsilon$  $C^0$-converge to $g$ for $\epsilon\to 0$.}

[V]  {\sl The metrics $g_\epsilon$ are equal to $g$ within distance $\geq \epsilon$  from $Y$.}

\vspace {1mm}

{\it Proof} Let $\varepsilon <<\epsilon$ be a (very small) positive number,  such that the  the $\varepsilon$-neighbourhood $U_\varepsilon$ of $Y$  
{\it normally splits},
$$U_\varepsilon = Y\times [0,\varepsilon],$$
where, by the definition of "normally",    the hypersurfaces $Y\times \{t\}$ are $t$-equidistant to $Y=Y\times \{0\}\in X$.
 
Let us replace the metric $g$ in $U_\varepsilon $ by the above $g_{[0, \varepsilon]}$ with $A_{old}=-\overrightarrow A_{old}$ and  $A_{new}=-\overrightarrow A_{new}$, where 
  the sign reversion is due to the inward direction of the field $\frac {d}{dt}$. 

Also observe that  the Riemannian form $g_\varepsilon$ on the hypersurface  $Y\times\{\varepsilon\}$ as well as second fundamental form of $Y\times\{\varepsilon\}$ with respect to  $g_\varepsilon$ are $\varepsilon$-close to the  $g$-related forms.  Therefore, according to  $\bigstar '$ from section 11.1  applied to sections of the bundle of  quadratic forms on the tangents to submanifolds $Y\times \{t\}\in X$, the form $g_\varepsilon$ extends to all of $X$ with the required properties except [III]; this can be guaranteed by $\bigstar '$ only up to an arbitrarily small error. Yet, this error can be taken care of by {\sf redistribution of the scalar curvature} from section 11.2.

(If our  "redistribution"   is obtained  with a  use of the  Fredholm + unique continuation properties  of the operator 
$f\mapsto Sc(f^2g)$,  then   $V$ needs  to be compact and also  
the condition  [V] suffer .  But none of this   happens if we rely on the  Local Sc-Lemma.\footnote{There is no circularity here, since  
 Local Sc-Lemma, albeit  being a special case of the   Bending Lemma, admits an independent (and obvious) proof.})

\vspace {1mm} 

{\sf \large Half Way  
 {$[^\varepsilon\leftrightarrows_\varepsilon]^{\hspace {-1mm}\Rightcircle}$-Rounding} and Making Doubles with $Sc\geq \sigma$.} If $X=\partial X$ has strictly  positive mean curvature, then the above applies to $A_{new}=0$ and allows  a "bending" (of the Riemannian  metric in)  $X$ near the boundary $Y=\partial X$,  which makes $Y$ {\it totally geodesic}  and such that   \vspace {1mm}

{\it the  scalar curvature of  $X_{bent}= (X, g_{new})$ is bounded from below by the  scalar 

curvature 
of the original metric} $g$. 

Then the  metric $ g_{new}$\&$g_{new}$ on the  double $(X, g_{new})\cup_Y(X, g_{new})$
is $C^1$-smooth and, by $\bigstar$, it can be $C^\infty$-smoothed keeping its scalar  curvature bounded from below
by $Sc( g_{new})$ which   itself is bounded from below by  $ S(g)$ according to the above {\sf Bending Lemma}.

\vspace {1mm}

\textbf {VI}$_{loc}$  {\it \large Localisation of Bending.} Let   the form  $A_{new}$  be equal to $A_{old}$ on a (now possibly non-empty) compact subset  $Y_0\subset Y=\partial X$ and the inequality   $trace(A_{new})< trace (A_{old})$ holds in  the complement $Y\setminus Y_0$. 

{\sl If  $X$ is connected and the complement $Y\setminus Y_0$ is non-empty, then  
  there exists a family $g^\star_\varepsilon$, which satisfies the above [I]-[IV] and where  [V] is strengthened to the following }
\vspace {1mm}

[V$_{\star}$] {\sl The metrics $g_{\varepsilon}$  are equal to $g$  within distance $\geq \varepsilon$  from $Y\setminus Y_0$.} 
  \vspace{1mm}

{\it About the Proof.} Construction of  $g^\star _{\varepsilon}$, which satisfies this condition (parallel to the condition $\bullet _5$ in 11.3 for bending  hypersurfaces with controlled mean curvatures) is achieved by an obvious smooth cut off of the above  $g _{\varepsilon}$, where the arising error term in $Sc(g^\star_{\varepsilon})$ can be compensated due to the  $\frac {1}{\varepsilon}$-contribution  in the positivity  of the scalar curvature.

\subsection {Minimal Hypersurfaces in Non-compact Manifolds.} Let us return to  $M$-bubbles from section 9:  these are   closed cooriented hypersurfaces in Riemannian manifolds, say  $Y_{min}\subset X$,  which  locally minimise the functional 
$$Y\mapsto vol_{n-1}(Y)-\int_ {U_-} M(x)dx,$$
where $M$ is a function defined in a neighbourhood of $U_\ast\supset Y_{min}$ and   where  $U_-
\subset U$ is the interior part of $U_\ast$.

Unconditionally, {\it the existence} of such a   $Y_{min}$, possibly a singular one if $n=dim(X)\geq 8$, is ensured 
only for {\it closed} manifolds $X$.

And   if $X$ is a compact manifold with boundary, then the existence of $Y_{min}$ follows from a  suitable bound on the mean curvature of the boundary  as in section 9, 
$$\mbox { $M_{|\partial _-X}\leq -mean.curv (\partial _-X)$ and $M_{|\partial _+X}\geq mean.curv (\partial _+X)$}$$
 where $\partial_\pm$ are the parts of the boundary of $X$ positioned in the interior/exterior  region of $X$  with respect to $Y_{min}$. (This makes sense since the decomposition $\partial X=\partial_-X\cup \partial_-X$
depends only on the homology class of $Y_{\min}$.)

Now we look for  this kind of condition for non-compact manifolds, where our motivations are twofold: \vspace {1mm}

(A)  {\sf Finding geometric conditions on   $X$, such that a given hypersurface $Y\subset X$ would  admit a  metric  with positive scalar curvature.}
\vspace {1mm}

(B)  {\sf Finding  least demanding  constraints on   a perturbation $g_\varepsilon$ of a  Riemannian metric $g$ on $X$,
such that a  hypersurface  $Y\subset X$ would admit an $\varepsilon$-close to it hypersurface with the mean curvature with respect to  $g_\varepsilon$ being close to $mn.curv_g(Y)$.}\vspace {1mm}

An obvious  (over-optimistic?) {\sf \large conjecture}  in the direction of  (A) is as follows.\vspace {1mm}

({\large\sf A \textbf {?})\hspace{-0mm}} {\sf Let $X$  be a complete Riemannian manifold of dimension $\geq 6$ and let
$M(x)$ be  a continuous function on $X$, such that 
$$ \frac {n}{n-1}M(x)^2 -2||dM(x)|| +Sc(X)(x)\geq 0.$$
Then every closed cooriented  hypersurface $Y_0\subset X$, for which  the inclusion homomorphism between the fundamental groups
$$ \pi_1(Y_0)\to \pi_1(X)$$
is injective,  admits a a metric with $Sc>0$.}

(This generalises Conjecture 1.24 in [Ros 2007] for $X=Y_0\times S^1$ and $M=0$ and a similar conjecture
 -- I recall seeing it in a paper by Rosenberg and/or Stolz --
for complete $X$ homeomorphic to $Y_0\times \mathbb R$.) 

\vspace{1mm}

{\it Motivating Example.} Let  let $X$ be  a {\it complete  two-ended} Riemannian manifold 
and let $Y_0$ be a smooth  closed   hypersurface, such that the following four conditions are satisfied.
\vspace{1mm}

[\hspace {1mm}{\large \textbf {\textdivorced}}\hspace {1mm}] \hspace {0.6mm}$Y_0$ is {\it connected} and it 
{\it separates the two ends of} $X$. Thus,  both components of the complement $X\setminus Y_0$ are infinite and both are one-ended. For instance,   $X$ is homeomorphic to $Y_0\times \mathbb R$. 

 (We agree that the only  "infinities of" $X\setminus Y_0$ which count   are parts of the   "infinity of" $X$.)

\vspace {1mm}

[{\Large $\circ$=$\circ$}]$_{H_1}$ The inclusion homology homomorphism  $H_1(Y_0)\to H_1(X)$ is an {\it isomorphism.} (If  the inclusion homomorphism between the fundamental groups 
$\pi_1(Y_0)\to \pi_1(X)$ is injective then the covering of $X$ with the fundamental group equal to the image $\pi_1(Y_0)\subset \pi_1(X)$ has this property.)  
\vspace {1mm}

[{\Large $\circ$$\hookrightarrow$$\circ$}]$_{\pi_1}$ The inclusion homomorphism between the fundamental groups 
$\pi_1(Y_0)\to \pi_1(X)$ is injective. 

(This implies [{\Large $\circ$=$\circ$}]$_{H_1}$ for the covering of $X$ with   fundamental group equal the image $\pi_1(Y_0)\subset \pi_1(X)$, but the roles played by  $H_1$
 and $\pi_1$ in the arguments below are different.)


 \vspace {1mm}

$[vol=\infty]$ Both connected components of the complement $X\setminus Y_0$, call them $X_\pm\subset X$, have {\it infinite} volumes.

 \vspace {1mm}

{[\sf spin]} Manifold $X$ is {\it spin}.

 \vspace {1mm}

$[\odot_{\lambda,\rho}]$ There exist constants $ \lambda,\rho>0$, such that   all balls in $X$ of radii $r \leq \rho$ are {\it $\lambda$-Lipschitz contractible in $X$} (but not necessarily within themselves), i.e. there exist  $\lambda$-Lipschitz  maps
$$ \phi= \phi_{x,r}: B_x(r)\times [0, r]\to X, \hspace {0.7mm} x\in X, \hspace {0.7mm}r\in [0, \rho],$$
where  the maps  $\phi( ..., 0)$  are  the original imbeddings $B_x(r)\subset X$   and   such that the maps $\phi(..., r): B_x(r)\to X$  {\it are constant.}

(Coverings of closed manifolds, obviously, satisfy this condition.)

\vspace {1mm}

[{\Large $\bullet $}]  {\sf If the scalar curvature of $X$ is everywhere  bounded from below by $\sigma_-\in (-\infty, +\infty)$ 
and $Sc(X)(x)\geq \sigma_+>0$ in the $\varepsilon$-neighbourhood of $Y_0\subset X$, such that 
$$ \sigma_+ \geq \frac {2}{\varepsilon} \sqrt  {\frac {|\sigma_-|(n-1)}{n}}.\leqno {[ \sigma_+>>|\sigma_-|]}$$
Then $Y_0$ admits a metric with $Sc>0$, provided $6 \leq dim(X)\leq 8$. }

\vspace {1mm}

Prior to explaining the proof, a few remarks  are in order.\vspace {1mm}
{\Large $(\ast)$} The spin condition  can be significantly relaxed and, possibly, fully removed. (We shall explain this in the course of the proof of [{\Large $\bullet $}].)

{\Large $(\ast) $} The condition $[\odot_{\lambda,\rho}]$ implies $[vol=\infty]$ as we shall see in the course of the proof of  [{\Large $\bullet $}]; however, the two play  opposite roles in the proof of [{\Large $\bullet $}]. But in any case, we would rather get rid of $[\odot_{\lambda,\rho}]$ altogether. \vspace {1mm}

{\Large $(\ast)$} The main function of the inequality $dim(X)\geq 6$  is to rule out $4$-manifolds $Y_0$, where there are topological  obstructions  to the existence of metrics with $Sc>0$ which have no counterparts for other dimensions, see [Ros 2007]   and references therein.

{\Large $(\ast )$} The inequality $ dim(X)\leq 8$ is, most likely, unnecessary -- it is  due to our inability to handle singularities of minimal hypersurfaces in manifolds of dimensions $\geq 9$.

\vspace {1mm}

{\it Proof of }[{\Large $\bullet $}].  {\large \sf Step 1}. The inequality $[ \sigma_+>>|\sigma_-|]$ implies that there exists 
a function $M(x) $, which vanishes on $Y_0$ which is positive on one of the components of $X\setminus Y_0$, say on $X_+$  and negative on $X_-$,   which is constant $2\varepsilon$-far from $Y_0$ and such that   
$$ \frac {n}{n-1}M(x)^2 -2||dM(x)|| +Sc(X)(x)> 0$$
at all points $x\in X$. 
\vspace {1mm}

{\large \sf Step 2.} Since $vol(X_\pm)=\infty$, the integrals of $M$ over both components $X_\pm$ of $X\setminus Y_0$   are $\pm$infinite.
Therefore, there exists a {\it compact connected} subset $V_0\subset  X$, which contains $Y_0$ and such that 
$$\int_{V_0\cap X_\pm}  |M(x)|dx> vol_{n-1}(Y_0).$$
Consequently, every closed hypersurface $Y'$ in the complement $X\setminus V_0$,  which  is  homologous to $Y_0$ satisfies  
$$vol_{n-1}(Y')- \int_{X'_+ \setminus \infty_+ }M(x)dx \geq vol_{n-1}(Y_0)-\int_{X_+ \setminus \infty_+}M(x)dx+\mu_0,  $$
where 

$\bullet $  $\infty_+$ denotes   a unspecified   subset  in $X_+ $, which has a sufficiently large  compact complement; in particular, $\infty_+$ doesn't intersect $V_0$,

$\bullet $  $X'_+ \subset X$ is the (infinite)  connected component of $X\setminus Y'$ which contains $\infty_+$.

$\bullet $ $\mu_0=\mu_0(V_0)$ is a positive constant. \vspace{1mm}

Since the difference   $$\int_{X'_+ \setminus \infty_+ }M(x)dx-\int_{X_+ \setminus \infty_+}M(x)dx  $$
does not depend on the choice of  $\infty_+\subset X_+$,    the above inequality is unambiguous.

It follows that

\hspace {32mm}{\sc   $Y_{min} $  can't escape from $V_0$}:

\vspace{1mm}

{\sl all closed hypersurfaces $Y'\subset X $, which are   homologous to $Y_0$ and which almost  (up to $\mu_0$) minimize the functional  
$$Y'\mapsto  vol_{n-1}(Y')- \int_{X'_+ }M(x)dx=_{def}vol_{n-1}(Y')- \int_{X'_+ \setminus \infty_+ }M(x)dx$$
{\it intersect} $V_0$.}

And since these $Y'$, because of the above [\hspace {1mm}{\large \textbf {\textdivorced}}\hspace {1mm}] \hspace {0.6mm}$Y_0$ and 
[{\Large $\circ$=$\circ$}]$_{H_1}$, are connected, the minimization process for this functional converges (in the sense of the geometric measure theory) to a 
 minimum, call it  $Y_{min}\subset X$: this is   a possibly infinite hypersurface in $X$  of {\it finite volume}, 
where this $Y_{min}$ {\it doesn't fully escape} $V_0$. 
Namely

  \vspace {1mm}

$\bullet $  $Y_{min}$ is 
{\it proper}, i.e. it is a  closed  as a subset in 
$X$, but it {\it may be non-compact}  (hence,  unbounded for complete  $X$);

  $$vol_{n-1}(Y_{min})\leq vol(Y_0)<\infty \leqno {\hspace {5.5mm}\bullet}$$

$\bullet$ The hypersurface $Y_{min}$ has {\it non-empty} intersection with $V_0$. 

(Notice that completeness of $X$ is non-essential at this point.) 
\vspace {1mm}

{\large \sf Step 3.} Let us  bring forth   the above  local Lipschitz  contractibility condition $[\odot_{\lambda,\rho}]$, 
 invoke  the {\it cone/filling  inequality} 3.4.C from [Gr 1983]. This shows that $[\odot_{\lambda,\rho}]$
 implies the following {\it filling inequalities $[\odot^{n-1}_{n-2}]$ and $[\odot^{n}_{n-1}]$}. \vspace {1mm}
 
$[\odot^{n-1}_{n-2}]$ {\sl All  $(n-2)$-cycles $\beta$ in all $\rho$-balls,    $\beta =\beta_{n-2} \subset B_x(\rho)\subset  X$, $n=dim(X)$,  bound $(n-1)$-chains $\alpha\subset X $, such that 
$$vol_{n-1}(\alpha)\leq c\cdot vol_{n-2}(\beta)^\frac {n-1}{n-2},$$
for some   constant $c=c(n, \lambda, \rho)$.

 $[\odot^{n}_{n-1}]$ All subdomains $U\subset B_x(\rho)$ satisfy
 $$vol_{n}(U)\leq c'\cdot vol_{n-1}(\partial U)^\frac {n}{n-1},$$
 for some    $c'=c'(n, \lambda, \rho)$.}
\vspace {1mm}
 
 Then these inequalities easily   yield the following   lower bound on the volume of the  balls in the
 minimizer $Y_{min}\subset X$,
 $$vol_{n-1} (Y_{min}\cap B_x(\rho))\geq \nu=\nu(n, c,c')=\nu(n,\lambda, \rho)>0$$
 for all $x\in Y_{min}$.\footnote{This is explained   in [Gr 1983], where  the  corresponding inequalities are formulated for the functional $Y\mapsto vol_{n-1}(Y)$  and  where  $X\supset Y$ is often  required to be compact.  But all (relevant) arguments from [Gr 1983] apply  in the present case.} 

Since   $Y_{min}$ is connected, the   above volume bound  $vol_{n-1}(Y_{min})\leq vol(Y_0)$ implies a bound on its diameter,
$$diam(Y_{min})\leq R=\frac { 2\rho\cdot  vol_{n-1}(Y_0)}{\nu }<\infty.$$
Then, by combining this inequality with the above intersection property for our minimizer,  $Y_{min}\cap V_0\neq \emptyset $,
we conclude that
\vspace {1mm}

 \hspace {7mm}{\sc  $Y_{min}\subset X$ is trapped\footnote {This means "contained".}  in the $R$-neighbourhood of $V_0$.}
\vspace {1mm}

 \hspace {-6mm} And since $X$ is complete, $Y_{min}$ is {\it compact.}
\vspace {1mm}

{\large \sf Step 4.} By the classical  regularity theorem(s) of  Simons-Federer-Almgren-Allard this $Y_{min}$ is a smooth hypersurface for $dim(X)\leq 7$;  if $n=8$,  then the regularity is achieved by a small perturbation of the Riemannian metric in $X$.\footnote {The argument from [NS 1993)] easily generalises  
to our $Y_{min}$. }

Then,  the inequality  $ \frac {n}{n-1}M(x)^2 -2||dM(x)|| +Sc(X)(x)> 0$, implies 
that the induced metric in $\tilde Y_{min}$ is conformal to a metric with $Sc> 0$  by the  $M$-version of the Schoen-Yau argument (see \S $5\frac {5}{6}$ in [Gr 1995]). \vspace {1mm}

{\large \sf Step 5.}  If $Y_{min}$ is smooth and  $dim(Y_{min})\geq 5$ 
then the kernel of the inclusion homomorphism 
$$\pi_1(Y_{min})\to\pi_1( X)=\pi_1(Y)$$
can be "killed" by 2-dimensional surgery that results in  another  hypersurface, say $Y\subset  X$, which is also homologous to $Y_0$ and which admits a metric with $Sc\geq 0.$

If, moreover,  the inclusion homomorphism $\pi_1(Y_0) \to \pi_1(X)$ is {\it an isomorphism,}  and $X$ is spin,  then $Y$ is spin-bordant to $Y_0$ in the classifying space of $\pi_1(Y_0)$ and the existence of metric with $Sc>0$ on $Y_0$, follows  from the  theorem 1.5 in [St 2001]. (Possibly -- I am not certain --  this theorem   also  covers the non-spin case.)

Finally,  since the homomorphism $\pi_1(Y_0) \to \pi_1(X)$  is injective according to   
[{\Large $\circ$$\hookrightarrow$$\circ$}]$_{\pi_1}$, all of the above  applies to the covering space   $\tilde X\to X$ with the fundamental group $\pi_1(\tilde X)-\pi)Y_0)$.
QED.

\vspace {1mm}

{\sc Discussion.}  The above argument can be  generalised, refined and made  effective, which would 
result in specific {\sl inequalities for the "relative size"} of  pairs $(X, M(x))$, where  $Sc(X)(x) $ is suitably  bounded from below in terms of $M(x)$ and  where certain hypersurfaces  $Y_0\subset X $ admit {\it no metrics} with $Sc>0$.

However, it remains unclear if the condition $[\odot_{\lambda,\rho}]$, 
or anything of this kind, is truly necessary. 

In fact, we do know that {\it no such condition is needed} for the bounds on the width   of  overtorical and similar band-shaped manifolds to which   
our symmetrization with point-wise 
 control of the scalar curvature  applies (see  sections 7 and 8).

Also, even the present  form of [{\Large $\bullet $}] remains problematic for manifolds $X$ with 
$dim(X)\geq 9$.

\subsection {Bounds on Widths of  Non-Compact Riemannian  Bands.}

The  band inequalities in sections 2, 4, 5 generalise to certain {\it non-compact complete} bands.

{\sc \textbf {Case I}} .  Let $V$ be a  {\it proper} (see section 2)  band  with {\it compact boundaries}  
$$\partial V=\partial_- \cup \partial _+$$

The concepts of "overtorical", "isoenlargeable", and "SYS" make sense for these $V$, where the instances  of such non-compact bands are, topologically, obtained from  compact ones by removing isolated  points from their interiors.

\vspace {1mm}

 {\sl proper   bands $V$  with compact boundaries which have $Sc(V)\geq \sigma>0$ satisfy the same width bound  as their compact counterparts in sections   2, 4, 5.}

\vspace {1mm}

{\it Proof.} Let $V'\subset V$ be a compact manifold obtained by cutting of the infinity of $V$ far away from $\partial V$.  Thus $V'$ has an extra set of components, call their union $\partial'$, where one can make  make 
 the distance $dist(\partial(V),\partial')$ as large as one wants.

Thus we  arrange such a $V'$,  where the the $\rho$-neighbourhood $U'_\rho\subset V'$ of $\partial'\subset V'$ for  a large $\rho$ does not disturb the essential topology of  $V$:\vspace {1mm}
 
 [{\huge $\star_\rho$}] the subset 
 {\sl $U'_\rho \subset V'$ does not intersect $\partial=\partial V$ and, moreover,  $V'\setminus  U'_\rho$ is in the same topological largeness  class as $V$, namely  overtorical, "isoenlargeable  or SYS  correspondingly.}

Clearly,  this property is inherited by minimal hypersurfaces constructed in the suitable homology classes as in  section 2, 4, 5, where the final stage of the symmetrization goes through if $\rho$ is sufficiently large.

For instance, $\rho \geq \frac  {4\pi}{\sqrt \sigma} $ is sufficient for this purpose but this bound doesn't seem  sharp.
 QED

\vspace {1mm}



\vspace {2mm}

{\sc \textbf {Case II}}. Let $V$ be a proper complete orientable  $n$-dimensional   band, such that 
 $V$ admits a proper 1-Lipschitz map $\psi: V\to \mathbb R^{n-1}$, such that the restriction of $\psi$
to $\partial_-$}  ({\sl and hence, on $\partial_+$})   has degree $d\neq 0$.

\vspace {1mm}

{\sl If the scalar curvature of $V$ is bounded from below by $Sc(V)\geq \sigma>0$, then   the width of $V$ is bounded as follows.
$$width(V)=dist (\partial_- , \partial _+)\leq 2\pi\sqrt{\frac { n-1}{\sigma n}}.$$}

{\it Proof.} Let $S^{n-3}(R)\subset \mathbb R^{n-1}$ be a very large codimension 2  sphere  say of  radius $\geq 
100^{n-1}$ and let  
$$ \Sigma =\partial V\cap  \psi^{-1}( S^{n-3}(R))\subset V,$$ where one may   assume if one wishes -- this is not truly necessary -- that  the map $\psi$ is  transversal to the sphere $S^{n-3}(R)$.

This  $ \Sigma$, which has codimension 2 in $V$, is meant to serve  as the boundary condition for a minimizing hypersurface, $Y$,
namely,  the boundary of $Y$ must be contained in $ \Sigma$ and the relative homology class of 
$Y$ in $H_{n-1}(V, \Sigma)$ must be equal to the homology pullback of the $(n-1)$-ball bounded by this sphere.


To avoid unnecessary non-compactness problems,\footnote{There is no such problem   
if $V$  has {\it locally bounded geometry}, e.g. if it is a covering space of a compact band, then the minimizing $Y$ is compact since locally bounded geometry implies uniform local Lipschitz contractibility ($[\odot_{\lambda,\rho}]$ from the previous section.} we cut off $V$
and thus $Y$, by taking a large $V'\subset V$  which contains  $\psi^{-1}( S^{n-3}(R))$, say the $\psi$-pullback of the  $R'$-neighbourhood $U=U(R,R')\subset  \mathbb R^{n-1}$ of  the sphere  $S^{n-3}(R)\subset \mathbb R^{n-1}$ for 
 $R'\geq 100 R$.

\vspace {1mm}

 This $Y$   serves as the   first step of the  inductive   symmetrization, except that unlike the original    
$V$ it is compact and has extra part of the boundary, namely, the $\psi$-pullback of the boundary of $U$.

Thus, in order to have a proper inductive scheme , one should drop the completeness assumption on $V$,  and require instead that 

{\sl $V$ admits a proper 1-Lipschitz map to the $R_n$-ball  in $\mathbb R^{n-1}$ for a sufficiently large $R_n$, say for $R_n\geq 1000^n$, such that the restriction of $\psi$
to $\partial_-$  has degree $\neq 0$.}

At this point, we invite the reader to fill in the details in the above argument.
\vspace {1mm}

{\sc  Corollary}: \textbf {Sharp Bound on Widths of  Iso-Enlargeable Bands}.

  {\sl \textbf Compact Iso-enlargeable bands $V$ (e.g., those homeomorphic to $V_0\times [-1,1]$
  where $V_0$ admits a metric with non-positive curvature) with their scalar curvatures $\geq \sigma$
  satisfy the 
    $\frac {2\pi}{n}$-Inequality.
     $$width(V)\leq    2\pi \sqrt{\frac {n-1}{\sigma n}}.$$}

{\it Proof.} It follows from the definition of iso-enlargeability of $V$ that there exist Riemannian bands $V_R$
for all $R>0$, such that

$\bullet_1$ there exist locally isometric band maps $V_R\to V$,

$\bullet_2$ there exist  proper 1-Lipschitz maps from $V_R$  to the $R$-balls  in $\mathbb R^{n-1}$,  such that the restrictions of $\psi$
to $\partial_\pm$  have (equal) degrees $\neq 0$.

Thus, the above applies and the proof follows.\vspace {1mm}

{\it Remark.} Minimal hypersurfaces with mixed
 boundary conditions in our argument are similar to those in   
\S 12 in [GL 1983], where, in fact, a non-sharp version of the above inequality is proven.

\vspace {1mm}

{\it On the Definition of Iso-enlargeability.} These $\bullet_1$ and $\bullet_2$ can be taken for the definition of iso-enlargeability of bands as it was mentioned in footnote 9 in section 4.

\subsection {Stable  Enlargeability and Stable Bounds on the Scalar Curvature.}

 Let $X$  be an  orientable Riemannian manifold of dimension $n$ and  $P$ an orientable pseudomanifold  of dimension $N$.

  Let us consider three different cases   indexed by $i=0,1,2$, each   associated with     a map
$$f =f_i:  X\times P\to B\subset \mathbb R^{n+N-i},  \hspace {1mm} i=0,1,2,$$
where $B$ is an open unit ball,  and where $f$ is  a smooth proper, such that   
\vspace {1mm}

{\sl  the restrictions of $f$ to the the  submanifolds $X_p=X\times \{p\}\subset P$, $p\in P$, 
are $\lambda$-Lipschitz for some constant $\lambda>0$} 
and  \vspace {1mm}

$\bullet_0$  if $i=0$, then  the map $f$ has {\it non-zero degree};

$\bullet_1$  if  $i=1$,  then the cap product of the pullback of the fundamental cohomology class  of $B$ with compact support with the fundamental class of $X\times P$,
$$f^\ast [B]^{n+N-1}\cap [ X\times P]_{N+n}\in H_1( X\times P),$$
{\it does not vanish} when taken with  $\mathbb Q$-coefficients, i.e. under     the homomorphism $H_1( X\times P)\to H_1( X\times P);$

$\bullet_2$  if $i=2$ then the cap product class\footnote{Geometrically,  this class is represented by (possibly singular) surfaces in  $X\times P$, which are the  $f$-pullbacks of  generic points $b\in B$.}  
$$f^\ast [B]^{n+N-2}\cap [ X\times P]_{N+n}\in H_2( X\times P)$$
is {\it aspherical}, i.e. it  is not contained in the image of the the Hurewicz homomorphism, or, equivalently, it does not lift to the universal covering of $X\times P$.

\vspace {1mm}

{\sl Then the infimum of the scalar curvature of $X$ is bounded by
$$\underline\kappa_X=\inf _{x\in X}Sc(X)(x)\leq \varepsilon \lambda^{-2}$$
where $\varepsilon>0$ depends only on $n$ and $N$.}\vspace {1mm}

{\it Proof}. First, let  $P$ be a manifold and let us endow it   with a very (arbitrarily) large metric which has  scalar curvature $\geq -\delta$ for  a given arbitrarily  small $\delta>0$. This makes  $Sc(X\times P)\geq \kappa-\delta$ and, at the same time,   the Lipschitz constant of $f$ as close to $\lambda$ as you wish.

 Then, the proof follows either by  arguing as in the proof of the spherical Lipschitz bound theorem in section 3  with a torical band in $B$ and a use of   width inequalities for compact bands  or doing it  more directly with open bands as in the previous section. (This directly  applies to $i=0,2$ and the case $i=1$ trivially reduces to $i=0$.)

In general, if  $P$ is a pseudomanifold, we take a similarly large  metric in $P$, where $Sc\geq -\delta$ on all  $N$-faces $Q$ of $P$.  Then, clearly, there exist \vspace {1mm}

(*) a face $Q$,

  (**) a (large) open subset $X'\subset X$,
  
  (***)  a  
(small) open   subball  $ B'=B'_{Q, X'}\subset B$ around some point $b\in B$,  \vspace {1mm}

\hspace {-6mm}such that the above applies to the restriction $f'$ of the  map $f$
to some the subproduct   $ X'\times Q'\subset f^{-1}(B')\cap (X'\times Q)$ for some open subset $Q'\subset Q$
$$f': X'\times Q'\to B'.$$

 \vspace {1mm}

{\it Complaint.} One can't help but to be annoyed     by the  {\it the dimension  being brought  up}  by   incorporating   {\it purely topological}  parameters $P$ into the {\it geometry} of $X\times P$, 
 only to be  immediately    {\it  brought   down} by constructing a decreasing   chain of minimal hypersurfaces.\footnote{If $X$ is {\it complete} and the   universal covering $\tilde X$ of $X$ is {\it spin}, then the bound  $\inf _{x\in X}Sc(X)(x)\leq \varepsilon \lambda^{-2}$ for $i=0,1$ (but not for $i=2$)  follows by applying a suitable index theorem for  the $P$-family of Dirac operators on $\tilde X$,  keeping the $P$-parameters at their proper place.}




\subsection {Mean Curvature  Stability of Polyhedral Domains.} Let $Y$ be a closed smooth cooriented hypersurface in a Riemannian manifold $X=(X,g)$ and let $g_\varepsilon $ be a family of smooth Riemannian metrics on $X$ which $C^0$-converge  to $g$ for $ \varepsilon\to 0$. It is shown in section 10.2 of  [Gr 2012]  that such small  perturbations $g_\varepsilon$ can be accompanied by small perturbations of $Y$ which only slightly change the mean curvature of $Y$. Namely we have the following perturbation stability property

\vspace {1mm}

[{\small   $\Circle_{\varepsilon} $}]. {\sl there exists a family of diffeomorphisms  $\psi_\varepsilon :X \to X$, $\varepsilon>0$, such that

(1)  the diffeomorphisms  $\psi_\varepsilon$, $ \varepsilon\to 0$, converge to the identity map $id: X\to X$ in the $C^0$-topology;

(2) the $(n-1)$-volumes of the hypersurfaces $Y_\varepsilon=\psi_\varepsilon(Y)\subset X$ with respect to $g_\varepsilon$ converge to 
$vol_{n-1}(Y)$  for $g$;

(3) the $g_\varepsilon$-mean curvatures of $Y_\varepsilon$ converge to $mn.curv_g(Y)$.

For instance, if $mn.curv(Y)> c$  for some number $c\in (-\infty,+\infty)$, then the hypersurfaces $Y_\varepsilon$ satisfy the same inequality for all sufficiently small $\varepsilon>0$.}

\vspace {1mm} 

The essential steps in  the   proof of this  can be seen as $\varepsilon$-miniaturised  versions of  the first three steps in the proof of     [{\Large $\bullet $}] in section 11.6, which now   take place in a small tubular  neighbourhood $U\supset Y$. 
\vspace {1mm} 

{\sf Step$_{\varepsilon}$ 1+2.} Take a suitable function $M_\varepsilon(x)$ in $U$   for which $Y$ strictly minimizes the functional  
$$Y\mapsto vol_{n-1}(Y)-\int_ {U_+}M_\varepsilon(x)dx,$$ 
where $U_+\subset U$ is the part of $U$ positioned {\it inward of} $Y$,
 such that the corresponding  minimizer $Y_\varepsilon$ of this functional for $g_\varepsilon$.
 (Such an $M$ on $Y$ must be equal to $mn.curv(Y)$)

\vspace {1mm} 

\hspace {3mm} {\sc can't fully   escape from the $\varepsilon$-neighbourhood $U_\varepsilon \subset U$ of $Y$.}
\vspace {1mm}

{\sf Step$_{\varepsilon}$ 3.} Observe that the filling inequalities $[\odot^{n-1}_{n-2}]$ and $[\odot^{n}_{n-1}]$
are stable under small perturbations $g_\varepsilon$ of $ g$ and conclude that \vspace {1mm} 

\hspace {31mm} {\sc $Y_\varepsilon$ is trapped 
in $U_\varepsilon$. }\vspace {1mm} 

(This $\varepsilon$ may be  slightly larger than the above  one.)
\vspace {1mm}

Then  the smoothness of $Y_\varepsilon$ ($n\geq 8$ included) follows from Almgren's optimal  isoperimetric inequality  (see  [Alm 1986]), which also allows a construction of diffeomorphisms $\psi_{\varepsilon} :X\to X$ which send $Y\to Y_\varepsilon$ (see [Gr 2012] for details).

\vspace {1mm}

{\it Warning.}  The hypersurfaces $Y_\varepsilon$ do not, in general, $C^1$-converge to $Y$ and,  conceivably, there are examples (I have not scrutinised the literature),  where there are  {\it no  diffeomorphisms} $\psi_\varepsilon $ having  the norms of their differentials (and/or of their inverses)  {\it   bounded} by $1+\varepsilon$.

On the other hand, one can control some  {\it H\"older norm} of $\psi_\varepsilon $ according to    {\it Reifenberg's topological disk theorem}.

Also {\it Reifenberg's flatness condition} implies   relative versions of the filling inequalities $[\odot^{n-1}_{n-2}]$ and $[\odot^{n}_{n-1}]$ from the previous section that is useful  for smoothing "intrinsic edges and corners" in manifolds with $Sc\geq \sigma$. (see  the next section and 
 [Gr 2014*]).

\vspace {1mm}

\hspace {41mm}\textbf {Localization of  [{\small   $\Circle_{\varepsilon} $}].}
\vspace {1mm}

 [{\small   $\Circle^{loc}_{\varepsilon} $}].} {\sl  If the metrics $g_\varepsilon$ are equal to $g$ on a neighbourhood $U_0\subset X$ of compact subset $X_0\subset X$ then the above 
diffeomorphisms $\psi_\varepsilon$ can be taken equal  to the  identity map on another (smaller) neighbourhood $U_\varepsilon\supset X_0$.}

\vspace {1mm}

{\it About the Proof.} This is achieved with suitable functions $M_\varepsilon$ defined on the complement of  
$X\setminus U_\varepsilon$. And here, as on other localization occasions, we leave  the actual proof to the reader, since the corresponding  localised properties are only  marginally used in the present paper.

\vspace {1mm}

\hspace {31mm}\textbf {Stability Relative to  $\partial  X$.}
\vspace {1mm}

In this paper, we need  a   relative version of [{\small   $\Circle_{\varepsilon} $}], and also of   [{\small   $\Circle^{loc}_{\varepsilon} $}],  where $X$ is a manifold {\it with boundary} and our hypersurface $Y\subset X$ has $\partial Y\subset \partial X$. 

The above  argument applies (almost)   word for word   to such $Y$, where,  additionally, one has to keep track of the dihedral angles between (the tangent spaces of) $Y$ and $\partial X$ along $\partial Y\subset \partial X$.

For instance, \vspace {1mm}

{\sl if these  angles  on  the  inward side of $Y$ satisfy $\angle_{in} (Y, \partial X) <\frac {\pi}{2}$, then   also 

$\angle_{in} (Y_\varepsilon , \partial X) <\frac {\pi}{2}.$}

\vspace {1mm}

This together with Bending Lemma (section 11.5) combined with $\varepsilon$-redistribution of curvature (section 11.2)  yield  the following.\vspace {1mm}

{\sc {\large $\Leftcircle\hspace {-1mm}_\varepsilon$}-Flattening Corollary.}
 Let $X$ be a   Riemannian manifold with smooth boundary 
and let $Y\subset X$ me a compact  smooth cooriented  hypersurface, where $\partial Y\subset \partial X$ 
and such that   $mn.curv (Y)>0$ and the inward (with respect to the coorientation of  $Y$) dihedral angles between $Y$ and $\partial X$ are everywhere $\leq \frac {\pi}{2}$.

{\sl Then there exists a family   of smooth Riemannian  metrics $g_\varepsilon$, $\varepsilon>0$, on $X$, such that 

$\bullet_1$ the hypersurface $Y$ is totally geodesic in $X$ with respect to $g_\varepsilon$ for all $\varepsilon>0$;

$\bullet_2$ $mn.curv_{g_\varepsilon}(\partial X)\geq mn.curv_g(\partial X)$ and the inward  dihedral angles with respect to $g_\varepsilon$ between $Y$ and $\partial X$  are everywhere $\leq \frac {\pi}{2}$  for all $\varepsilon>0$;

$\bullet_3$ the scalar curvatures of $ g_\varepsilon$ are bounded from below at all points $x\in X$ by $Sc(g)(x)$;

$\bullet_4$ the metrics $g_\varepsilon$ are equal to $g$ outside the $\varepsilon$-neighbourhood (with respect to $g$) of the union $Y\cup \partial X$;

$\bullet_5$ the  quadratic  forms  $g_\varepsilon-g$ are positive semidefinite   for all $\varepsilon>0$. }\vspace {3mm}

\hspace {30mm}{\Large $\bullet_?$  On Convergence  $g_\varepsilon\to g$.}

Ideally, one would like to have $C^0$-convergence $g_\varepsilon\to g$
 for $ \varepsilon\to 0$  but all we are able to show is that  only the {\it negative part} of the difference 
  $g_\varepsilon-g$ tends to zero, which, by  small perturbations of $g_\varepsilon$,  allows one to  make $g_\varepsilon-g$  positive semidefinite.
  
To see where the problem resides, let us return  to a closed smooth cooriented hypersurface $Y$   in a Riemannian manifold $X=(X,g)$  at the beginning of this section, 
and a family of smooth Riemannian metrics  let $g_\varepsilon $ on $X$ which $C^0$-converge  to $g$ for $ \varepsilon\to 0$. 

 Assume, for simplicity's sake, that $mn.curv(Y)=0$  and let $Y_\varepsilon$ be the  perturbations of  $Y$ which have  $mn.curv(Y_\varepsilon)\to 0$ and which themselves uniformly converge to $Y$.
 
 We do know that these $Y_\varepsilon$ are diffeomorphic to $Y$   but we do not expect their Riemannin metrics $h_\varepsilon=g_{| Y_\varepsilon}$  induced from $g_\varepsilon $ to converge to $h= g_{|Y}$.  All we can say is that  the $g$-normal projections $Y_\varepsilon \to Y$, besides being homotopic to diffeomorphisms, are 
  $(1+o(1))$-Lipschitz. 
  This eventually transforms   to the above $\bullet_5$. 
 
 And albeit our construction of $g_\varepsilon$ can't deliver the $C^0$-convergence $g_\varepsilon \to g$,  probably,  one can show that     the distance functions $dist_{g_\varepsilon} $ on $X\times X$  uniformly converge to $dist_g$ for $\varepsilon\to 0$.

\vspace {1mm}

\hspace {10mm}{\sc On stability of piecewise smooth hypersurfaces.}\vspace {1mm}

Here, as everywhere in this paper, "piecewise smooth" hypersurfaces $Y$  in $X$ are, by definition,  locally  diffeomorphic  to  polyhedral hypersurfaces in $\mathbb R^n$, $n=dim(X)$.

If $X$ comes with a Riemannian metric $g$, these  $Y$,  if they are cooriented,  are characterised by the mean curvatures of their ${n-1}$-faces  and the dihedral angles between these faces. 

Apparently, as indicated  (without proof)  in section  4.8 in [Gr 2014*] these $Y$ satisfy a piecewise-smooth version of the above [{\small   $\Circle_{\varepsilon} $}].

 This  means  that, \vspace {1mm}

{\sl given an $\varepsilon$-family of smooth Riemannian metrics on $X$, which $C^0$-converge to a  metric $g$  as in}  [{\small   $\Circle_{\varepsilon} $}],{\sl
 $$g_\varepsilon\underset {\varepsilon\to 0} \to g, $$} 
 {\sl there exists a family of piecewise smooth homeomorphisms   $\psi_\varepsilon :X \to X$, $\varepsilon>0$,  such that 
 
$\bullet$ $\psi_\varepsilon$ are smooth on the $(n-1)$-faces of $Y$     away from the    $(n-3)$-faces,

$\bullet$  $\psi_\varepsilon$  converge to the identity map $id: X\to X$ in the $C^0$-topology;

$\bullet$ the $(n-1)$-volumes of the hypersurfaces $Y_\varepsilon=\psi_\varepsilon(Y)\subset X$ with respect to $g_\varepsilon$ converge to 
$vol_{n-1}(Y)$ for $g$;

$\bullet$ the $g_\varepsilon$-mean curvatures of the faces of  $Y_\varepsilon$ converge to $mn.curv_g$ of the corresponding  faces of $Y$;

$\bullet$ the $g_\varepsilon$-dihedral angles between the faces of $Y_\varepsilon$ converge to the  $g$-dihedral angles between the corresponding faces of $Y_\varepsilon$.}

\vspace {1mm}

The relative version of this   in manifolds $X$  with corners (i.e. with piecewise smooth  boundaries),  which  also seems to follow by the available techniques\footnote{ The local version of this must be also true}, 
 would automatically yield  smoothing of  metrics on manifolds obtained by reflections of manifolds with corners with no decrease of their scalar curvatures (see section 4.8 in [Gr 2014*]). 

On the other  hand, one can prove the existence of such  smoothings  in the essential cases by some  roundabout argument as we shall explain  in the next section. \footnote{An unpleasant technical difficulty in the proof of the mean curvature stability in the piecewise smooth case, which  one has to {?} go around,  is the absence {?} of $C^1$ regularity theorem for minimal  $Y\subset X$  at the singular boundary points of $X$.} 

\subsection {Flattening of Faces and Regularisation of Reflections.}

\textbf {$\square$-Flattening Lemma.}   Let $X = (X,g)$ be a Riemannian manifold with corners, where all faces of the boundary $\partial X$ have positive mean curvatures and where
the dihedral angles between the  pairs of $(n-1)$-faces in $\partial X$, wherever they meet,  are $\leq \frac {\pi}{2}$.\vspace {1mm}

Then there exists  a family of smooth  metrics $g_\delta$, $\delta>0$, such that   \vspace {1mm}

(1) {\sl  the faces of $\partial X$ are totally geodesic with respect to $g_\delta$  for all $\delta>0$};

(2) {\sl all dihedral angles in $ \partial X$ with respect to $g_\delta$ are   $\frac {\pi}{2}$ for all $\delta>0$};

(3) {\sl the scalar curvatures of $g_\delta$ are bounded from below at all $x\in X$ by $Sc{(g)}(x)$};

(4) {\sl the metrics $g_\delta$ coincide with $g$ outside the $\delta$-neighbourhood of the boundary $\partial X$};


(5)   {\sl  the differences  $g_\delta-g$ are positive semidefinite for all $\delta>0$.}

\vspace {1mm}

This is a special case of the  Approximation/Reflection Lemma in section 4.9 in [Gr 2014*] where this  is proven not only  for $g$ itself, but for metrics $g_\varepsilon$ which $C^0$-converge to $g$.\footnote{ This lemma is formulated in [Gr 2014*] only for $Sc(g)>0$ and without formulating the above   (4) and (5). However the proof of  this lemma indicated in  [Gr 2014*]  automatically  deliver these properties .}

\vspace {1mm}

Granted that, one can reflect $X$ around the faces (see \textbf {2} in  section 11.1) smoothly the   metric in the resulting manifold $\tilde X$  by using   Cut-off Homotopy   Lemma  ($\bigstar$ in 11.1) and applying a corresponding inequality for manifolds (bands) without  corners.\vspace {1mm}

Thus, for instance,  one shows in [Gr 2014*] that \vspace {1mm}

{\sl Riemannian manifolds $X$ with $Sc(X) > 0$ can't contain mean curvature convex (e.g. convex) cubical domains $Q$, where all dihedral angles,  are non-obtuse.}\vspace {1mm}

Below is another { \sc Example}. \vspace {1mm}

{\sf \large Sub-Rectangular  $\frac {2\pi}{n}$-Inequality}. Let $X$ be a Riemannian $n$-manifold, let 
$Q\subset X$  be a   domain diffeomorphic to   the  $n$-cube $[-1,+1]^n$  and let $Q_i^\pm\subset \partial Q\subset Q$, $ i=1,...,n$,   denote the pairs of opposite codimension 1 faces in $Q$ which correspond to 
such pairs  in the cube.

{\sl Let

(i) the  faces $Q^{\pm}_i$  for  $i= 1,...,n-1$, are {\it mean curvature convex}, i.e.  
$$mn.curv(Q^{\pm}_i)\geq 0,$$

(ii)  the   dihedral angles   $ \angle_{\pm i, \pm j}  =\angle (Q^{\pm}_i, Q^{\pm}_j)$ between these faces are {\it non-obtuse}  at all points in the $(n-2)$-"edges"  where these faces meet,
$$ \angle_{\pm i, \pm j}\leq \pi/2, \mbox {  for all } i,j=1,...,n-1,  i \neq j,$$

(iii) the scalar curvature of $X$ satisfies $Sc(X)> n(n-1)$.\vspace {1mm}

Then the distance between the two remaining opposite faces satisfy 
$$dist_\pm=dist_X(Q^+_n, Q^-_n)< \frac  {2\pi}{n}.\leqno{\mbox {$\left[\square_\pm<\frac {2\pi}{n}\right ]$}}$$}

{\it Proof.} 
 Reflect $Q$, but now only  only in the faces $Q^\pm_i$ with $i<n$. Thus we construct  a torical band with $Sc> \sigma$ and apply  $[\circledcirc_\frac {2\pi}{\sqrt\sigma}]$ to this band.

\vspace {1mm} 

\hspace {26mm}{\sc  Flattening and  Gluing with $Sc>0$.} \vspace {1mm}

The proof of the  Approximation/Reflection Lemma  in [Gr 2014*] proceeds by    induction on the number of faces  with  application  of  Reifenberg's flatness property of minimal varieties at each  step of induction. 

But if the boundary of $X$ consists of  {\it only  two } (possibly  disconnected) $(n-1)$-faces then  the $\square$-Flattening Lemma  follows from the  $\Leftcircle\hspace {-1mm}_\varepsilon$-Flattening Corollary  from the previous section.

Indeed if $Y\subset $X is a {\it  totally geodesic} hypersurface {\it normal}  to $\partial X$, then the intrinsic bending delivered by the proof of the Bending Lemma from section 11.5 does not disturb $Y$ and  we have both $Y$ and $\partial X$ totally geodesic as well as mutually normal.
 
This also applies if     the set of $(n-1)$-faces one can be   divided into  {\it two} subsets of mutually disjoint ones, i.e. if 
the incidence  graph  
$\cal I$ between the faces  is 
{\it bipartite}.

 And  assuming no three $(n-1)$-faces   meet,  $\cal I$    can be artificially  made bipartite   by  subdividing the faces. 
 
 This is done by creating   "new narrow"  $(n-1)$-faces positioned close to  the  $(n-2)$-faces on one side of them. ( The $(n-2)$ faces  are  assumed  {\it two-sided} , i.e.  coorientable in $\partial X$.)
 
 Observe that the  the dihedral angles  between the  "new" $(n-1)$-faces and the "old"  ones are equal to  $\pi$, which, however,   leads to no problem due to the possible    localization of bending (section 11.5 ) and of [{\small   $\Circle_{\varepsilon} $}]  (section 11.7).

\vspace {1mm}

Now let us explain how one can get rid of "higher  order corners"   in  $X$  by paying the  price of a  change of their   topologies. In fact this "price"  limits the application of such gluing to $dim(X)=3$

Given an $n$-dimensional   cosimplicial manifold $X$ with corners, e.g. diffeomorphic to the $n$-cube $[0,1]^n$, one may {\it double} it over the set of its (preliminarily cut off) {\it highest  order corners}, where the {\it maximal numbers}, say  $m$,  of  the $(n-1)$-faces in $\partial X$ meet.

For instance, if $X= [0,1]^n$  then $m=n$  where these highest corners are the ordinary  vertices in  $[0,1]^n$.

The resulting double, call it $X^{[1]}$ has a natural corner structure where the highest corners have order $m^{[1]}=m-1$.
Thus we can continue unless we arrive at $X^{[m]}$  with smooth boundary (see  section 1.1 in [Gr 2014*]).

Geometrically, if $X$ is the  ordinary cube $[0,1]^n\subset \mathbb R^n$, (almost but not quite  exactly)  this can be achieved by cutting the  faces of this cube  of dimensions $\leq n-2 $ by hyperplanes and then by reflecting the resulting polyhedron around its (new as well as old) $(n-1)$-faces. 

Now we want to stop at $X^{[m]}$, for $m\leq n-2$,   such that the only singularities are $(n-2)$-dimensional  faces/"edges", where pairs of 
$(n-1)$-faces meet.
\vspace {1mm}

{\sc Gluing/Surgery Lemma.} Let $g$ be  a Riemannian metric on  $X$, for which all faces have $mn.curv_g>0$ and all dihedral angles are $ \leq  \alpha \leq \pi$.
\vspace {1mm}

{\sl Then  $X^{[m\leq n-2]}$ admits a smooth  metric $g^{[m]}$,  where the  faces also have $mn.curv_{g^{[m]}}>0$ and and all dihedral angles are   $\leq \alpha $, and such that the scalar curvature of $g^{[m]}$ is bounded from below (in a natural sense) by that of $g$.

For instance, if $Sc(g)>0$, then also $Sc(g^{[m]})>0$.}
 \vspace {1mm}
 
 {\it Sketch of the Proof.} Let first $X$ be diffeomorphic to the $3$-cube $[0,1]^3$. By consecutively  applying the  Cut-off Homotopy   Lemma at a  vertex $p\in [0,1]^3$  to the three $2$-faces at $p$, one can "infinitesimally straighten" these faces at $p$, i.e.  make them  geodesic at $p$, while   keeping $mn.curv\geq 0$. 
 
If  there are two 3-manifolds,  $X$ and $X'$ with such vertices, one can arrange the connected sum  
$$(X,p) \#  (X,p')$$
by "gluing" them  with an arbitrarily small decrease   of their scalar curvatures by the same construction as it is done for the ordinary connected sum in [GL 1980*], where one can preliminarily enlarge the scalar curvature of the two at the points $p$ and $p'$ by means of the  Cut-off Homotopy   Lemma  (as it is done on  p. 111 in [Gr 1986] in a similar context.)

Then  the required double of our  cubical  $X$ is achieved by "infinitesimally straightening" the faces at all  eight vertices in $X$  and  then taking the double  with a  "glueing metric"  at  all vertices.

The above  applies to 0-faces (vertices) $p$ of all  $X$  for all $n=dim(X)$.
In general, if   $P$ is an $m$-face for  $0<m\leq n-2$,   one  "infinitesimally straighten" the $(n-1)$-faces {\it normally to} $P$, where a face $P^+\supset P$  is regarded infinitesimally straight normally to $P$ if the second fundamental form of $P\subset X$ {\it vanishes} on the tangent subspaces $T_x\subset T_x(P^+)$  
{\it normal to $P\subset P^+$} at all points $x\in P$, where, observe, $rank (T_x)=  dim (P^+)-dim(P)\geq 2$.

\vspace {1mm}
  
{\it Clarifying Topological Remark.} To see what makes the difference between dimensions $n=3$ 
and $n\geq 4$, let $V$ be  an $n$-dimensional manifold $X$   minus an open  tubular neighbourhood of a closed  submanifold $Y$. 

If  $X$ carries a metric with $Sc>0$ and  $codim(Y)\geq 3$, then the double $W=V\cup_{\partial V}V$ also carries such metric. 

Now,  if $X$ is overtorical and the submanifold $Y$ can be  homotoped to a single point in $X$, e.g. $Y$ is a finite union of points in $X$,  then $W$ is also overtorical.

But if, for instance, $X=\mathbb T^n$ and $Y$ is the union of $n$ disjoint circles  which generate 
$H_1(X)$, then $W$ is {\it not} overtorical. In fact, the corresponding $W$ in this case does carry a metric with positive scalar curvature.

\subsection {Non-existence Results and Conjectures.}

Following a suggestion by the referee, we briefly  overview in this section a few   constrains,  some  of which  proved in the main body of this article and some  conjectural,   on  the topology of manifolds which carry metrics with $Sc>0$.\footnote{Unsolved problems on $Sc>0$ are collected in [Gr 2017] and [Gr 2017*].}
\vspace {1mm}

 Say that a closed oriented  manifold $X$ is SYS (Schoen-Yau-Schick) {\it over a cohomology class $h$ in a  topological space} $K$, where 
$$h\in H^{n-2} (K;\mathbb Z), \hspace {1mm} n=dim(X),$$
if there is a continuous map $A: X\to K$, such that the $2$-dimensional  homology class in $X$ which is  the  Poincar\'e dual of the cohomology pullback $A^\ast (h)\in   H^{n-2}(X;\mathbb Z)$    of $h$,
$$ A^\perp(h)=PD(A^\ast (h))\in H_2(X),$$
 is {\it non-spherical},  i.e.  it  is {\it not contained} in the image of the {\it Hurewicz homomorphism} $\pi_2(X)\to H_2(X)$.
 \vspace {1mm}

For example, "SYS  over the fundamental cohomology class of the  torus $\mathbb T^{n-2}$"  is  the same as "SYS" defined in section 5.

 \vspace {1mm}

Let us generalise the Schoen-Yau theorem ([SY - Str 1979] for $n\leq 7$ and 
[SY 2017] for all $n$) on {\sl non-existence of  metrics with positive scalar curvatures
on SYS-manifolds}  and   (some case of)  theorem 13.8 in   [GL 1983] as follows. 
  \vspace {1mm}

{\textbf{$\cal SYSE$-Non-existence Theorem.} $K$, be  a  manifold which  admits {\it a  complete metric with non-positive sectional curvature} and let   $X$ be  a closed manifold  $X$  which  is {\it SYS over an  integer  cohomology class $h\in H^{n-2}(K) $  which does not vanish in $h\in H^{n-2}(K;
\mathbb Q) $}.    \vspace {1mm} 

{\sl Then $X$ admits no metric with positive scalar curvature.}
  \vspace {1mm}

{\it Proof.}    Start with the  case where $K$ is compact oriented of dimension  $n-2=dim(X)-2$ and 
$h\in H^{n-2}(K;\mathbb Z)$ is the fundamental class of $K$.  Let $\tilde K$ be the universal covering of $K$ and $\tilde X$ be  the covering of $X$ induced by the above map $A: X\to K$. 

Since $\tilde K$ has non-positive curvature, it admits proper $\lambda$-Lipschitz maps $F_\lambda$ to the unit ball
$B=B^{n-2}(1)\subset 
\mathbb R^{n-2}$ of degree  1 for all $\lambda>0$ and  then, when  $\lambda\to 0$, the proof follows from 
the stable  SYSE-bound on the scalar curvature from section 6  applied   to the composed maps 
$$ \tilde X\overset {\tilde A}\to \tilde K\to B$$
in the role of $f=f_\lambda$ and    $X=X\times \{p\}$, where $\{p\}$  a single point space for $P$.

Now let us turn to the case of  a complete orientable  $K$ of dimension $m\geq  n-2$ and let $h\in H^{n-2}(K;\mathbb Z)$. Assume $K$ 
is parallelizable, otherwise, pass to the total space of the normal (for an embedding $K\to \mathbb R^M$)  vector bundle $T^\perp(K)\to K$, which (by an easy argument) also carries a metric with non-positive sectional curvatures. 

In this case,  there exists  continuous proper maps  $F_\lambda: \tilde K\times K\to B=B^m(1)$, $\varepsilon>0$,  which are
 $\lambda$-Lipschitz  (diffeomorphisms)   on all "slices "  $\tilde K\times k$, $k\in K$: these maps are constructed  as above with the use of inverse exponential maps in $\tilde K$ at the points over $k\in K$, where the  tangent spaces $T_k(K)$ are identified with 
a single $\mathbb R^m$ with a use of  a frame in $T(K)$ as it is done in section 13 in [GL 1983]).

Let $P\subset K$ be a subpseudomanifold of codimension $n-2$, which represents the Poincare dual of $h$  (which is a homology class  with, a priori,  infinite support)  and let 
$f_\varepsilon: \tilde X\times P\to B$ be the maps obtained by composing $\tilde A:\tilde X\to \tilde K$ and the restriction of $F_\lambda$ to $\tilde X \times P,$ 
$$f= f_\lambda(\tilde x, p)=F_\lambda(\tilde A(\tilde x),p).$$

Clearly, this map $f:X\times P\to B$, satisfies the assumptions on $f$ in section 11.6, and 
  stable  SYSE-bound on the scalar curvature from 11.6  applies. QED.

  \vspace {1mm}

 {\sc Conjecture \textbf {A}.}  {\sf If a closed manifold $X$ is SYS over  a non-torsion  cohomology class  $h$ in an aspherical space $K$, then  $X$  admits no metrics with $Sc>0$. }

 \vspace {1mm}
 
 This conjecture may be  unrealistically strong but the special case of this, where
 
  \hspace {-3mm}$K$ {\sl admits a complete  (possibly singular) metric  with non-positive curvature} 
 
   \hspace {-6mm}seems within reach.
 
 \vspace {2mm}
 
 Recall (see section 10) that a closed oriented manifold $X$ is called 
 
 \hspace {-6mm}[{$\tilde \uparrow 0$}]-oversymplectic  if

 $\bullet$  a multiple  of the  fundamental cohomology class  of $X$ decomposes into product of one and two dimensional classes,
$$k\cdot [X]^\circ=h_1\smile...\smile h_m,  $$
and 

 $\bullet $ $X$ 
 the  classes $h_i$   vanish in the cohomology of the  universal covering $\tilde X$. 
 
Also recall (see section 10) that 
 
 [{$\tilde \uparrow 0$}]-{\it oversymplectic manifolds, the universal covers of which are spin, carry no metrics with $Sc >0$.} 

Probably, the spin condition is redundant; moreover, one may merge this with the above \textbf {A} as follows.

\vspace {1mm}

{\sc Conjecture \textbf {B}.} {\sf Products of manifolds $X$ as in the above \textbf {A} by [{$\tilde \uparrow 0$}]-oversymplectic
ones admit no metrics with $Sc>0$.}

\subsection {A Few Geometric Problems and Conjectures.}

{\sc Conjecture \textbf {C}.} {\sf If a closed manifold $V_0$  of dimension $n-1\geq 5$,  admits no metric with $Sc>0$  then Riemannian bands $V$  diffeomorphic to $V_0\times [-1,1]$ which have $Sc(V)\geq \sigma>0$, satisfy the sharp width inequality,
$$width(V)\leq 2\pi\sqrt{\frac { n-1}{\sigma n}}.$$}

{\sc Conjecture \textbf {D}.} Let $ g_0$ stands for  the standard Riemannian metric on the unit sphere  $S^n$ with the sectional curvature 1. 

{\sf If a  RIemannian metric $g$ on $S^n$ minus a point satisfies
$$ g\geq g_0\mbox {  and } Sc(g)\geq Sc(g_0)=n(n-1),$$ 
then $g_0=g$.} 

  (If $g$ is complete, this  follows,  by  the relative  index theorem for the Dirac operator, see  [Ll 1998]  and (8) in section 10)
\vspace {1mm}

{\sc Conjecture \textbf {D'}.} Let $X$ be closed $n$-manifold, such that $X$ minus a point admits {\it no  complete metric with} $Sc>0$.

Let $V$ be obtained by removing a  small open $n$-ball from $X$, i.e. $V=X\setminus B_{x_0}(\varepsilon)$, and let $g$ be a metric on $V$
with $Sc(g)\geq \sigma>0$. 
 If the $\rho$-neighbourhood  with respect to $g$  of the boundary sphere $S^{n-1}= \partial V=\partial B_{x_0}(\varepsilon)$  is homeomorphic to $S^{n-1}\times [0,1]$, then  
$$\rho\leq \frac {20}{\sqrt \sigma}.$$

(If $X$ is a $SYS$-manifold, then  metrics $g$  with $Sc(g)\geq \sigma$ on  $V$  do satisfy this inequality as it follows by     Schoen-Yau's kind of   argument adapted to manifolds with boundaries as in  section 11.6. 

On the other hand,  the conjecture must be  vacuous for {\it simply connected} manifolds $X$, since $X\setminus \{x_0\}$ for such an $X$  contracts to an $(n-2)$-subpolyhedron in $X\setminus \{x_0\}$, which,  most probably, implies that $X\setminus \{x_0\}$  admits a complete  metric with $Sc>0$.) \vspace {1mm}

{\sc Conjecture \textbf {E}.} {\sf Let $V$ be a Riemannian manifold homeomorphic to $\mathbb T^2\times [-1.1]$ with sectional curvature everywhere  $\geq 1$. Then} (this was already mentioned in section 3)
$$width (V) =dist(\partial_-(V),\partial_+(V))\leq \pi/2,$$ 
(The  simplest unsettled  case is where $V$ is a domain in   $S^3$ or, more generally if it admits an isometric embedding or immersion into $S^3$.
\vspace {1mm} 

Below  is a far reaching generalisation of the  spherical case of   { \textbf {E}.}\vspace {1mm}

{\sc Conjecture \textbf {E$_{+\infty}$}.} {\sf  Let $Y\subset S^N$, $N=n,n+1, ...   ,\infty$,   be a submanifold homeomorphic to  the product of $n$ closed manifolds of dimensions $\geq 1$, e.g.  homeomorphic to the $N$-torus
 and let $U\supset Y$  be a neighbourhood of $Y$ in $S^N$ which admits a retraction  to $Y$.
Then 
$$ dist(Y, \partial U)\le \arcsin \frac {1}{\sqrt n}.$$}

But, in reality,  one has no estimate for this distance even  for  high codimensional tori in  spheres, which suggests the following conjecture opposite to  {\sc \textbf {E$_\infty$}.}

{\sc Conjecture \textbf {$-$E}.} {\sf Every compact smooth manifold $Y^n$ of dimension $n$ admits a smooth embedding to the sphere  $S^{2n}$ such that all principal curvatures of the  image satisfy
$$ curv (Y^n \subset  S^{2n})\leq const< \infty,$$
say for $const=1000$.}\vspace {1mm}

It is hard to believe in the validity of either {\sc  \textbf {E$_{+\infty}$}}  or {\sc  \textbf {-E}},
but something in between may be true, e.g. the following. \vspace {1mm}

{\sc Conjecture \textbf {E$_{\mathbb T^n}$}.}
{\sf The minimal constant $\beta$, such that the $n$-torus admits an smooth  immersion to $S^{n+1}$
with principal curvatures $\leq \beta$ 
is  asymptotic,   for   $n\to \infty$,  to 
$$ const \cdot n^\beta\mbox{   for some $\beta > 1$ }. $$}

{\sc Problem  \textbf {F}.} Let $U\subset \mathbb R^N$ and $V\subset \mathbb R^n$, $n\leq N$, be open subsets, e.g. balls $B^N(r)$ and $B^n(R)$, $R\geq r$.  Evaluate the minimal $\beta =\beta (U,V)>0$, such that  $V$ admits a smooth {\it locally expanding  immersion/embedding}\footnote{A 
 a smooth map $f:V\to U$  is {\sl locally expanding} if the differential  $Df:T(V)\to T(U)$ 
 doesn't decrease the norms of the tangent vectors, $||Df(\tau)||\geq ||\tau ||$ for all $\tau \in T(V)$.}   
 to $U$ with the  principal  curvatures $\leq \beta$.

\vspace {1mm}

{{\it In Conclusion.} The  above  {\sc  \textbf {A}} - {\sc  \textbf {F}} are   only  tips of the  iceberg of what we don't know about the scalar curvature and nearabouts.

(An outline of this "iceberg" is given  in   [Gr 2017].)

\section {References.}

\hspace {4mm} [Al 1985]   Sebastiao Almeida, Minimal Hypersurfaces of a Positive Scalar 

Curvature.
Math. Z. 190, 73-82 (1985).

\vspace {1mm}

[Alm 1986] F.J. Almgren Jr.
Optimal isoperimetric inequalities, Indiana Univ. Math. J.
35 (1986), 451-547.\vspace {1mm}

 [AndMinGal 2007]  Lars Andersson, Mingliang Cai, and Gregory J. Galloway, Rigidity and positivity of mass for asymptotically hyperbolic manifolds, Ann. Henri Poincar\'e 9 (2008), no. 1, 1-33. 
\vspace {1mm}

 [At 1976] M. F. Atiyah. Elliptic operators, discrete groups and von Neumann algebras.
Ast\'erisque 32-3 (1976), 43-72.\vspace {1mm}

[Bam 2016]  R. Bamler . A Ricci flow proof of a result by Gromov on lower bounds for scalar curvature. 
Mathematical Research Letters
Volume 23 (2016). Number 2, Pages 325 - 337.\vspace {1mm}

[BER 2017]      Boris Botvinnik, Johannes Ebert, Oscar Randal-Williams, Infinite loop spaces and positive scalar curvature, Inventiones mathematicae
Volume 209, Issue 3, 2017, pp 749-835.  \vspace {1mm}

[BH 2009] M. Brunnbauer, B. Hanke, Large and small group homology , J. Topology 3 (2010) 463-486.\vspace {1mm}

[BHMM 2015] Jean-Pierre Bourguignon, Oussama Hijazi, Jean-Louis Milhorat, Andrei Moroianu and Sergiu Moroianu A Spinorial Approach to Riemannian and Conformal Geometry,EMS Monographs in Mathematics 2015.\vspace {1mm}

 [BMN  2010]  S Brendle, F.C. Marques A. Neves .Deformations of the hemisphere that increase scalar curvature,
 	arXiv:1004.3088 [math.DG]

 [DFW 2003] A. N. Dranishnikov, Steven C. Ferry, and Shmuel Weinberger,   Large Riemannian manifolds which are flexible,  Annals of Mathematics
 157 (2003), 919-938. \vspace {1mm}

[Dr 2000] A. Dranishnikov, Asymptotic topology
, Russian Math. Surveys 55:6 2000), 71-116.\vspace {1mm}

[Dr 2006 ] A. N. Dranishnikov,
On hypereuclidean manifolds  Geom. Dedicata 117 (2006), 215-231.\vspace {1mm}

[FCS 1980] D. Fischer-Colbrie, R. Schoen, The structure of complete stable minimal surfaces in 3-manifolds of nonnegative scalar curvature, Comm. Pure
Appl. Math. 33 (1980) 199-211.\vspace {1mm}

 [GL 1980] M. Gromov, B. Lawson,  Spin and Scalar Curvature in the Presence of a Fundamental Group I
Annals of Mathematics, 111 (1980), 209-230. \vspace {1mm}

[GL 1980*] M. Gromov and H. B. Lawson,  The classification. of simply connected. manifolds of positive scalar curvature. 111 (1980), Annals of Mathematics, 11, pp  423-434 (1980).

[GL 1983] M. Gromov and H. B. Lawson,
Positive scalar curvature and the Dirac operator on complete Riemannian
manifolds, Inst. Hautes Etudes Sci. Publ. Math.58 (1983), 83-196.  \vspace {1mm}

 [Gr 1983] M. Gromov.   Filling Riemannian manifolds. J. Differential Geom. 18, no. 1, 1--147. (1983).

 [Gr 1986]  M. Gromov, Partial differential relations, Springer 1986  
 \vspace {1mm}
 
  [Gr 1996] M. Gromov. Positive curvature, macroscopic dimension, spectral gaps and higher signatures.
In Functional analysis on the eve of the 21st century, Vol. II (New Brunswick, NJ, 1993) ,volume 132 of
Progr. Math., pages 1-213, Birkh\"auser, 1996.
\vspace {1mm}

  [Gr 2012] M. Gromov. Hilbert volume in metric spaces. Part 1.
Cent. Eur. J. Math.
 10(2):371-400, 2012.

 [Gr 2014] Gromov M., Plateau-Stein manifolds, Cent. Eur. J. Math., 2014, 12(7), 923?951.

 [Gr 2014*] M. Gromov, Dirac and Plateau billiards in domains with corners, Central European Journal of Mathematics, Volume 12, Issue 8,   2014, pp 1109-1156.
\vspace {1mm}

[Gr 2017]   101 Questions, Problems and Conjectures
around Scalar Curvature.

\url{ http://www.ihes.fr/~gromov/PDF/101-problemsOct1-2017.pdf}\vspace {1mm}

[GS 2001]  S. Goette, U. Semmelmann, Spin$^c$ Structures and Scalar Curvature Estimates, Annals of Global Analysis and Geometry
 Volume 20, Issue 4, pp 301-324, 2001
\vspace {1mm}

[GS 2002]  S. Goette and U. Semmelmann, Scalar curvature estimates for compact symmetric spaces.
Differential Geom. Appl.  16(1):65-78, 2002.

\vspace {1mm}

[Han 2011] B. Hanke. Positive scalar curvature, K-area and essentialness, Global Differential Geometry pp 275-302,
 (2011)

  \vspace {1mm}

[HaSch 2006]  B. Hanke, T. Schick, Enlargeability and index theory,J. Differential Geom. 74 (2) (2006), 293-320.\vspace {1mm}

[HaPS 2015] Hanke, Bernhard; Pape, Daniel; Schick, Thomas Codimension two index obstructions to positive scalar curvature, Annales de l'institut Fourier, Volume 65 (2015) no. 6 , p. 2681-2710.\vspace {1mm}

[Hit 1074] N. Hitchin, Harmonic spinors, Advances in Math. 14 (1974), 1-55.

\vspace {1mm}

 [Ll 1998] M. Llarull, Sharp estimates and the Dirac operator, Mathematische AnnalenJanuary 1998, Volume 310, Issue 1, pp 55-71.\vspace {1mm}

[LM 1984] H.B., Jr Lawson;,M.-L. Michelsohn,
Approximation by positive mean curvature immersions: frizzing.
Inventiones mathematicae Volume: 77, page 421-426 (1984).\vspace {1mm}

[Loh 1999] J. Lohkamp, Scalar curvature and hammocks, Math. Ann. 313, 385-407, 1999.
\vspace {1mm}

[Loh 2006]  J. Lohkamp,  The Higher Dimensional Positive Mass Theorem I, arXiv math.DG/0608795.\vspace {1mm}

[Loh 2008]  J. Lohkamp, Inductive Analysis on Singular Minimal Hypersurfaces, arXiv:0808.2035.\vspace {1mm}

[Lok 2016] J. Lohkamp,  The Higher Dimensional Positive Mass Theorem II
arXiv:1612.07505   \vspace {1mm}

[McFSzk 2012]  Donovan McFeron, G\'abor Sz\'ekelyhidi. On the positive mass theorem for manifolds with corners,
Communications in Mathematical PhysicsVolume 313, Issue 2, 2012 pp 425-443

\vspace {1mm}

[Mia 2002]  P.  Miao.  Positive  mass  theorem  on  manifolds  admitting  corners  along  a
hypersurface. Adv. Theor. Math. Phys., 6(6):1163-1182, 2002.\vspace {1mm}

[Min-Oo 1988] M. Min-Oo. Scalar curvature rigidity of certain symmetric spaces. In
Geometry, topology, and dynamics (Montreal, PQ, 1995), volume 15 of
CRM Proc. Lecture Notes , pages 127- 136
\vspace {1mm}

[Min-Oo 1989]  M.  Min-Oo,  Scalar  curvature  rigidity  of  asymptotically  hyperbolic  spin
manifolds. Math. Ann. 285, 527- 539 (1989)\vspace {1mm}

[Min-Oo 2002]  M.  Min-Oo,   K-Area, mass and asymptotic geometry.

\url{http://ms.mcmaster.ca/minoo/mypapers/crm_es.pdf}\vspace {1mm}

[NS 1993)] N. Smale,
Generic regularity of homologically area minimizing hyper
surfaces in eight-dimensional mani-
folds, Comm. Anal. Geom. 1, no. 2 (1993), 217-228.\vspace {1mm}

[Ros 2007]  J. Rosenberg, Manifolds of positive scalar curvature:  a progress report
, in:  Surveys on Differential Geometry, vol. XI: Metric and Comparison Geometry, 
International Press 2007.\vspace {1mm}

  [Sch 1998] Thomas Schick, A counterexample to the (unstable) Gromov-Lawson-Rosenberg conjecture,
Topology 37 (1998), no. 6.  \vspace {1mm}

   [SY - Inc 1979] R.  Schoen  and  S.  T.  Yau,  Existence  of  incompressible  minimal  surfaces
and  the  topology  of  three  dimensional  manifolds  of  non-negative  scalar
curvature, Ann. of Math. 110 (1979), 127-142.\vspace {1mm}
  
 [SY - Str 1979]  R. Schoen and S. T. Yau, On the structure of manifolds with positive scalar
curvature, Manuscripta Math. 28 (1979), 159-183.
  
  \vspace {1mm}

   [SY 2017] Positive Scalar Curvature and Minimal Hypersurface Singularities
Richard Schoen, Shing-Tung Yau, arXiv:1704.05490 [math.DG]

[St 2001] S. Stolz.
Manifolds of a positive scalar curvature, in: T. Farrell et
al. (eds.),
Topology of high dimensional manifolds
ICTP Lect. Notes, vol. 9, 665-706. 1.1, Trieste (2001).




\end{document}